\documentclass[11pt]{article} 

\usepackage[colorlinks=true,linkcolor=blue,urlcolor=blue]{hyperref}
\usepackage{dsfont}

\usepackage{amsfonts,latexsym,amstext} 
\usepackage{amsmath} 
\usepackage{amssymb} 
\usepackage{comment} 
\usepackage{amsthm} 
\usepackage[latin1]{inputenc} 
\usepackage{color}

\DeclareMathOperator{\supp}{supp}                
 
\theoremstyle{plain} 
\newtheorem{theorem}{Theorem}[section] 
\newtheorem{lemma}[theorem]{Lemma} 
\newtheorem{proposition}[theorem]{Proposition} 
\newtheorem{corollary}[theorem]{Corollary} 
 \def\R{\rz R}

\def\N{\rz N}
\def\E{\rz E}

\def\P{\rz P}

\def\shf{{\cal F}}

\def\shl{{\cal L}}

\newtheorem{definition}[theorem]{Definition} 
\newtheorem{hypothesis}[theorem]{Hypothesis} 
\theoremstyle{remark} 
\newtheorem{remark}[theorem]{Remark} 
\newtheorem{example}[theorem]{Example}

\allowdisplaybreaks 
 
\makeatletter 
\@addtoreset{equation}{section} 
 
\makeatother 
 
\def\qed{{\hfill\hbox{\enspace${ \square}$}} \smallskip} 
\def\sqr#1#2{{\vcenter{\vbox{\hrule height .#2pt \hbox{\vrule 
 width .#2pt height#1pt \kern#1pt \vrule 
width .#2pt} \hrule height .#2pt}}}} 
\def\square{\mathchoice\sqr54\sqr54\sqr{4.1}3\sqr{3.5}3}

\def\ds{\begin{displaystyle}} 
\def\eds{\end{displaystyle}} 
 
\def\<{\langle } 
\def\>{\rangle }

\def\R{\mathbb R} 
\def\N{\mathbb N}

\def\E{\mathbb E} 
\def\P{\mathbb P} 
 
\def\F{\mathbb F}

\newcommand{\sper}[1]{\mathbb{E} \left[ #1 \right]}                               
 
\DeclareMathAlphabet{\mathonebb}{U}{bbold}{m}{n}                           %
\newcommand{\one}{\ensuremath{\mathonebb{1}}}                               

\title{Weak Dirichlet processes and generalized martingale problems}

\author{ 
Elena Bandini 
\thanks{Universit\`a di Bologna, Dipartimento di Matematica, Piazza di Porta S. Donato 5, 40126 Bologna, Italy; e-mail: elena.bandini7@unibo.it.} 
\and Francesco Russo 
\thanks{UMA, ENSTA Paris, Institut Polytechnique de Paris, 
828, boulevard des Mar\'echaux, F-91120 Palaiseau, France; e-mail: francesco.russo@ensta-paris.fr.} 
} 
\date{}
 
\begin{document} 
  
\allowdisplaybreaks 
\maketitle 
 
\begin{abstract} 
In this paper we explain how the notion of {\it weak Dirichlet process}
is the suitable generalization of the one of semimartingale with jumps.
For such a process we provide a  unique decomposition which is new also for semimartingales: 
in particular we introduce {\it  characteristics} for weak Dirichlet processes. We also introduce a weak concept (in law) of finite quadratic variation.   We investigate a set of new useful chain rules 
and we discuss a general framework of (possibly path-dependent with jumps) 
martingale problems with a set of examples 
of SDEs with jumps driven by a distributional drift.

\end{abstract} 

{\bf Key words:} 
Weak Dirichlet processes; c\`adl\`ag semimartingales; jump processes;
martingale problem; singular drift; random measure.
 \\
 {\small\textbf{MSC 2020:}  
60H10; 60G48; 60G57
}

\section{Introduction}
The central notion of this work is the one of 
{\it weak Dirichlet process} with jumps and the related
{\it martingale problem}. 
In this work we want in particular to convince the reader
that the concept of weak Dirichlet process plays  a
similar 
central role as the one of semimartingale.
An $(\shf_t)$-weak Dirichlet process $X$
is the
sum of an $(\shf_t)$-local martingale $M$ and an $(\shf_t)$-martingale orthogonal process $\Gamma$, generally 
fixed to vanish at  zero. When self-explanatory, the filtration will be omitted.  A {\it martingale orthogonal process} $A$
 has the property that $[A,N] = 0$ for every continuous
 martingale $N$. In particular a purely discontinuous martingale  is a  martingale orthogonal process.  $A$ substitutes the usual bounded variation
 component $V$ when $X$ is a semimartingale. As a matter of fact, any bounded variation process is $(\shf_t)$-martingale orthogonal, see Proposition 2.14 in \cite{BandiniRusso1}.

When $X$ is a continuous process, the notion of  weak Dirichlet process was introduced in \cite{er}
and largely investigated in \cite{gr}.
In \cite{BandiniRusso1}, we extended the concept of weak Dirichlet 
jump process with related calculus. In particular, generalizing the notion of special semimartingale, we introduced the one of special weak Dirichlet process  $X = M + \Gamma,$
where $M$ is a possibly (càdlàg) local martingale
and $\Gamma$ is a predictable; 
 in that case the decomposition $X = M + \Gamma$ is unique.
In fact this was earlier introduced by  \cite{cjms}
and partially studied, omitting the mention ''special''.

An important feature of the calculus beyond semimartingales
is the one of stability, which often constitutes a generalization of It\^o formula. If $X$ is a c\`adl\`ag semimartingale and $f\in C^2(\R)$, we recall that
$f(X)$ is again a semimartingale by a direct application of It\^o formula.
However, if $f$ is only of class $  C^1(\R)$, $f(X)$ is generally 
no more a semimartingale, nevertheless it
remains a  finite quadratic variation process,
see \cite{rv4}
when $X$ is a continuous process. For instance, if $X$ is a Brownian motion and $f$ is not the difference of convex functions, then it is well-known that $f(X)$ is not a semimartingale. On the other hand, if $f$ is bounded and of class $C^1$ and $X$ is a Poisson process, then $f(X)$ is a special semimartingale.

A c\`adl\`ag process $X$ will be called  \emph{finite quadratic variation process}  (see \cite{rv95, BandiniRusso1})
whenever the u.c.p. limit (which will be denoted by $[X,X]$)     of
$ [X,X]^{ucp}_{\varepsilon}$ exists, where
\begin{equation} 
\label{Appr_cov_ucpI}	 
[X,Y]^{ucp}_{\varepsilon}(t):= \,\int_{]0,\,t]}\,
\frac{(X((s+\varepsilon)\wedge t)-X(s))(Y((s+\varepsilon)\wedge t)-Y(s))}{\varepsilon}\,ds.
\end{equation} 
By Lemma 2.10 in \cite{BandiniRusso1}, we know that
\begin{equation} \label{QVC}
[X,X] = [X, X]^c  + \sum_{s \leq \cdot} |\Delta X_s|^2,
\end{equation}
where $[X,X]^c$ is the continuous component.
The covariation of two c\`adl\`ag 
processes $X$ and $Y$ was defined in \cite{rv95} as the u.c.p. limit
of $[X,Y]^{ucp}_{\varepsilon}$, whenever it exists.

The notion of quadratic variation for a non-semimartingale $X$
was introduced  in \cite{fo} by means of discretizations (instead
of regularizations, that we denote here by $[X]$). 
One also proved that  if $f \in C^1(\R)$ and $[X]$ exists,  then
$[f(X)]$ also exists.
A natural extension of the notion of the semimartingale  is the one of
$(\shf_t)$-Dirichlet process introduced in \cite{FolDir} still 
in the discretization approach, which is the (unique) sum
of an $(\shf_t)$-local martingale and a zero quadratic variation process
(vanishing at zero).
For us, an $(\shf_t)$-Dirichlet process will be the analogous
concept in the regularization approach.
Let $X$ be such a Dirichlet process, which is obviously in
particular an $(\shf_t)$-(even special) weak Dirichlet process. Let
$X = M + \Gamma$ its unique decomposition.
Observe that the notion of Dirichlet process  does not
naturally fit the
jump case; indeed $[\Gamma,\Gamma] = 0$ implies
that $\Gamma$ is continuous by \eqref{QVC}, therefore predictable.
On the other hand, since
$X$ is a finite quadratic variation process, if $f \in C^1(\R)$ then
$f(X)$ is also a finite quadratic variation process but it is not necessarily
a Dirichlet process, see 
\cite{BandiniRusso_DistrDrift}.

For applications (for instance to control problems and BSDEs theory, see e.g. \cite{gr1, BandiniRusso2, FuhrmanTessitore}), it is useful to investigate stability for functions 
 $f \in C^{0,1}([0,\,T] \times \R; \R)$. 
Let  $f \in C^{0,1}([0,\,T] \times \R; \R)$ and $X$ be a finite quadratic variation process. In general we cannot expect that $f(t,X_t)$ to be a finite quadratic variation process: for instance a very irregular function $f$ not depending on $x$ may not be of finite quadratic variation.

If $X$ is a continuous weak Dirichlet process with finite quadratic variation, it is known that $f(t, X_t)$ is a weak Dirichlet process, see
Proposition 3.10 in \cite{gr}.
This stability result was extended to the case where $X$ is a discontinuos
process, provided a specific relation between the jumps of $X$ and $f$, see \cite{BandiniRusso1}. Under these conditions, $f(t, X_t)$ is a special weak Dirichlet process. Recently, an interesting generalization in the continuous framework has been provided in \cite{BLT}, where $f$ is a $C^{0,1}$-path-dependent functional in the sense of horizontal-vertical Dupire derivative.
A former work including $C^{1,2}$-chain rules for a significant class of path-dependent  processes with jumps is \cite{cordoni}.
The family of weak Dirichlet processes has also interesting connections
with  the so called  stochastically controlled processes,
see \cite{GOR}, and also the related interesting recent reference \cite{Hocquet}.

As mentioned earlier, the second central notion of the present
paper is the one of  martingale problem, whose classical notion  is due to Stroock and Varadhan,
see e.g. \cite{Stroock_Varadhan}. In general, one says that a process
$X$ is a solution to the martingale problem with respect to a probability
$\P$ (on some probability space), to some domain
${\mathcal D} $ and to a time-indexed family of operators $L_t$ if
$$ f(X_t) - f(X_0) - \int_0^t L_sf(X_s) ds, $$
is a $\P$-local martingale.
One also says that $(X,\P)$ is a solution to the martingale problem related to
${\mathcal D}$ and $L_t$.
In the classical martingale problem in \cite{Stroock_Varadhan}
one takes ${\mathcal D} = C^2(\R^d)$ and $L_t$ is a
second order PDE operator. One also knows that
the solution of  Stroock-Varadhan martingale problem is equivalent
to the one of an SDE in law (or weak).
In more singular situations, the notion of SDE seems
difficult to exploit and define, and in that case it
is substituted by the more flexible notion of martingale problem.
This is the case for instance when one investigates
the notion of {\it SDEs with distributional drift}, i.e.
of the
type
$$ X_t = X_0 + \int_0^t \sigma(X_s) dW_s + \textup{``}\int_0^t b(X_s) ds\textup{''},$$
and $b$ is a Schwartz distribution. In this case it is more comfortable
to express those processes as solutions to a martingale problem
with respect to some domain $ {\mathcal D}   $ which is a suitable
subset of $C^1(\R)$, where
the formal map  $L_t f(x) := \frac{1}{2} \sigma^2(x) f''(x) + b(x) f'(x)$
is well-defined (independently of $t$).
We remark that this is not the case in general for $f \in C^\infty_0(\R)$.
This was the objet of \cite{frw1,frw2,rtrut}.
Later extensions to the multi-dimensional case were performed,
see \cite{issoglio, diel, cannizzaro}.
Generalizations to the path-dependent case were done
by \cite{ORT1_PartI, ORT1_Bessel}.
In all these cases the solutions are constituted by  continuous processes which are not necessarily semimartingales.
In the literature of jump processes the notion of
martingale problem has followed two different routes.
The first one is a new formulation of martingale problem
 given  in \cite{JacodBook} where the
formulation makes essentially use of the notion
of characteristics. This approach is particularly natural  in the purely discontinuous framework, see e.g. \cite{BandiniCalviaColaneri}.
The second one  continues the Stroock-Varadhan approach, see e.g.
\cite{LaachirRusso, barrasso3}.

We describe now the main contributions of the paper.
\begin{enumerate}
\item First we formulate a unique decomposition
  for a weak Dirichlet process $X = X^c + A$,
  where $X^c$ is a continuous local martingale
  and $A$ a martingale orthogonal process  vanishing at zero, 
  see Proposition \ref{P:uniqdec}.
  Until now unique decompositions of weak Dirichlet processes
  were established when $X$ has the special weak Dirichlet property,
  in particular if $X$ is a special semimartingale.
  Even when $X = M + V$ is a general c\`adl\`ag semimartingale,
  no unique decomposition was available because
  the property of $V$ to have  bounded variation was not
  enough to determine it. We recall that often purely
  discontinuous martingales have bounded variation.

\item  	

 Let $X$ be a c\`adl\`ag process  satisfying
\begin{equation}\label{int_small_jumps}
		\sum_{s \leq \cdot} |\Delta X_s|^2 
		< \infty \,\,\,\,\textup{a.s.}
\end{equation}
Notice that condition \eqref{int_small_jumps} is equivalent to ask that 
$ (1 \wedge |x|^2) \star \mu^X \in \mathcal A_{\textup{loc}}^+$
(see Proposition \ref{P:new}), where $\mu^X$ is the jump measure related to $X$ defined in \eqref{jumpmeasure}.
\\
	In Corollary \ref{C:3.21} we prove that if $X$ is 
 a weak Dirichlet process then it  
 is a special weak Dirichlet process if and only if 
		\begin{equation}\label{int_big_jumps_INTRO}
		  x \,\one_{\{|x| > 1\}} \star \mu^X \in \mathcal{A}_{\textup{loc}}.
		\end{equation}
					This result in particular extends the classical   characterization of a special semimartingale, see Proposition 2.29, Chapter II, in \cite{JacodBook}. We recall that if $X$ is a special weak Dirichlet process and a semimartingale, then it is a special semimartingale, see Proposition 5.14 in \cite{BandiniRusso1}.
\\
 More generally, let $v: \R_+ \times \R \rightarrow \R$  continuous. 
  In Theorem \ref{T:spw} we 
	  give a necessary and sufficient condition on the weak Dirichlet process $Y_t=v(t, X_t)$ to be a special weak Dirichlet process, namely
	  \begin{equation*}
	 (v(s,X_{s-}+x)-v(s,X_{s-})) \,\one_{\{|x| >1\}} \star \mu^X \in \mathcal{A}_{\textup{loc}}.
		\end{equation*}
		 Notice  that in the literature only sufficient conditions were available for $Y$ to be a special weak Dirichlet process, see for instance  Theorem 5.31 in \cite{BandiniRusso1}.
		 
\item   In Theorem  \ref{T:2.11} we provide a first chain rule expanding a process $v(s, X_s)$, where $v$ has no regularity at all, 
 that  extends  a similar chain rule established in \cite{BandiniRusso1} but only for purely jump processes. Indeed, under the conditions  \eqref{cond5.36} and \eqref{G1_cond}, if $v(t,X_t)$ is an $\F$-weak Dirichlet process with  unique continuous martingale component $Y^c$, then we have 
 \begin{align}\label{dec_wD1_intro}
Y &= Y^c + (v(s,X_{s-} + x)-v(s,X_{s-}))
\frac{k(x)}{x}\,\star (\mu^X- \nu^X) \notag\\
&+ \Gamma^{k}(v) + (v(s,X_{s-} + x)-v(s,X_{s-}))
\frac{x- k(x)}{x}\,
\star \mu^X, 
\end{align}
with
$ \Gamma^{k}(v)$ a predictable and   martingale orthogonal  process.
Our result constitutes an important tool to solve the identification problem  for   a BSDE driven by a random measure and a Brownian motion, see Remark \ref{R:previouschainrule}.

\item 
  We relax the notion of finite quadratic variation process, by giving the notion of weakly finite quadratic variation, see Definition \ref{D:weakfin}.
  This is a notion which is more related to the convergence in law of subsequences. The u.c.p convergence of \eqref{Appr_cov_ucpI} is replaced by
  the fact that, for every $T>0$,  $[X,X]_{0 < \varepsilon \le \varepsilon_0}^{ucp}(T)$ are tight
  for $\varepsilon_0 >0$ small enough.
  A classical example of weakly finite quadratic variation
  process comes up when 
  $$  \sup_{0 < \varepsilon \le \varepsilon_0}
  [X,X]_\varepsilon^{ucp}(T) < \infty \,\,\textup{a.s.},$$
  see Remark \ref{R:twoitems}-(i).
  Another example is given when $X$
  has finite energy, see Remark \ref{R:twoitems}-(ii).
  It is not difficult to exhibit
  a process with finite energy which has no quadratic variation, see Example \ref{E31}.
Notice that  condition \eqref{int_small_jumps} holds if $X$ is a finite quadratic variation process, and it is also valid under the weaker condition of $X$ being a weak finite quadratic variation process, see Proposition \ref{P:3.12}.

	\item Let $v \in C^{0,1}(\R_+ \times \R)$. If $X$ is a weakly finite quadratic variation process,  then Theorem
  \ref{P:3.10} states that 
  the process $Y = v(\cdot,X)$ is a weak Dirichlet process,
  and  identifies its unique  continuous local martingale component
  $Y^c$ of $Y$, as  
  	$$
	Y^c = Y_0 + \int_{0}^{\cdot}\partial_x v(s, X_{s})\,dX^c_s.
$$
The combined results in Theorems   \ref{T:2.11}  
 and  \ref{P:3.10} constitute indeed a $C^{0,1}$-type chain rule for   weak Dirichlet processes 
 which generalizes Theorem 5.15  
 in \cite{BandiniRusso1}.
As far as special weak Dirichlet processes are concerned, the corresponding $C^{0,1}$-type chain rule is given in  Corollary \ref{C:new} generalizing Theorem 5.31 in \cite{BandiniRusso1}, 
 see also Remark \ref{R:3.10bis}. 
Notice that in \cite{BandiniRusso1} 
  $X$ was also
   supposed to be a finite quadratic variation process. 

\item 
We introduce the notion of  {\it   characteristics} $(B^k, C, \nu)$ for weak Dirichlet processes (see Definition \ref{D:genchar}), that extends the corresponding one for semimartingales. 
We remark that, if $X$ is a weak Dirichlet process, then   $\Gamma^k(Id)=B^k\circ X$ with $ \Gamma^{k}(v)$  defined in \eqref{dec_wD1_intro}, see Corollary \ref{C:genchar} and  Remark \ref{R:cor}. Given the characteristics $(B^k, C, \nu)$ of a  weak Dirichlet process $X$, a natural question is to determine the characteristics of a process $h(\cdot, X)$, where $h \in C^{0,1}$. In fact, it is possible to  provide the second and third characteristic of $h(\cdot, X)$ in terms of  $C$ and $\nu$, see Remark \ref{R:332}, while it is a challenging problem to evaluate the first   characteristic. Nevertheless, we are able to solve this problem in the case when $h$ is bijective and time-homogeneous, see  Remark \ref{R:Rred}.
\item We introduce a  notion of martingale problem, which applies in a general framework  including   possibly non-Markovian jumps processes and non semimartingales, by generalizing the classical Stroock-Varadhan martingale problem with respect to some domain  $\mathcal D_{\mathcal A}\subseteq C^{0,1}$ (replacing $\mathcal D$) and operator $\mathcal A$ (replacing $\partial_ t + L_t$), see Definition \ref{D:mtpb1}. 
Moreover, here the Lebesgue measure  $dt$ can be substituted by some random kernel. 
\\
Let $X$ be a c\`adl\`ag weakly finite quadratic variation process. The fact that $X$ is a solution to some martingale problem in the sense of Definition \ref{D:mtpb1} does not imply that it is a semimartingale (indeed, we are in particular interested in the case where $X$ is not a semimartingale). 
Among others, if $X$ is a solution of a martingale problem with respect to $\mathcal D_{\mathcal A}$ and $\mathcal A$, with $\mathcal D_{\mathcal A}$ dense in $C^{0,1}$,  we even do not  know if $X$ is a weak Dirichlet process.  
Corollary \ref{C:3.30bis} provides some necessary and sufficient condition conditions under which $X$ is weak Dirichlet. 
To get those conditions, Theorem \ref{T:new4.4} together with Proposition \ref{T:3.30_bis}  give some crucial preparatory stochastic calculus tools.

\item Section \ref{S:mtpb_hom} relates our inhomogeneous formulation
of martingale problem with the (more classical) time-homogeneous
expression. Knowing that a process $X$ solves some martingale
problem (in the time-homogeneous sense), the fact that it also
solves a non-homogeneous martingale problem corresponds
to some general chain rule.

\item  In Section \ref{SSExamples} we discuss
five classes of examples of martingale problems. The first two are respectively the case of general semimartingales  and the case  where there is a bijective function $h$ in $C^{0,1}$ such that $h(t, X_t)$ is a  semimartingale.   
  The third one concerns discontinuous processes
  solving martingale problems with distributional drift.
  For this, existence and uniqueness is
  discussed systematically in the companion
  paper \cite{BandiniRusso_DistrDrift}.
  The fourth one is about continuous path-dependent
  problems involving distributional drifts.
  The latter one is about the martingale problem
  solved by a piecewise deterministic Markov  process. 
\end{enumerate}

 \section{Preliminaries and notations}\label{SPrelim}

 In the sequel we will consider the space of functions
$ 
u: \R_+ \times \R \rightarrow \R$, $(t,x)\mapsto u(t,x)$, 
 which are of class  $C^{0,1}$ or  $C^{1,2}$.
 $C^{0,1}_b$  (resp. $C^{1,2}_b$) stands for the class of bounded
 functions which belong to  $C^{0,1}$ (resp. $C^{1,2}$).
 $C^{0,1}$ is equipped with the topology of uniform convergence on each compact
 of $u$ and $\partial_x u$. 
$C^0$ (resp. $C^0_b$) will denote the space of continuous functions (resp. continuous and bounded  functions) on $\R$ equipped with the topology of uniform convergence on each compact (resp. equipped with the topology of uniform convergence).  $C^1$ (resp. $C^2$) will be the space of continuously differentiable (twice continuously differentiable) functions $u:\R\rightarrow \R$.  $C^{1}_b$  (resp. $C^{2}_b$) stands for the class of bounded
 functions which belong to  $C^{1}$ (resp. $C^{2}$). 
 $D(\R_+)$ will denote the space of real c\`adl\`ag functions on $\R_+$.  

Let $T>0$ be a finite horizon. 
$C^{0,1}([0,T]\times \R)$ will denote the space of functions in  $C^{0,1}$ restricted to $[0,T]\times \R$. In the following, $D(0,\,T)$ (resp.
 $D_{-}(0,\,T)$, $C(0,\,T)$, $C^1(0,\,T)$) 
will indicate  the space of real c\`adl\`ag
(resp. c\`agl\`ad,  continuous, continuously differentiable) functions on $[0,\,T]$. 
These spaces   
are equipped with the uniform norm.
We will also indicate by $||\cdot||_{\infty}$ the essential supremum norm and by $||\cdot||_{var}$ the total variation norm. 
 Given a topological space $E$, in the sequel $\mathcal{B}(E)$ will denote 
the Borel $\sigma$-field associated with $E$.

A stochastic basis  $(\Omega, \mathcal F, \mathbb F, \P)$ is fixed
throughout the section. We will  suppose that
$\mathbb F=(\mathcal F_t)$ satisfies the usual conditions.
By convention, any c\`adl\`ag process (or function) defined on $[0,\,T]$ is   extended to
$\R$ by continuity.  
A similar convention is made for random fields (or functions)  on $[0,T] \times \R$. 
Related to   $\mathbb F$,
the symbol $\mathbb{D}^{ucp}$  
will denote  the  space of all adapted c\`adl\`ag  
 processes endowed with the  u.c.p. (uniform convergence in probability) 
topology on each compact interval.

$\mathcal{P}$ (resp. $\mathcal{\tilde{P}}:=\mathcal{P}\otimes \mathcal{B}(\R)$) will denote the predictable $\sigma$-field on $\Omega \times \R_+$ (resp. on $\tilde{\Omega} := \Omega \times \R_+\times \R$).
For a random field $W$, the simplified notation  $W \in\mathcal{\tilde{P}}$ means that $W$ is $\mathcal{\tilde{P}}$-measurable.

A process $X$ indexed by $\R_+$  will be said to be with integrable variation if the expectation of its total variation is finite. 
$\mathcal{A}$ (resp. $\mathcal{A}_{\textup{loc}}$) will denote  the collection of all adapted processes with   integrable variation (resp.  with locally integrable variation), and    $\mathcal{A}^+$ (resp $\mathcal{A}_{\textup{loc}}^+$)  the collection of all adapted integrable increasing (resp. adapted locally integrable)  processes. 
The significance of locally is the usual one which refers 
to  localization by stopping times, see e.g. (0.39) of  
\cite{jacod_book}. 

The concept of random measure   
will be extensively used 	throughout   the paper. 
For a  detailed  discussion on this topic  and the unexplained  notations,
we refer to
Chapter I and Chapter II, Section 1, in \cite{JacodBook}, Chapter III in \cite{jacod_book},  and  Chapter XI, Section 1, in \cite{chineseBook}.
In particular, if $\mu$ is a random measure on $[0,\,T]\times \R$, for any measurable real function $H$ defined on $\Omega \times [0,\,T]$, one denotes $H \star \mu_t:= \int_{]0,\,t] \times \R} H(\cdot, s,x) \,\mu(\cdot, ds \,dx)$,
when the stochastic integral in the right-hand side is defined
(with possible infinite values).

We  recall that a transition kernel $Q(e, dx)$ of
a measurable space $(E, \mathcal E)$ into another measurable space $(G,\mathcal G)$ is a family
$\{Q(e, \cdot): e \in E\}$ of positive measures on $(G,\mathcal G)$, such that $Q(\cdot ,C)$ is $\mathcal E$-measurable for each $C \in \mathcal G$, see
for instance in Section 1.1, Chapter I of \cite{JacodBook}.

 Let $X$ be an adapted c\`adl\`ag process. 
 We set the 
 corresponding  jump measure  $\mu^X$ by
\begin{equation}\label{jumpmeasure}
\mu^X(dt\,dx)= \sum_{s >0} \one_{\{\Delta X_s \neq 0\}}\, \delta_{(s, \Delta X_s)}(dt\,dx). 
\end{equation}
We denote by  $\nu^X= \nu^{X, \P}$ the compensator of $\mu^X$,
         see \cite{JacodBook} (Theorem 1.8, Chapter II).
         The dependence on $\P$ will be omitted when self-explanatory.
 For any random field
$W$, 
we set 
\begin{align*} 
	\hat{W}_t = \int_{\R} W_t(x)\,\nu^X(\{t\} \times  dx), 
\quad  
	\tilde{W}_t = \int_{\R} W_t(x)\,\mu^X(\{t\}\times dx)
	  -\hat{W}_t, 
\end{align*}
whenever they are well-defined. 
	We also define
$$
C(W) := |W- \hat W|^2\star \nu^X + \sum_{s \leq \cdot}|\hat{W}_s|^2\,(1-\nu^X(\{s\}\times \R)), 
$$	
and, for every $q \in [1,\,\infty[$,   the linear spaces 
\begin{align*}
	\mathcal{G}^q(\mu^X)=\Big\{W \in \mathcal{\tilde{P}}
	:\,\, \forall s \geq 0 \,\int_\R |W(s,x)|\,\nu^X(\{s\} \times dx)< \infty, \,\,
	\Big[\sum_{s\leq \cdot}|\tilde{W}_s|^2\Big]^{q/2}\in \mathcal{A}^+\Big\},\quad 
	\\ 
	\mathcal{G}^q_{\textup{loc}}(\mu^X)=\Big\{W \in  \mathcal{\tilde{P}} 
	: \,\,   \forall s \geq 0 \,\int_\R |W(s,x)|\,\nu^X(\{s\} \times dx)< \infty, \,\,
	\Big[\sum_{s\leq \cdot}|\tilde{W}_s|^2\Big]^{q/2}\in \mathcal{A}_{\textup{loc}}^+\Big\}.\nonumber
\end{align*} 
For a random field $W$ on $[0,T] \times \R$ we set the norms 
$
||W||^2_{\mathcal{G}^2(\mu^X)}
:=\sper{C(W)_T}$,  
$||W||_{\mathcal{L}^2(\mu^X)}:=\E[  |W|^2 \star\nu_T]$, 
and the  space 
$\mathcal{L}^2(\mu^X):=\{W \in \tilde{\mathcal{P}} 
	:\,\, ||W||_{\mathcal{L}^2(\mu^X)}< \infty\}
	$.

If $W \in \mathcal G_{\textup{loc}}^1(\mu^X)$, we call stochastic integral with respect to $\mu^X-\nu^X$ and we denote it by $W \star(\mu^X-\nu^X)$, any purely discontinuous local martingale $X$ such that $\Delta X$ and $\tilde W$ are indistinguishable, see Definition 1.27, Chapter II, in \cite{JacodBook}. 
We recall that, if  $W \in \tilde {\mathcal P}$ such that $|W|\star \mu^X \in \mathcal{A}_{\textup{loc}}^+$, then $W \in \mathcal G_{\textup{loc}}^1(\mu^X)$ and 
	$
	W \star (\mu^X - \nu) = W \star \mu^X - W \star \nu^X$, see Theorem 1.28, Chapter II, in  \cite{JacodBook}. Moreover, by Theorem 11.21, point 3) in  \cite{chineseBook},
the following statements are equivalent:
\begin{enumerate}
	\item $W \in \mathcal{G}^2_{\textup{loc}}(\mu^X)$;
	\item 
$
C(W) \in \mathcal{A}_{\textup{loc}}^+
$;
\item $W \star (\mu^X - \nu^X)$ is a square integrable local martingale.
\end{enumerate}  
In this case 
$\langle W \star (\mu^X-\nu^X), W \star (\mu^X-\nu^X)\rangle = C(W)$.
Finally, if $W \in \mathcal{L}^2(\mu^X)$ 
then $W \in \mathcal{G}^2(\mu^X)$, and 
$C(W)= |W|^2 \star\nu^X- \sum_{s \leq\cdot}|\hat{W}_s|^2$.
In this case
	$||W||^2_{\mathcal{G}^2(\mu^X)} \leq ||W||^2_{\mathcal{L}^2(\mu^X)}$.

  \section{Weak Dirichlet processes: the suitable generalization of semimartingales with jumps}
  \label{S3}
  
\subsection{A new unique decomposition}

\label{S31}

A stochastic basis  $(\Omega, \mathcal F, \mathbb F, \P)$ is fixed
throughout the section. Sometimes the dependence on $\F$ will be omitted.
 Given an adapted (c\`adl\`ag) process $X$ on it,  we will denote by    $\mu^X$ its  jump measure given in \eqref{jumpmeasure} and by  $\nu^X$ the corresponding   compensator.

We recall that an $\mathbb F$-weak Dirichlet process is a process of the type 
\begin{equation} \label{GenDec}
X= M+ \Gamma,
\end{equation}
 where $M$ is an $\mathbb F$-local martingale and $\Gamma$ is an $\mathbb F$-orthogonal process vanishing at zero, while a special weak Dirichlet process is a weak Dirichlet process $X=M+\Gamma,$ where  $\Gamma$ is in addition predictable, see Definitions 5.5. and  5.6  in \cite{BandiniRusso1}.
For complementary results, the reader can consult Section 5 in \cite{BandiniRusso1}.  
\begin{remark}\label{R:decomp_mart}
	Any local martingale $M$ can be uniquely decomposed as the sum of a continuous local martingale   $M^c$ and a purely discontinuous local martingale $M^d$ such that $M^d_0 =0$,
see Theorem 4.18, Chapter I, in \cite{JacodBook}\end{remark}

The decomposition \eqref{GenDec} is not unique, but the result below proposes a particularly natural one, which is unique.

\begin{proposition}\label{P:uniqdec}
	Let $X$ be a c\`adl\`ag $\mathbb F$-weak Dirichlet process. Then there is a unique   continuous $\mathbb F$-local martingale  $X^c$ 
   and a unique  $\mathbb F$-martingale orthogonal process $A$ vanishing at zero, such that 
	\begin{equation}\label{E:dec}
		X=X^c + A.
\end{equation}
\end{proposition}
\proof
\noindent \emph{ Existence.} 
Since $X$ is an $\mathbb F$-weak Dirichlet process, by \eqref{GenDec} it is a process of the type $X= M+ \Gamma$, with $M$  an $\mathbb F$-local martingale and $\Gamma$  an $\mathbb F$-martingale orthogonal process vanishing at zero. 
Recalling Remark \ref{R:decomp_mart}, it follows that $X$ admits the decomposition 
\begin{equation}\label{weakDir}
		X= M^c + M^d + \Gamma,
	\end{equation}
	that provides \eqref{E:dec} by setting $A := M^d + \Gamma$ and $X^c:= M^c$.

\medskip

\noindent \emph{Uniqueness.} 
Assume that $X$ admits  the two decompositions
\begin{align*}
	X&= M^1 + A^1, \quad 
	X = M^2 + A^2, 
\end{align*}
with $M^1, M^2$ continuous  $\mathbb F$-local martingales and $A^1, A^2$  $\mathbb F$-martingale orthogonal processes vanishing at zero. 
So we have $0 = M^{1}-M^{2}  + A^1-A^2$. Taking the covariation of previous equality with $M^{1}-M^{2}$,
we get $  [M^{1}-M^{2},M^{1}-M^{2}] \equiv 0$. Since $M^1 - M^2$ is a continuous martingale vanishing at zero
we finally obtain $M^1 = M ^2$ and so  $A^1 = A^2$.
\endproof
\begin{remark} \label{R32}
	Notice that decomposition \eqref{weakDir} of the weak Dirichlet $X$  is not unique. 
\end{remark}

\begin{remark} \label{Edec}
  The unique  decomposition \eqref{E:dec} holds (and it is new)
  in particular for semimartingales.
  On the other hand, when $X$ is a semimartingale, the notion of $X^c$ (as unique continuous martingale component)
  was introduced in Definition 2.6, Chapter II, in \cite{JacodBook} (in that case it was fixed $X_0^c=0$). There, a unique decomposition was provided
only  after having fixed a truncation.
\end{remark}

\begin{proposition}\label{P:2.14}
		Let $X$ be an $\mathbb F$-semimartingale. Then $[X, X]^c = \langle X^c, X^c\rangle$.
\end{proposition}
\proof
By Proposition  \ref{P:uniqdec} we have the unique decomposition $X= X^c + A$. On the other hand, being $X$ a semimartingale, we also have that $X= M+ V$ with $V$  a bounded variation process, and $M$ a local martingale. Since the unique decomposition $M=M^c + M^d$ holds (see  Remark \ref{R:decomp_mart}),
we get that $X^c=M^c$  and 
 \begin{align}\label{GammaDef}
A = M^d + V.	
 \end{align}
 Indeed  $M^d $ and $V$ are $\F$-martingale orthogonal processes.
 Now, by \eqref{QVC} and \eqref{E:dec}
 \begin{align}\label{cov1}
[X,X] =[X, X]^c  + \sum_{s \leq \cdot} |\Delta A_s|^2.
\end{align}
On the other hand, the bilinearity of the covariation gives 
 \begin{align}\label{cov2}
[X,X] = [X^c, X^c] + [A, A]
\end{align}
also taking into account that $A$ is an
$\F$-martingale orthogonal process.
To conclude, comparing \eqref{cov1} and \eqref{cov2}, we have to show that 
$
[A, A]=\sum_{s \leq \cdot} |\Delta A_s|^2
$.
Now, by \eqref{GammaDef},  Proposition 5.3 and Proposition 2.14 in \cite{BandiniRusso1},
\begin{align*}
[A, A] &= [M^d, M^d]+ [V, V] + 2 [M^d, V]=\sum_{s \leq \cdot} (|\Delta M^d_s|^2 + |\Delta V_s|^2 + 2 \Delta M^d_s \Delta V_s)\\
&= \sum_{s \leq \cdot} |\Delta M^d_s + \Delta V_s|^2= \sum_{s \leq \cdot} |\Delta A_s|^2.
\end{align*}
\qed
\begin{remark}
  The equality $[X,X]^c = \langle X^c, X^c\rangle$ valid for semimartingales (see Proposition \ref{P:2.14}) may be not true for weak Dirichlet processes, even if they are of finite quadratic variation. Indeed, let $W$, $B$ be two canonical
  Brownian motions, and $\F$ be the canonical filtration associated with $W$ and $B$. We set 
	$$
	X_t = \int_0^t B_{t-s}\, dW_s. 
	$$
	By Proposition 2.10 of \cite{er2}, $X$ is an $\F$-weak Dirichlet and $\F$-martingale orthogonal. By Proposition \ref{P:uniqdec} it follows that  $X^c\equiv 0$.  On the other hand, by Remark 2.16-(2) in \cite{er2}, $[X, X]_t = \frac{t^2}{2}$ which is different from zero.  
\end{remark}

	\subsection{Fundamental chain rules}

From here on we will denote by $\mathcal K$ the set of truncation functions, namely
$$
\mathcal K:=\{k:\R \rightarrow \R \textup{ bounded with compact support: } k(x) =x \textup{ in a neighborhood of } 0\}.
$$
 A typical choice of $k(x)$ will be $k(x) = x \one_{\{|x|\leq 1\}}$. In this case, $\frac{x-k(x)}{x}= \one_{\{|x| >1\}}$. 
       We will make use  the following assumption on
a pair $(v, X)$, with  $v:\R_+ \times \R \rightarrow \R$ locally bounded and
$X$  a  c\`adl\`ag process.
\begin{hypothesis}\label{H:3.7-3.8}
\begin{align}
	 	&v(t, X_t) \,\,\textup{is a c\`adl\`ag process, and  for every}\,\, t \in \R_+, \,\,\Delta v(t, X_t)= v(t, X_{t}) - v(t, X_{t-});\label{cond5.36}\\
&\textup{$\exists \,k \in \mathcal K$ \,\,such that}\,\,\,\,\,(v(s,X_{s-} + x)-v(s,X_{s-}))\,\frac{k(x)}{x}\in {\mathcal G}^1_{\textup{loc}}(\mu^X). \label{G1_cond}
\end{align}
\end{hypothesis}
	\begin{remark}\label{R:cont}
	\begin{itemize}
		\item [(i)]
	If $v$ is continuous, then the pair $(v, X)$ obviously fulfills
 \eqref{cond5.36}. For a more refined condition on $(v, X)$ to guarantee the validity of \eqref{cond5.36} we refer to Hypothesis 5.34 in \cite{BandiniRusso1}. 
 \item[(ii)] Assume the validity of \eqref{cond5.36}. If there is $a \in\R_+$ such that  $\sum_{s \leq \cdot }|\Delta v(s, X_s)|\one_{\{|\Delta  X_s|\leq a\}} <\infty$ a.s., then \eqref{G1_cond} is verified. This is trivially verified if $v(\cdot, X)$ is a bounded variation process. 
 \end{itemize}
	\end{remark}
	\begin{proposition}\label{R:th_chainrule}
Let	$X$ be an adapted c\`adl\`ag process satisfying \eqref{int_small_jumps}.
Let $v:\R_+ \times \R \rightarrow \R$  be  a function of class $C^{0,1}$. Then  Hypothesis \ref{H:3.7-3.8}  holds true. 
\end{proposition}
\proof
Condition \eqref{cond5.36} holds true being $v$ continuous, see Remark \ref{R:cont}-(i). 
On the other hand, by Proposition \ref{P:2.10}, 
$$
|v(s,X_{s-} + x)-v(s,X_{s-})|^2\,\frac{k^2(x)}{x^2}\star \mu^X \in \mathcal{A}^+_{\rm loc}, \quad \forall k \in \mathcal K.
$$
 In particular,  $(v(s,X_{s-} + x)-v(s,X_{s-}))\,\frac{k(x)}{x} \in \mathcal G^2_{\textup{loc}}(\mu^X)$, that in turn  implies condition \eqref{G1_cond}, being $\mathcal G^2_{\textup{loc}}(\mu^X) \subseteq \mathcal G^1_{\textup{loc}}(\mu^X)$. 
\endproof

\begin{theorem}\label{T:2.11}
	Let $X$  be a c\`adl\`ag and  $\F$-adapted process satisfying \eqref{int_small_jumps}.
Let  $v:\R_+ \times \R \rightarrow \R$ be a locally bounded function such that $(v, X)$ satisfies Hypothesis \ref{cond5.36}.
Let $Y_t= v(t,X_t)$ be an $\F$-weak Dirichlet process with  
continuous martingale component $Y^c$. 
  Then,
for every $k \in \mathcal K$, one can write the decomposition 
\begin{equation}\label{dec_wD1}
Y = Y^c + M^{k,d} + \Gamma^{k}(v) + (v(s,X_{s-} + x)-v(s,X_{s-}))
\frac{x- k(x)}{x}\,
\star \mu^X, 
\end{equation}
with
\begin{align}\label{Mdk}
M^{k,d}&:= (v(s,X_{s-} + x)-v(s,X_{s-}))
\frac{k(x)}{x}\,\star (\mu^X- \nu^X)
\end{align}
and $ \Gamma^{k}(v)$ a predictable and  $\F$-martingale orthogonal  process.
\end{theorem}

\begin{remark}\label{R:previouschainrule}
\begin{itemize}
\item[(i)] Sufficient conditions for $Y$ to be weak Dirichlet are given in Theorem \ref{P:3.10}.
	\item [(ii)] Theorem \ref{T:2.11} drastically generalizes the results given in Proposition 5.37 in \cite{BandiniRusso1}, where we considered the case of an $\F$-martingale orthogonal process $Y_t= v(t, X_t)$  and such that $\sum_{s \leq \cdot} |\Delta Y_s| \in \mathcal A_{\textup{loc}}^+$; in particular, $\sum_{s \leq T} |\Delta Y_s| < \infty$ a.s. Notice that, by Remark \ref{R:cont}-(ii), Hypothesis \ref{H:3.7-3.8} is verified. 
In that case $Y^c=0$
by Proposition \ref{P:uniqdec}.
\end{itemize}
\end{remark}
\begin{remark}\label{R:linmapGamma}
Let 
$\mathcal D$ be the set of functions $v\in \R_+ \times \R \rightarrow \R$  locally bounded such that $v(\cdot, X)$ is a weak Dirichlet process and   $(v, X)$ satisfies Hypothesis \ref{H:3.7-3.8}.
We observe that 
$v \mapsto \Gamma^k(v)$ in Theorem \ref{T:2.11} can be seen as a linear map from  $\mathcal D$
 to the space of c\`adl\`ag  adapted processes.
			In the sequel  $\mathcal D$ will also denote a similar set of functions  defined on $[0,T]\times \R$ instead of $\R_+ \times \R$, and  related to a process $X=(X(t))_{t \in [0,T]}$. 
 \end{remark}

\noindent \emph{Proof of Theorem \ref{T:2.11}.}  
Let $k \in \mathcal K$. We claim that  
\begin{equation}\label{Yk}
Y^k := Y 
-(v(s,X_{s-} + x)-v(s,X_{s-}))
\frac{x- k(x)}{x}
\star \mu^X
\end{equation}
is an $\mathbb F$-special weak Dirichlet process. 
Indeed, we notice that,
	\begin{equation}\label{eqAbis}
		Y^k=   Y^c + A^Y -(v(s,X_{s-} + x)-v(s,X_{s-}))
\frac{x- k(x)}{x}\star \mu^X=: Y^c + \tilde A^k,
	\end{equation}
	where $\tilde A^k$  is an $\mathbb F$-martingale orthogonal process.
By condition \eqref{G1_cond} and Lemma \ref{R:B}, $(v(s,X_{s-} + x)-v(s,X_{s-}))
\frac{k(x)}{x} \in \mathcal G^1_{\textup{loc}}(\mu^X)$. 
Then the process $M^{k, d}$ in \eqref{Mdk}, is a 
   purely discontinuous
  $\mathbb F$-local martingale, see Section \ref{SPrelim}. 
	We rewrite \eqref{eqAbis} as
	\begin{equation}\label{tot_jumps}
		Y^k= Y^c + M^{k,d} + \Gamma^{k}(v),
	\end{equation}
	with 
	\begin{equation}\label{Gammak}
			\Gamma^{k}(v):= \tilde A^k- M^{k,d}.
	\end{equation}
	This gives  
	$
	\Delta Y^k_t = \Delta  M^{k,d}_t + \Delta \Gamma^{k}_t(v)
	$.
	Now,  by \eqref{Yk} and \eqref{Mdk}
we get
	\begin{align*}
		\Delta Y^k_t &= \int_{\R} (v(s,X_{s-}+x)-v(s,X_{s-})) \,\frac{k(x)}{x} \,\mu^X(\{t\}\times\,dx),\\
		\Delta M^{k,d}_t &= \int_{\R} (v(s,X_{s-}+x)-v(s,X_{s-})) \,\frac{k(x)}{x} \,\mu^X(\{t\}\times\,dx)\\
		&- \int_{\R} (v(s,X_{s-}+x)-v(s,X_{s-})) \,\frac{k(x)}{x} \,\nu^X(\{t\}\times\,dx),
	\end{align*}
so that, by \eqref{tot_jumps},
	\begin{align*}
		\Delta \Gamma^{k}_t(v)&= \int_{\R} (v(s,X_{s-}+x)-v(s,X_{s-})) \,\frac{k(x)}{x} \,\nu^X(\{t\}\times\,dx).
	\end{align*}
	Recall  that an adapted c\`adl\`ag process  is predictable if and only if its jump process is predictable, see Remark \ref{R:pred}-2. below. We conclude that $\Gamma^{k}(v)$ is an $\F$-predictable process; moreover, by \eqref{Gammak},  $\Gamma^k(v)$ is also an $\F$-martingale orthogonal process.
This yields  decomposition \eqref{dec_wD1}.
\endproof
\begin{remark}\label{R:pred}
\begin{enumerate}
	\item 

Any  c\`agl\`ad process is locally bounded, see the lines above Theorem 15,
Chapter IV, in 
\cite{protter}.

\item 	Let $X$ be a c\`adl\`ag adapted process. Recalling that  $\Delta X_s = X_{s} - X_{s-}$, we have the following. 
	\begin{itemize}
		\item $X$ is locally bounded  if and only if $\Delta X$ is locally bounded. As a matter of fact,  $(X_{s-})$ is a c\`agl\`ad process and therefore locally integrable.
		\item $X$ is predictable if and only if $\Delta X$ is predictable. Indeed,  $(X_{s-})$ is a predictable process being adapted and left-continuous.
	\end{itemize} 
	\end{enumerate}
\end{remark}
Taking  $k(x) = x \one_{\{|x|\leq a\}}$ in Theorem \ref{T:2.11} we get the following result. 
        \begin{corollary}\label{P:boundedjumps}
          Let $X$ be a c\`adl\`ag  $\F$-adapted  process
         satisfying \eqref{int_small_jumps}. Let $v:\R_+ \times \R \rightarrow \R$ be a locally bounded function, 
          such that            $(v,X)$ satisfies 
          Hypothesis \ref{cond5.36}.
          Assume moreover that, for some $a \in \R_+$, 
\begin{equation}\label{deltaXa}
		|\Delta X_t |\leq a, \quad \forall  t \in \R_+. 
\end{equation}
Then, if  $Y_t = v(t,X_t)$ is an $\F$-weak Dirichlet process, then it is an  $\F$-special weak Dirichlet process.
\end{corollary}

The particular case of Corollary \ref{P:boundedjumps} with $v \equiv \textup{Id}$ is stated below.
\begin{corollary}\label{C:gencharbis}
	Let  $X$ be a c\`adl\`ag and  $\F$-adapted process satisfying \eqref{int_small_jumps} and \eqref{deltaXa}. Then, if   $X$ is an $\F$-weak Dirichlet process, it is an $\F$-special weak Dirichlet process.
\end{corollary}
\begin{theorem}\label{T:spw}
	Let $X$ be c\`adl\`ag and $\F$-adapted  process satisfying \eqref{int_small_jumps}.  Let $v:\R_+ \times \R \rightarrow \R$ be a locally bounded function such that $(v, X)$ satisfies Hypothesis \ref{cond5.36}.
	Set $Y_t = v(t,X_t)$, and assume that $Y$ is an $\F$-weak Dirichlet process. 
	Then $Y$ is an $\F$-special weak Dirichlet process if and only if 
	\begin{equation}\label{intY}
\textup{$\exists \,a \in \R_+$ s.t.}\,\,\,	 (v(s,X_{s-}+x)-v(s,X_{s-})) \,\one_{\{|x| >a\}} \star \mu^X \in \mathcal{A}_{\textup{loc}}.
	\end{equation}
\end{theorem}

\begin{remark}\label{R:boundspecfirst}
	If $v$ is bounded and $X$  is a c\`adl\`ag and  $\F$-adapted, 
	 then condition \eqref{intY} is satisfied because of  Lemma \ref{L:c}.  
\end{remark}

\noindent \emph{Proof of Theorem \ref{T:spw}.}
Let $a \in \R_+$ and set 
$$
\tilde Y:= \sum_{s \leq \cdot} \Delta Y_s \,\one_{\{|\Delta X_s| >a\}}= (v(s,X_{s-}+x)-v(s,X_{s-})) \,\one_{\{|x| >a\}} \star \mu^X. 
$$
By Theorem \ref{T:2.11} with $k(x)= \one_{\{|x| \leq a\}}$, 
$
Y - \tilde Y
$
is an $\mathbb F$-special weak Dirichlet process. It follows that $Y$  is an $\mathbb F$-special weak Dirichlet process if and only if  
$
\tilde Y
$
is an $\mathbb F$-special weak Dirichlet process. 
$\tilde Y$ has bounded variation,  so  it is a semimartingale, therefore it is a special semimartingale,   see Proposition 5.14 in \cite{BandiniRusso1}. 
 This can be shown to be equivalent to  condition \eqref{intY} by making use of the  first three equivalent items of 
 Proposition 4.23, Chapter I,   in \cite{JacodBook}.
\qed

 Corollary \ref{C:3.21} below follows from Theorem \ref{T:spw}
taking  $v \equiv \textup{Id}$.  
It extends a characterization stated in Proposition 5.24 in \cite{BandiniRusso1}: thereby, $X$ was supposed to belong to  a particular class of weak Dirichlet processes.

\begin{corollary}\label{C:3.21}
Let $X$ be an $\F$-weak Dirichlet process satisfying \eqref{int_small_jumps}.  Then $X$ is an $\F$-special weak Dirichlet process if and only if 
	\begin{equation}\label{int_big_jumps}
	\exists \,a  \in \R_+ \,\,\textup{s.t.}\,\,  x \,\one_{\{|x| > a\}} \star \mu^X \in \mathcal{A}_{\textup{loc}}.
	\end{equation}
      \end{corollary}
The following result  is the analogous of Theorem \ref{T:2.11} for  special weak Dirichlet processes.

\begin{theorem}\label{T:2.11bis}
	Let $X$  be a c\`adl\`ag and  $\F$-adapted process satisfying \eqref{int_small_jumps}.
Let  $v:\R_+ \times \R \rightarrow \R$ be a locally bounded function such that $(v, X)$ satisfies Hypothesis \ref{cond5.36}.
Let $Y_t= v(t,X_t)$ be an $\F$-weak Dirichlet process  
with continuous martingale component $Y^c$.
Assume moreover that the pair $(v, X)$
 satisfies condition \eqref{intY}. Then  $Y_t = v(t,X_t)$ is an  $\F$-special weak Dirichlet process, admitting the unique decomposition 
 \begin{align}\label{Spec_weak_chainrule1}
Y&= Y^c + (v(s,X_{s-} + x)-v(s,X_{s-}))\,\star (\mu^X- \nu^X) +\Gamma(v), 
\end{align}
 with 
 $\Gamma(v)$ a predictable and $\F$-martingale orthogonal process.
 Moreover,  \begin{equation}\label{GammaY}
	\Gamma(v) = \Gamma^k(v) + 
	(v(s,X_{s-} + x)-v(s,X_{s-}))
\frac{x- k(x)}{x}\,\,\star  \nu^X,
\end{equation}
with $\Gamma^k(v)$ the predictable and $\F$-martingale orthogonal process  appearing in \eqref{dec_wD1}.
\end{theorem}
\proof
Thanks to condition  \eqref{intY},  by Theorem \ref{T:spw} $Y$ is an $\F$-special weak Dirichlet process. 
By  Theorem \ref{T:2.11},
\begin{align}\label{3.9}
Y &= Y^c +  (v(s,X_{s-} + x)-v(s,X_{s-}))
\frac{k(x)}{x}\,\star (\mu^X- \nu^X) + \Gamma^{k}(v)\notag\\
& + (v(s,X_{s-} + x)-v(s,X_{s-}))
\frac{x- k(x)}{x}\,
\star \mu^X.  
\end{align}
We add and subtract in \eqref{3.9}  the term 
$
(v(s,X_{s-} + x)-v(s,X_{s-}))
\frac{x- k(x)}{x}\,\,\star  \nu^X$. 
Recalling that, for every random field $W$,   $|W|\star \mu^X \in \mathcal{A}_{\textup{loc}}^+$ implies $W \in \mathcal G_{\textup{loc}}^1(\mu^X)$ and 
	$
	W \star (\mu^X - \nu^X) = W \star \mu^X - W \star \nu^X$ (see Section \ref{SPrelim}), 
we get decomposition \eqref{Spec_weak_chainrule1} with $\Gamma(v)$ provided by \eqref{GammaY}. 
\endproof

\begin{remark}\label{R:3items}
\begin{enumerate}
\item 
It directly follows  from \eqref{Spec_weak_chainrule1} that 
	$$
	\Delta \Gamma_s(v) = \int_\R	(v(s,X_{s-} + x)-v(s,X_{s-}))
\star  \nu^X(\{s\}\times dx).
	$$

\item 
Let  $X$  be a c\`adl\`ag and  $\F$-adapted process satisfying \eqref{int_small_jumps}. Let $v: \R_+ \times \R \rightarrow \R$ be  such that $(v, X)$ satisfies Hypothesis \ref{cond5.36}.
If $v$ is moreover a bounded function,  by Remark \ref{R:boundspecfirst} condition  \eqref{intY} is satisfied as well,  and all the assumptions of Theorem \ref{T:2.11bis}  are verified.  
	\end{enumerate}
\end{remark}

\begin{remark}\label{R: BSDE}
Theorem \ref{T:2.11bis} is  in particular useful to solve the so-called identification problem for BDSEs. 
Proposition 5.37 in \cite{BandiniRusso1} allowed to solve the identification problem when the BSDE was  exclusively driven by  the random measure, see Theorem 3.14 in \cite{BandiniRusso2}. 

Let $\zeta$ be a non-decreasing,  adapted and continuous process, and  $\lambda$ be a predictable random measure on $\Omega \times [0,T] \times \R$. 
We consider a  BSDE driven by a random measure $\mu-\nu$ and a continuous martingale $M$ of the type
\begin{align}\label{BSDE}
	Y_t &= \xi + \int_{]t, T]}\tilde g(s, Y_{s-}, Z_s) d\zeta_s + \int_{]t, T]\times \R}\tilde f(s, Y_{s-},  U_s(e)) \lambda(ds\,de) \notag\\
	&- \int_{]t, T]}Z_s dM_s - \int_{]t, T]\times \R}U_s(e)(\mu-\nu)(ds\,de),
\end{align}
whose solution is a triple of processes $(Y, Z, U(\cdot))$.
If $Y_t=v(t, X_t)$ for some function $v$ and some adapted c\`adl\`ag process $X$, the identification problem consists in expressing $Z$ and $U(\cdot)$ in terms of $v$. 
\begin{itemize}
\item [(i)]
 Being $Y_t = v(t, X_t)$  a solution to a BSDE, it is a special weak Dirichlet process
  (even a special semimartingale),
   and  therefore $(v, X)$  satisfies  condition \eqref{intY}. Therefore, if Hypothesis \ref{cond5.36} holds for $(v, X)$, then   Theorem \ref{T:2.11bis} allows  to identify $U(\cdot)$. More precisely, this provides 
	$$
	U(e)\star (\mu-\nu)= (v(s,X_{s-} + x)-v(s,X_{s-}))\,\star (\mu^X- \nu^X), \quad \textup{a.s.}
	$$ 
	From now on  suppose $\mu=\mu^X$, even though this can be generalized, 	see Hypothesis 2.9 in \cite{BandiniRusso2}. This yields
	$$
	H(x)\star (\mu^X- \nu^X)=0, \quad \textup{a.s.}, 
	$$ 
	with $H(x):=U(x)- (v(s,X_{s-} + x)-v(s,X_{s-}))$.
	If $H \in \mathcal G^2_{\textup{loc}}(\mu^X)$, then the predictable bracket of  $H(x)\star (\mu^X- \nu^X)$ is well-defined and equals $C(H)$. Since $C(H)_T=0$, we get (see Proposition 2.8 in \cite{BandiniRusso2}) that there is a predictable process $(l_s)$ such that 
	$$
	H_s(x)= l_s \one_K(s)\quad d\P\, \nu^X(ds\,dx) \,\,\textup{a.e.},
	$$
	where  $K:=\{(\omega, t):\,\,\nu^X(\omega, \{t\} \times \R)=1\}$.
	 	\item [(ii)] 
In order to identify the process $Z$ we need that $v$ belongs to $C^{0,1}([0,T]\times \R)$ and that $X$ is a  weakly finite quadratic variation process, and this will be discussed in Remark \ref{R:BSDEid2}. 

\end{itemize}
\end{remark}

Theorems \ref{T:2.11} and \ref{T:2.11bis} with $v \equiv \textup{Id}$ give in particular the following result.
\begin{corollary}\label{C:genchar}
	Let  $X$ be an $\F$-weak Dirichlet process satisfying condition \eqref{int_small_jumps}.
	Let $X^c$ be the continuous martingale part of $X$. Then, the following holds.  
	\begin{itemize}
		\item[(i)]	
 Let $k \in \mathcal K$. Then $X$ can be decomposed as
	\begin{equation}\label{dec_X}
X = X^c + k(x)\,\star (\mu^X- \nu^X) + \Gamma^{k}(Id) +
(x- k(x))\,
\star \mu^X,
\end{equation} 
with  $\Gamma^{k}$  the  operator introduced in Theorem \ref{T:2.11}.
\item[(ii)]If \eqref{int_big_jumps} holds, then the $\F$-special weak Dirichlet process $X$ admits the decomposition 
	\begin{equation}\label{dec_X_spec}
X = X^c + x\,\star (\mu^X- \nu^X) + \Gamma 
\end{equation} 
with 
$\Gamma := \Gamma^{k}(Id) + 
	x \,
\one_{\{|x| >1\}}\,\star  \nu^X$.
	\end{itemize}
\end{corollary}

\subsection{The notion of  characteristics for weak Dirichlet processes} \label{S:genchar}

	 In the present section we provide a generalization of the concept of characteristics for semimartingales (see Appendix \ref{A:D}) in the case  of weak Dirichlet processes. 
\begin{remark}\label{R:cor}	 
Let  $X$ be   an $\F$-weak Dirichlet process with jump measure $\mu^X$ satisfying condition \eqref{int_small_jumps}. Given $k \in \mathcal K$, by Corollary \ref{C:genchar}-(i), we have that 
 $
 X^k= X- \sum_{s \leq \cdot} [\Delta X_s - k(\Delta X_s)]
 $
  is an $\F$-special weak Dirichlet process with unique decomposition
\begin{equation}\label{dec_Xk}
	X^k= 
	X^c + k(x) \star (\mu^X- \nu^X) + B^{k, X}, 
\end{equation}
where  
\begin{itemize}
\item  $X^c$ is the unique continuous $\mathbb F$-local martingale part of $X$ introduced in  Proposition \ref{P:uniqdec};
\item   $B^{k, X}:=\Gamma^{k}(Id)$, which  is in particular a predictable and  $\mathbb F$-martingale orthogonal process.  
\end{itemize}
\end{remark}
\begin{remark}
When $X$ is  a semimartingale, 	$B^{k, X}$ is a bounded variation process, so in particular $\mathbb F$-martingale orthogonal. 
\end{remark}

From here on we denote by $\check \Omega$ the canonical  space of all c\`adl\`ag functions $\check \omega: \R_+ \rightarrow \R$, namely $\check \Omega= D(\R_+)$, and   by $\check X$ the canonical process defined by  $\check X_t(\omega)= \check \omega(t)$. We also set $\check {\mathcal  F}= \sigma(\check X)$, and $\check {\mathbb  F}= (\check {\mathcal  F}_t)_{t \geq 0}$. Let $\mu$ be the jump measure of $\check X$ and $\nu$ the compensator of $\mu$ under the law $\mathcal L(X)$ of $X$. 
\begin{definition}\label{D:genchar}
We call \emph{characteristics} of $X$, associated with $k \in \mathcal K$, 
  the triplet $(B^k,C= \langle \check X^c, \check X^c\rangle,\nu)$ on $(\check \Omega, \check {\mathcal F}, \check {\mathbb F})$ obtained from the unique decomposition \eqref{dec_Xk} for $\check X$ under $\mathcal L(X)$. 
  In particular, 
		 \begin{itemize}
	\item[(i)] $B^k$ is a predictable and $\check{\mathbb F}$-martingale orthogonal process, with $B_0^k=0$;
	\item[(ii)] $C$ is an  $\check{\mathbb F}$-predictable and  increasing process, with $C_0=0$; 	
	\item[(iii)] $\nu$ is an $\check{\mathbb F}$-predictable random measure on $\R_+ \times \R$.
	\end{itemize}
	\end{definition}
	\begin{remark}\label{R:3.27}
	\begin{itemize}
	\item [a)]  $\mu \circ X$ is the jump measure of $X$, and its compensator  under $\P$ is $\nu \circ X$, see Proposition 10.39-b) in \cite{jacod_book}.
\item [b)] It is not difficult to show that  $\check X^c \circ X$ is a continuous local martingale under $\P$. 
\item [c)] Let $Y$ and $Z$ be two processes on $(\check \Omega, \check {\mathcal F}, \check {\mathbb F})$ such that $[Y, Z]$ exists under $\mathcal L(X)$. Then $[Y, Z]\circ X= [Y \circ X, Z \circ X]$ under $\P$. In particular, $B^k \circ X$	is an $\mathbb F$-martingale orthogonal process. 
\item [d)] By previous items  we have a new decomposition of the process $X^k$: 
$$
X^k = \check X^c \circ X + k \star (\mu \circ X-\nu \circ X) + B^k \circ X.
$$
\item [e)] By the uniqueness of decomposition \eqref{dec_Xk} we have  $X^c= \check X^c \circ X$ and $B^{k, X}=B^k \circ X$.
	\item [f)]$C^X:= C \circ X=  \langle X^c, X^c\rangle$ by c) and e).
\end{itemize}
	\end{remark}

\begin{remark}
Assume that $X$ admits  characteristics $(B^{k_0}, C, \nu)$ depending to a given  truncation function $k_0$. Then, if we choose another truncation function $k$, the corresponding characteristics will be $(B^{k},  C,  \nu)$ with 
	$$
	B^{k} \circ X= B^{k_0} \circ X + (k - k_0)\star (\nu \circ X),
	$$
	where  $(k - k_0)\star (\nu \circ X) \in \mathcal{A}_{\textup{loc}}$, see Lemma \ref{L:D4}. 
\end{remark}

\begin{remark}
\begin{itemize}
\item[a)]	Given a process  $X \in \mathcal A_{\textup{loc}}^+$, Theorem 3.17, Chapter I,  in \cite{JacodBook} shows that there is a unique predictable process $X^p \in \mathcal A_{\textup{loc}}^+$, called \emph{compensator} of $X$, such that $X-X^p$ is a local martingale. In particular, $X$ is a special semimartingale.
\item [b)]The notion of compensator can be naturally extended. Given a process  $X$, we denote by $X^p$ a process verifying the following conditions: 
\begin{enumerate}
	\item[(i)] $X^p$ is predictable;
	\item[(ii)] $X^p$ is martingale orthogonal;
	\item[(iii)] $X- X^p$ is a local martingale.
\end{enumerate}
Obviously, such $X^p$ exists if and only if $X$ is a special weak Dirichlet process, and $X^p=\Gamma$, with $\Gamma$ the process in Corollary \ref{C:genchar}-(ii), so it is uniquely determined. 
\item[c)] The two notions of $X^p$ in a) and b) coincide when $X\in \mathcal A_{\textup{loc}}^+$. Indeed, a special semimartingale is a special weak Dirichlet process. 
\end{itemize}
\end{remark}

\subsection{A  weak notion of finite quadratic variation and a  stability theorem}
	\label{S32}

	In the following we will use the notation (for $t \geq 0$)
	\begin{align}
	[X,X]_{\varepsilon}^{ucp}(t)&:=\int_0^t \frac{(X_{(s+ \varepsilon) \wedge t}- X_s)^2}{\varepsilon} ds, \quad \varepsilon>0,\label{ucpbrack}\\
	C_\varepsilon(X,X)(t) &:=\int_0^t \frac{(X_{s+ \varepsilon}- X_s)^2}{\varepsilon} ds, \quad \varepsilon>0. \label{contbrack}
\end{align}
From now on $T>0$ will denote a fixed maturity time. 
\begin{definition}\label{D:weakfin}
	A c\`adl\`ag process $X=(X(t))_{t \in [0,T]}$ is said to be a weakly finite quadratic variation process if there is $\varepsilon_0 >0$ such that the laws of the random variables 
	$[X,X]_{\varepsilon}^{ucp}(T)$, $0<\varepsilon\leq \varepsilon_0$,
are tight. 
\end{definition}
Below, $\varepsilon>0$ will mean $0<\varepsilon \leq \varepsilon_0$ for some $\varepsilon_0$ small enough. 
For instance, a family $(Z_\varepsilon)_{\varepsilon >0}$  of random variables will indicate a sequence  $(Z_\varepsilon)_{0<\varepsilon \leq \varepsilon_0}$ for some $\varepsilon_0$ small enough. 
\begin{remark}
	A finite quadratic variation process is a weakly finite quadratic variation process. Indeed, if $\int_0^\cdot \frac{(X_{(s+ \varepsilon) \wedge \cdot}- X_s)^2}{\varepsilon} ds$ converges u.c.p., the random variable  $[X,X]_{\varepsilon}^{ucp}(T)$ converges in probability, and so it also converges in law.  
\end{remark}

\begin{proposition}\label{P:criteriumtight}
	Let  $(Z_\varepsilon)_{\varepsilon>0}$ be a nonnegative sequence of random variables.
	Suppose that one of the two items below hold. 
	\begin{itemize}
		\item [(i)]
	$\sup_{\varepsilon
	 >0}Z_\varepsilon< \infty$ a.s.
		\item [(ii)]  $\sup_{\varepsilon
	 >0}\E[Z_\varepsilon]< \infty$.
	\end{itemize}   
	 Then the family of distribution of $(Z_\varepsilon)_{\varepsilon >0}$ is tight. 
\end{proposition}
\proof
Let  $\delta>0$.  In order to prove the result, we need to find $M >0$ such that 
$$
\P\left(Z_\varepsilon\notin [0, M]\right)\leq \delta, \quad \forall \varepsilon>0.
$$

\noindent (i)   We choose $M=M(\delta)>0 $ such that 
$
\P\left(\sup_{\varepsilon
	 >0} Z_\varepsilon>M\right)\leq \delta$.

\noindent (ii) For every $M$, by Markov-Chebishev inequality, 
$$
\P\left( Z_\varepsilon>M\right)\leq \frac{\E[ Z_\varepsilon]}{M}\leq \frac{1}{M}\sup_{\varepsilon >0} \E[ Z_\varepsilon].
$$
We choose $M$ so that the previous upper bound is smaller or equal then $\delta$.
\endproof
\begin{remark}\label{R:twoitems}
It follows from Proposition \ref{P:criteriumtight} that a c\`adl\`ag process $X$ is a weakly finite quadratic variation process  when one of the following conditions holds:
\begin{enumerate}
	\item [(i)] $\sup_{\varepsilon
	 >0}\int_0^T \frac{(X_{s+ \varepsilon}- X_s)^2}{\varepsilon} ds < \infty$ a.s. 	 \item [(ii)]$\sup_{\varepsilon
	 >0}\E \Big[\int_0^T \frac{(X_{s+ \varepsilon}- X_s)^2}{\varepsilon} ds\Big]  < \infty$. A process $X$ fulfilling this condition was called a finite energy process, see  \cite{cjms} (in the framework of F\"ollmer's discretizations) and \cite{rv4}.
\end{enumerate}
\end{remark}

\begin{example}\label{E31}
Let $X$ be a square integrable  process with weakly stationary increments. Set $V(a):=\E \left[ (X_{a}- X_0)^2 \right] $, $a \in \R_+$. Let  $V(\varepsilon)=O(\varepsilon)$, $\varepsilon>0$, i.e., there is a constant  $C>0$ such that $V(\varepsilon) \leq C \varepsilon$ for a small $\varepsilon>0$.  Then condition (ii) of Remark \ref{R:twoitems} is verified. A classical example of a process that satisfies that condition is a weak Brownian motion of order $\kappa=2$, see \cite{fwy}.  In this case $V(a)=a$. Indeed the bivariate distributions of such a process are the same as those of Brownian motion. In general, such a process is not  a semimartingale and not even   a finite quadratic variation process.
\end{example}
 
\begin{proposition}\label{P:3.12}
 Suppose that $X$ is a weakly finite quadratic variation process. Then condition \eqref{int_small_jumps} holds true. 
\end{proposition}
\begin{remark}\label{R:equiv}
By Proposition \ref{P:new}, condition \eqref{int_small_jumps} is equivalent to $(1 \wedge |x|^2 )\star \nu^X \in  \mathcal {A}_{\textup{loc}}^+$.  The validity of latter condition  was known for semimartingales,  see  Theorem  2.42, Chapter II, in \cite{JacodBook}.
\end{remark}
\begin{proof}
Let $\gamma$ 
be a constant. 
For each fixed $\omega$, let $\tau_{0}=\tau_{0, \gamma}:=0$ and set
$$
\tau_{i}(\omega)=\tau_{i, \gamma}(\omega):=\inf_{t > \tau_{i-1, \gamma}(\omega)} \{t: \,\,|\Delta X_t|>\gamma\}, \quad i \in \N, 
$$
with the convention that it is $+ \infty$ if the set is empty. 
Notice that, almost surely, there exists a finite number of jumps of $X$ greater than $\gamma$.  
Let thus  $N= N(\omega)$ 
be such that $\tau_{N, \gamma}$ is the maximum of those jump times. 
 Let $\varepsilon >0$, and  define 
 $
 \Omega_{\varepsilon,\gamma}^0=\{\omega \in \Omega: \tau_{i}(\omega) - \tau_{i-1}
 (\omega) >\varepsilon, \,\,i \in \N\}$,
 with the convention that $\infty - \infty = \infty$. We have that 
 \begin{equation}\label{unionOmega}
\cup_\varepsilon \Omega_{\varepsilon, \gamma}^0=\Omega,
 \end{equation}
   up to a null set. 
We  set 
$$
J_A(\varepsilon):= \sum_{i=1}^N \frac{1}{\varepsilon} \int_{\tau_i - \varepsilon}^{\tau_i} (X_{(s + \varepsilon) \wedge T}- X_s)^2 ds.
$$
On $\Omega_{\varepsilon, \gamma}^0$ we have
\begin{equation}\label{estbrackucp}
[X,X]_{\varepsilon}^{ucp}(T) \geq J_A(\varepsilon). 
\end{equation}
We also know, by Lemma 2.11 in \cite{BandiniRusso1}, that 
\begin{equation}\label{JAconv}
 J_A(\varepsilon)\underset{\varepsilon \rightarrow 0}{\rightarrow} \sum_{i=1}^N \one_{]0,T]}(\tau_i) |\Delta X_{\tau_i}|^2=\sum_{t \leq T}|\Delta X_t|^2 \one_{\{|\Delta X_t| >\gamma\}}\,\,\,\,\textup{a.s.}
\end{equation}
Let $\kappa >0$. Since $X$ is a weakly finite quadratic variation process, there exists $\ell=\ell(\kappa)>0$ such that, for every $\varepsilon$ (small enough), $\P(\Omega_{\varepsilon, \ell}^c) < \kappa$ with
$
	\Omega_{\varepsilon, \ell}:=\{\omega \in \Omega:\,\,[X,X]_{\varepsilon}^{ucp}(T) \leq \ell\}
$.
We get 
\begin{align}\label{Psumjumps}
	&\P\Big(\sum_{t \leq T}|\Delta X_t |^2 \one_{\{|\Delta X_t |>\gamma\}} > K\Big)\leq \kappa + \P\Big(\sum_{t \leq T}|\Delta X_t |^2 \one_{\{|\Delta X_t |>\gamma\}} > K; \Omega_{\varepsilon,\ell}\Big)\notag\\
& \leq \kappa +\P\Big(|J_A(\varepsilon)|  > \frac{K}{2}; \Omega_{\varepsilon,\ell}\Big)+ \P\Big(\Big|\sum_{t \leq T}|\Delta X_t |^2\one_{\{|\Delta X_t |>\gamma\}}-J_A(\varepsilon)\Big|  >\frac{K}{2}; \Omega_{\varepsilon,\ell}\Big)\notag\\
& \leq \kappa +\P\Big(
|J_A(\varepsilon)|  > \frac{K}{2}
; \Omega_{\varepsilon,\ell}\cap \Omega^0_{\varepsilon,\gamma}
\Big)+\P\Big(|J_A(\varepsilon)|  > \frac{K}{2}; \Omega_{\varepsilon,\ell} \setminus \Omega^0_{\varepsilon,\gamma}\Big)\notag\\
&+ \P\Big(\Big|\sum_{t \leq T}|\Delta X_t |^2\one_{\{|\Delta X_t |>\gamma\}}-J_A(\varepsilon)\Big| >\frac{K}{2}\Big)\notag\\
& \leq \kappa +\P\Big(
[X,X]_{\varepsilon}^{ucp}(T)
  > \frac{K}{2}
; \Omega_{\varepsilon,\ell}\cap \Omega^0_{\varepsilon,\gamma}
\Big)+\P\Big( (\Omega^0_{\varepsilon,\gamma})^c\Big)\notag\\
&+ \P\Big(\Big|\sum_{t \leq T}|\Delta X_t |^2\one_{\{|\Delta X_t |>\gamma\}}-J_A(\varepsilon)\Big| >\frac{K}{2}\Big),
\end{align}
where we have used   
\eqref{estbrackucp} in the latter inequality.
The first probability in the right-hand side of  latter inequality   equals 
$$
\P\Big(
[X,X]_{\varepsilon}^{ucp}(T)
\wedge \ell  > \frac{K}{2}
; \Omega_{\varepsilon,\ell}\cap \Omega^0_{\varepsilon,\gamma}
\Big). 
$$
Choosing $K=K(\ell, \kappa)$ so that $\frac{K}{2} \leq \ell$,  this probability  is zero. Therefore,
applying the $\limsup_{\varepsilon \rightarrow 0}$ in \eqref{Psumjumps}, and taking into account  \eqref{unionOmega}  and \eqref{JAconv}, we get 
\begin{align}\label{f:gamma}
\P\Big(\sum_{t \leq T}|\Delta X_t |^2 \one_{\{|\Delta X_t |>\gamma\}} > K\Big) \leq \kappa.
\end{align}
Notice that $\sum_{t \leq T}|\Delta X_t |^2 \one_{\{|\Delta X_t |>\gamma\}}$ converges increasingly to $\sum_{t \leq T}|\Delta X_t |^2$, a.s. when $\gamma$ converges to zero.
Consequently, letting $\gamma \rightarrow 0$ in \eqref{f:gamma}, we obtain
$$
\P\Big(\sum_{t \leq T}|\Delta X_t |^2  > K\Big) \leq \kappa.
$$
Finally, 
\begin{align*}
	\P\Big(\sum_{t \leq T}|\Delta X_t |^2  = \infty\Big)\leq \P\Big(\sum_{t \leq T}|\Delta X_t |^2  > K\Big) \leq \kappa,
\end{align*}
so the conclusion follows. 
\end{proof}

Below we give a significant generalization of Proposition 3.10 in \cite{gr}, where the result was proven when $X$ is continuous and of finite quadratic variation. 
 When $X$ is c\`adl\`ag, even
  in the case when $X$ is a finite quadratic variation process,
  the result is new. 
Crucial tools to prove the result are the canonical  decomposition  stated in Proposition \ref{P:uniqdec} and Proposition \ref{P: weakfinconv}.

\begin{theorem}\label{P:3.10}
	Let $X$ be an $\mathbb F$-weak Dirichlet process with weakly finite quadratic variation. Let $v\in C^{0,1}([0,T] \times \R)$.  Then  $Y_t = v(t, X_t)$ is an $\mathbb F$-weak Dirichlet with continuous martingale component 
	\begin{align}\label{Yc}
	Y^c = Y_0 + \int_{0}^{\cdot}\partial_x v(s, X_{s})\,dX^c_s.
\end{align}
\end{theorem}

Theorem \ref{P:3.10},   together with Theorems \ref{T:spw} and  \ref{T:2.11bis} (recall Proposition \ref{R:th_chainrule}), provides the following result. 
\begin{corollary}\label{C:new}
	Let $X$ be an $\mathbb F$-weak Dirichlet process with weakly finite quadratic variation. Let $v\in C^{0,1}([0,T]\times \R)$   
such that $(v, X)$ satisfies 
  condition \eqref{intY}. Then $Y_t = v(t,X_t)$ is an $\F$-special weak Dirichlet process,  admitting the unique decomposition 
 \begin{align}\label{Spec_weak_chainrule2}
Y&= Y_0 + \int_{0}^{\cdot}\partial_x v(s, X_{s})\,dX^c_s + (v(s,X_{s-} + x)-v(s,X_{s-}))\,\star (\mu^X- \nu^X) +\Gamma(v), 
\end{align}
 with 
 $\Gamma(v)$ a predictable and $\F$-martingale orthogonal process.
\end{corollary}
\begin{remark}\label{R:3.10bis}
Corollary 	\ref{C:new} extends  the  chain rules   previously given in this framework in  Theorems 5.15 and 5.31 in \cite{BandiniRusso1}. More precisely, we have the following.
\begin{itemize}
	\item 	In Theorem  5.15 in \cite{BandiniRusso1} we already proved that, if $X$ is a $\F$-weak Dirichlet process  of  finite quadratic variation, and $v$ is of class $C^{0,1}$, then $Y=v(\cdot, X)$ is again a weak Dirichlet process. However   the decomposition established therein  was  not unique. Here instead in Theorem \ref{T:2.11}, together with Theorem \ref{P:3.10}, we are able to provide  an explicit form of its unique decomposition.
	\item Theorem 5.31 in \cite{BandiniRusso1} focused on sufficient conditions on $(v, X)$ so that $Y$ is  special weak Dirichlet. 
  Here,  by Corollary \ref{C:new},
    we are able to give the necessary and sufficient condition \eqref{intY} on $(v, X)$ such that $Y$ is  
a special weak Dirichlet, and we provide its unique  decomposition.  
\end{itemize}

\end{remark}
\begin{remark}\label{R:BSDEid2}
Let us consider the BSDE \eqref{BSDE} introduced in 
 Remark \ref{R: BSDE}. 
If $Y_t = v(t, X_t)$ is  a solution to \eqref{BSDE}, it is a special weak Dirichlet process, and  therefore $(v, X)$  satisfies  condition \eqref{intY}. Then, if 
$v\in C^{0,1}([0,T] \times \R)$, then 
Corollary  \ref{C:new} allows also  to identify $Z$.
More precisely, we get 
$$
Z_t =\partial_x v(t, X_t) \frac{d \langle X^c, M\rangle_t }{\langle M\rangle_t}, \quad d\P \,d \langle M \rangle_t\textup{-a.e.}
$$ 
\end{remark}

\noindent \emph{Proof of Theorem \ref{P:3.10}.}
We aim at proving that, for every $\mathbb F$-continuous local martingale $N$,  
\begin{align}\label{toproveY}
	[v(\cdot, X), N]_t = \int_0^t \partial_x v(s, X_s)\, d[X^c, N]_s, \quad t \in [0,\,T].
\end{align}
Indeed, if \eqref{toproveY} were true, it would imply that 
$
A(v) := v(\cdot, X) - Y^c
$
 is martingale orthogonal,
 and therefore by additivity $v(\cdot, X)$ would be a weak Dirichlet process. Then  \eqref{Yc} would follow by the uniqueness of the continuous martingale part of $Y$. 

Let us thus prove \eqref{toproveY}.
The approximating sequence of the left-hand side of \eqref{toproveY} is
\begin{align*}
	\int_0^t [v({(s+ \varepsilon) \wedge t}, X_{(s+ \varepsilon) \wedge t})-v(s, X_s)]\, \frac{N_{(s+ \varepsilon) \wedge t}- N_s}{\varepsilon} ds= 
I_1(t, \varepsilon) + I_2(t, \varepsilon),
\end{align*}
with
\begin{align*}
	I_1(t, \varepsilon) &:= \int_0^t [v({(s+ \varepsilon) \wedge t}, X_{(s+ \varepsilon) \wedge t})-v((s+ \varepsilon) \wedge t, X_s)]\, \frac{N_{(s+ \varepsilon) \wedge t}- N_s}{\varepsilon} ds,\\
	I_2(t, \varepsilon) &:= \int_0^t [v({(s+ \varepsilon) \wedge t}, X_{s})-v(s, X_s)]\, \frac{N_{(s+ \varepsilon) \wedge t}- N_s}{\varepsilon} ds.
\end{align*}
Concerning $I_2(t, \varepsilon)$,  by stochastic Fubini theorem we get  
\begin{align}\label{boudaryterm}
	I_2(t, \varepsilon) 
& =\frac{1}{\varepsilon}\int_0^t [v({(s+ \varepsilon) \wedge t}, X_{s})-v(s, X_s)]\,\int_s^{(s + \varepsilon) \wedge t} d N_u\notag\\
& =\int_0^t d N_u  \int_{(u-\varepsilon)^+}^u [v({(s+ \varepsilon) \wedge t}, X_{s})-v(s, X_s)] \frac{ds}{\varepsilon}.
\end{align}
Since 
$$
\int_0^T d [N,N]_u \Big( \int_{(u-\varepsilon)^+}^u [v({s+ \varepsilon}, X_{s})-v(s, X_s)] \frac{ds}{\varepsilon}\Big)^2 \underset{\varepsilon \rightarrow 0}{\rightarrow} 0\quad \textup{in probability},
$$
by Problem 2.27,  Chapter 3, in \cite{ks}, this is enough to conclude that $I_2(\cdot, \varepsilon) \underset{\varepsilon \rightarrow 0}{\rightarrow} 0$  u.c.p.   (since $N$ is continuous, it is clear that we can neglect the  ``$\wedge t$"  in \eqref{boudaryterm}).
Concerning $I_1(t, \varepsilon)$, we have 
\begin{align}\label{I1a+I1b}
	&I_1(t, \varepsilon) 
	= \int_0^t [v({(s+ \varepsilon) \wedge t}, X_{(s+ \varepsilon)\wedge t})-v((s+ \varepsilon) \wedge t, X_s)]\, \frac{N_{(s+ \varepsilon) \wedge t}- N_s}{\varepsilon} ds\notag\\
	&= \int_0^t \int_0^1 \partial_x v((s+ \varepsilon) \wedge t, X_s +  a (X_{(s+ \varepsilon) \wedge t}-X_s))\, da\,  (X_{(s+ \varepsilon) \wedge t}-X_s)(N_{(s+ \varepsilon) \wedge t}- N_s)\frac{ds}{\varepsilon}\notag \\
	&= \int_0^t  \partial_x v(s, X_s)\,  (X_{(s+ \varepsilon) \wedge t}-X_s)(N_{(s+ \varepsilon) \wedge t}- N_s)\frac{ds}{\varepsilon}\notag \\
	&+\int_0^t \int_0^1 [\partial_x v((s+ \varepsilon) \wedge t, X_s + a (X_{(s+ \varepsilon) \wedge t}-X_s))-\partial_x v(s, X_{s})]\, da\,  (X_{(s+ \varepsilon) \wedge t}-X_s)(N_{(s+ \varepsilon) \wedge t}- N_s)\frac{ds}{\varepsilon} \notag\\
&=:	 I_{1,a}(t, \varepsilon)+ 	 I_{1,b}(t, \varepsilon).
\end{align}
We set
$
\tilde I_{1,a}(t, \varepsilon) := \int_0^t  \partial_x v(s, X_{s})\,  (X_{s+ \varepsilon}-X_s)(N_{s+ \varepsilon}- N_s)\frac{ds}{\varepsilon}$.
By Proposition \ref{P: weakfinconv} with $g(s)= \partial_x v(s, X_{s-})$, we get	
$$
	\tilde I_{1,a}(\cdot, \varepsilon)\underset{\varepsilon \rightarrow0}{\rightarrow} \,\,\int_0^\cdot\partial_x v(s, X_s)\, d[X^c, N]_s, \quad \textup{u.c.p}.
	$$
This also shows the convergence of $ I_{1,a}(\cdot, \varepsilon)$ to the same limit,  
since $ I_{1,a}(\cdot, \varepsilon) - \tilde  I_{1,a}(\cdot, \varepsilon) \underset{\varepsilon \rightarrow 0}{\rightarrow} 0$ u.c.p., being $N$ continuous.

It remains to prove that $I_{1,b}(\cdot, \varepsilon)\underset{\varepsilon \rightarrow 0}{\rightarrow} 0$ u.c.p. 
 Let  $(\varepsilon_n)_n$ be a sequence converging to zero as $n$ goes to infinity. We fix  $\kappa >0$. Recall that $X$ and $N$ have weakly finite quadratic variation. Then,  there exists $\ell= \ell(\kappa)>0$ such that, for every $n$ (big enough), there is an event $\Omega_{n, \ell}$ such that
$\P(\Omega_{n, \ell}^c) \leq \kappa$, and for $\omega \in \Omega_{n, \ell}$, 
\begin{align}\label{ellest}
	&\Big(\sup_{s \in [0,T]}|X_s(\omega)|^2+\int_0^T (X_{(s+\varepsilon_n) \wedge T}(\omega)-X_s(\omega))^2\frac{ds}{\varepsilon_n}  \Big)
	\notag\\
	&+\Big(\sup_{s \in [0,T]}|N_s(\omega)|^2+\int_0^T (N_{(s+\varepsilon_n) \wedge T}(\omega)-N_s(\omega))^2\frac{ds}{\varepsilon_n}\Big) \leq \ell. 
\end{align}
We provide now some estimates which are valid only for  $\omega \in \Omega_{n, \ell}$. For this, we proceed as in the proof of Proposition 5.18 in \cite{BandiniRusso1} (estimate of the term $I_{13}$).  We enumerate the jumps of $X(\omega)$ on $[0,T]$ by $(t_i) _{i \geq 0}$, and  
\begin{equation}\label{KX}
	\mathbb K^X= \mathbb K^X(\omega) \quad \textup{is smallest convex compact set including $\{X_t(\omega): t \in [0,\,T]\}$}.
\end{equation}
 We fix $\gamma >0$, and we choose $M= M(\gamma, \omega)$ such that 
$
\sum_{i=M+1}^\infty |\Delta X_{t_i}|^2 \leq \gamma^2$. 
For $\varepsilon>0$ small enough and depending on $\omega$, we introduce 
$
	B(\varepsilon,M) = \bigcup_{i=1}^{M} \, ]t_{i-1},t_i - \varepsilon]
$
and we decompose 
$
I_{1,b}(t, \varepsilon_n)=J^A(t,\varepsilon_n)+ J^B(t,\varepsilon_n),
$
with
\begin{align*} 
J^A(t,\varepsilon_n)&=\sum_{i = 1}^M \int_{t_i - \varepsilon_n}^{t_i}\frac{ds} 
{\varepsilon_n}\, 
\one_{]0,\,t]}(s)\,(X_{(s+ \varepsilon_n)\wedge t}-X_s)(N_{(s+ \varepsilon_n)\wedge t}-N_s)\cdot \\ 
&\quad \cdot\int_0^1 (\partial_x v((s+ \varepsilon_n)\wedge t,\,X_s + a (X_{(s+ \varepsilon_n)\wedge t}-X_s))-\partial_x v((s+ \varepsilon_n)\wedge t,\,X_s))\,da,\nonumber\\ 
J^B(t,\varepsilon_n)&= \frac{1} 
{\varepsilon_n}\int_{]0,\,t]} 
(X_{(s+ \varepsilon_n)\wedge t}-X_s)(N_{(s+ \varepsilon_n)\wedge t}-N_s)\,R^B(\varepsilon_n,s,t,M)\,ds, 
\end{align*} 
where 
$$ 
R^B(\varepsilon_n,s,t, M)= \one_{B(\varepsilon_n,M)}(s)\int_0^1 [\partial_x v((s+ \varepsilon_n)\wedge t,\,X_s + a (X_{(s+ \varepsilon_n)\wedge t}-X_s))-\partial_x v((s+ \varepsilon_n)\wedge t,\,X_s)]\,da. 
$$ 
 Let $\delta$ denote the modulus of continuity. 
 By Remark 3.12 in \cite{BandiniRusso1} we have for every $s, t \in [0,T]$, 
\begin{align*} 
R^B(\varepsilon_n,s,t,M) &\leq \delta\Big( \partial_x v \bigg|_{[0,\,T]\times \mathbb{K}^X} 
,\,\sup_{l} \sup_{\underset{|r-a| \leq \varepsilon_n}{r, a \in [t_{i-1},\,t_{i}]}} |X_{a}-X_{r}|\Big), 
\end{align*} 
so that Lemma 2.12 in \cite{BandiniRusso1} applied successively to the intervals $[t_{i-1},\,t_{i}]$ implies 
\begin{align} \label{Bill}
R^B(\varepsilon_n,s,t,M) &\leq \delta\big( \partial_x v \big|_{[0,\,T]\times \mathbb{K}^X} 
,3\gamma\big). 
\end{align} 
This concludes the estimates restricted to $\omega \in \Omega_{n,\ell}$. 

Since $N$ is continuous, by \eqref{ellest} (we remember that  $\ell$ is fixed), 
\begin{equation}\label{convJa}
	\sup_{t \leq T}|J^A(t,\varepsilon_n)| \one_{\Omega_{n,\ell}} \leq  \sqrt\ell \,\delta(N(\omega),\varepsilon_n) \,M(\gamma, \omega) \sup_{(t,x) \in [0,\,T] \times \mathbb K^X(\omega)} |\partial_x v|\underset{n \rightarrow \infty}{\rightarrow} \,\,0, \quad \textup{a.s.}
      \end{equation}
On the other hand,
we remark that 
\begin{align}\label{Xen}
	\int_0^t (X_{(s+\varepsilon_n) \wedge t}(\omega)-X_s(\omega))^2\frac{ds}{\varepsilon_n}  
	&= \int_0^{t-\varepsilon_n} (X_{s+\varepsilon_n}(\omega)-X_s(\omega))^2\frac{ds}{\varepsilon_n}  +\int_{t-\varepsilon_n}^t (X_{t}(\omega)-X_s(\omega))^2\frac{ds}{\varepsilon_n} \notag\\
	&\leq  \int_0^{T} (X_{(s+\varepsilon_n) \wedge T}(\omega)-X_s(\omega))^2\frac{ds}{\varepsilon_n} + \sup_{s \in [0,T]} |X_s(\omega)|^2. 
\end{align}

A similar estimate holds replacing $X$ with for $N$.   Consequently, recalling \eqref{Bill},
we have 
\begin{align}\label{3.34}
	&\sup_{t \in [0,\,T]} |J^B(t,\varepsilon_n)| \one_{\Omega_{n,\ell}}\notag\\
	&\leq \delta(\partial_x v|_{[0,\,T] \times \mathbb K^X(\omega)}, 3 \gamma)\sup_{t \in [0,\,T]}\sqrt{\int_0^t (X_{(s+\varepsilon) \wedge t}(\omega)-X_s(\omega))^2\frac{ds}{\varepsilon_n}  \int_0^t (N_{(s+\varepsilon) \wedge t}(\omega)-N_s(\omega))^2\frac{ds}{\varepsilon_n} }\notag\\
	&\leq \delta(\partial_x v|_{[0,\,T] \times \mathbb K^X(\omega)}, 3 \gamma)\cdot\notag\\
	&\cdot \sqrt{\Big(\sup_{s \in [0,T]}|X_s(\omega)|^2+\int_0^T (X_{(s+\varepsilon_n) \wedge T}(\omega)-X_s(\omega))^2\frac{ds}{\varepsilon_n}  \Big)\!\!\Big(\sup_{s \in [0,T]}|N_s(\omega)|^2+\int_0^T (N_{(s+\varepsilon_n) \wedge T}(\omega)-N_s(\omega))^2\frac{ds}{\varepsilon_n}\Big)}\notag\\
	&\leq  \ell\,\delta(\partial_x v|_{[0,\,T] \times \mathbb K^X(\omega)}, 3 \gamma),
\end{align}
where in the second inequality we have used   \eqref{Xen}.

We continue now the proof of the u.c.p. convergence of $I_{1,b}(\cdot, \varepsilon_n)$. Let  $K>0$.  By \eqref{I1a+I1b} we have 
\begin{align*}
	&\P\Big(\sup_{t \in [0,\,T]}|I_{1,b}(t, \varepsilon_n)| >K\Big)\leq \P\Big(\Omega_{n, \ell}^c\Big) + \P\Big(\sup_{t \in [0,\,T]}|I_{1,b}(t, \varepsilon_n)| >K, \Omega_{n, \ell}\Big)\\
	&\leq \kappa + \P\Big(\sup_{t \in [0,\,T]}|J^A(\varepsilon_n, t)| >\frac{K}{2}, \Omega_{n, \ell}\Big)+ \P\Big(\sup_{t \in [0,\,T]}|J^B(\varepsilon_n, t)| >\frac{K}{2}, \Omega_{n, \ell}\Big)\\
	& \leq \kappa + \P\Big(\sup_{t \in [0,\,T]}|J^A(\varepsilon_n, t)| >\frac{K}{2}, \Omega_{n, \ell}\Big) + \P\Big(\frac{K}{2} < \ell\,\delta(\partial_x v|_{[0,\,T] \times \mathbb K^X(\omega)}, 3 \gamma)\Big),
\end{align*}
where in the latter inequality we have used \eqref{3.34}. 
This holds true for fixed   $\kappa$, $\ell(\kappa)$, $\gamma$, $K$. 
So, taking into account \eqref{convJa}, 
\begin{align*}
	\limsup_{n \rightarrow \infty} \P\Big(\sup_{t \in [0,\,T]}|I_{1,b}(t, \varepsilon_n)| >K\Big) \leq \kappa  + \P\Big(\frac{K}{2} < \ell\,\delta(\partial_x v|_{[0,\,T] \times \mathbb K^X(\omega)}, 3 \gamma)\Big).
\end{align*}
We let now $\gamma \rightarrow 0$. We get $\delta(\partial_x v|_{[0,\,T] \times \mathbb K^X(\omega)}, 3 \gamma) \underset{\gamma \rightarrow 0}{\rightarrow} 0$ a.s., so that 
$$
\P\Big(\frac{K}{2} < \ell\,\delta(\partial_x v|_{[0,\,T] \times \mathbb K^X(\omega)}, 3 \gamma)\Big) \underset{\gamma \rightarrow 0}{\rightarrow} 0.
$$
Consequently, 
\begin{align*}
\limsup_{n \rightarrow \infty} 	\P\Big(\sup_{t \in [0,\,T]}|I_{1,b}(t, \varepsilon_n)| >K\Big)&\leq \kappa,
\end{align*}
and since $\kappa$ is arbitrary, this concludes the proof.
\qed

The result below follows from 
Theorem 
\ref{P:3.10}, together with Remark \ref{R:3items}-2. and Propositions \ref{P:3.12} and \ref{R:th_chainrule}.
\begin{corollary}\label{C:spw2}
	Let $X$ be an $\mathbb F$-weak Dirichlet with weakly finite quadratic variation.
Let $v\in C^{0,1}([0,T]\times \R)$ and bounded,  and set $Y_t = v(t,X_t)$.
Then  $Y$ is  a $\F$-special weak Dirichlet process.
\end{corollary}
	
\begin{remark}\label{R:pushforward}
	Let $X$ be a c\`adl\`ag process, $v:[0,T] \times \R \rightarrow \R $ a continuous function, and set $Y=v(\cdot, X)$. 
	It is well-known that, for fixed $\omega \in \Omega$,  $\mu^Y(\omega, \cdot)$  is the push forward of $\mu^X(\omega, \cdot)$ through  
 	  	$\mathcal H_{\omega}: (s,x) \mapsto v(s,X_{s-}(\omega)+x)-v(s,X_{s-}(\omega))$:
 	  \begin{align}
\mu^Y(]0,\,t]\times A)&=\int_{ ]0,\,t]\times \R}\one_{A \setminus 0}  (v(s,X_{s-}+ x)- v(s,X_{s-})) \, \mu^X(ds\, dx),\label{mubar}
 \end{align}
for all $A \in \mathcal B(\R)$.
In particular,
$$
 \Delta Y_t = \int_{  \R}y \, \mu^Y(\{t\}\times dy) =\int_{  \R}  (v(t, X_{t-}+ x) -v(t, X_{t-}))  \,\mu^X(\{t\} \times dx). 
 $$
Taking the predictable projection in identity \eqref{mubar}, we get that
 	  \begin{align}
\nu^Y(]0,\,t]\times A)&=\int_{ ]0,\,t]\times \R}\one_{A \setminus 0}  (v(s,X_{s-}+ x)- v(s,X_{s-})) \, \nu^X(ds\, dx),\label{nubar}
 \end{align}
for all $A \in \mathcal B(\R)$.
\end{remark}

\begin{remark}\label{R:332}
Let $X$ be a weak Dirichlet process of weakly finite quadratic variation with  given characteristics $(B^k,\nu, C)$, and $v \in C^{0,1}([0,T] \times \R)$. 
\begin{itemize}
\item[(i)]
By  Theorem \ref{P:3.10}, $Y_t= v(t, X_t)$ is a weak Dirichlet process, so   it admits  characteristics $(\bar B^{k}, \bar C, \bar \nu)$.
Moreover, again by  Theorem \ref{P:3.10}, recalling Remark \ref{R:3.27}-f), 
 \begin{align}
 	\bar C \circ Y = \langle Y^c, Y^c \rangle 
 	&= \int_{]0,\,\cdot]}|\partial_x v(s, X_{s})|^2 d \langle X^c, X^c \rangle_s = \int_{]0,\,\cdot]}|\partial_x v(s, X_{s})|^2 d (C\circ X)_s.\label{Cbar}
 	 \end{align}
 If moreover $v(s, \cdot)$ is bijective for every $s$, from \eqref{Cbar}-\eqref{nubar} we get the explicit form of $\bar C$ and $\bar \nu$: 
 \begin{align}
 	&(\bar C \circ Y)_t 
 	= \int_{]0,\,t]}|\partial_x v(s, v^{-1}(s,Y_{s}))|^2 d (C \circ v^{-1}(\cdot,Y))_s,\label{Cbardef}\\
 	&(\bar\nu \circ Y)(]0,\,t]\times A)\notag\\
 	&=\int_{ ]0,\,t]\times \R}\one_{A \setminus 0} \, (v(s,v^{-1}(s,Y_{s-})+ x)- v(s,v^{-1}(s,Y_{s-}))) \, (\nu \circ v^{-1}(\cdot,Y))(ds\, dx),\label{nubardef}
  \end{align}
  for all $A \in \mathcal B(\R)$.
   \item[(ii)]
 In general, instead,  it does not look possible to give explicitly the characteristic $\bar B^k$ of $Y$ in terms of the corresponding characteristic of $X$, even when $\bar B^k$ is a bounded variation process. This can however be done for instance when $v$ is a    bijective homogeneous function of class $C^2$,
   see  Remark \ref{R:Rred}-(ii). 
 \end{itemize}
\end{remark}

\begin{remark}\label{R:Rred}
 Let $X$ be a weak Dirichlet process with weakly finite quadratic variation.
Let $h: \R \rightarrow \R$ in $C^1$ and bijective.
By Theorem \ref{P:3.10},    $Y_t := h(X_t)$ is a weak Dirichlet process.
\begin{itemize}
\item[(i)] 
By Theorem \ref{T:2.11}.
  we have 
\begin{align}\label{firstex}
Y&= Y^c +  \Gamma^{k}(h)+  (h(X_{s-}+x)-h(X_{s-}))\frac{x-k(x)}{x}\star \mu^X \notag\\
&+  (h(X_{s-}+x)-h(X_{s-}))\frac{k(x)}{x}\star (\mu^X-\nu^X),  
\end{align}
with $ \Gamma^{k}(h)$ a predictable and  $\F$-martingale orthogonal  process.
The characteristic $\bar B^{k}$ of $Y$ can be determined in terms of $\nu^X$ and of the map $\Gamma^k(h)$.

As a matter of fact,  
by Corollary \ref{C:genchar}-(i) together with Remark \ref{R:cor}, recalling \eqref{mubar}-\eqref{nubar}, we have
\begin{align}\label{secondex}
Y&= Y^c +  \bar B^{k,Y}+  (y-k(y))\star \mu^Y +  k(y)\star (\mu^Y-\nu^Y)\notag\\
&=Y^c +  \bar B^{k,Y}+  (h(X_{s-}+x)-h(X_{s-})-k(h(X_{s-}+x)-h(X_{s-})))\star \mu^X\notag\\
& +  k(h(X_{s-}+x)-h(X_{s-}))\star (\mu^X-\nu^X)\notag\\
&=Y^c +  \bar B^{k,Y}+  \Big[(h(X_{s-}+x)-h(X_{s-}))\frac{k(x)}{x}-k(h(X_{s-}+x)-h(X_{s-}))\Big]\star \mu^X\\
&+(h(X_{s-}+x)-h(X_{s-}))\frac{x-k(x)}{x}\star \mu^X +  k(h(X_{s-}+x)-h(X_{s-}))\star (\mu^X-\nu^X).\notag
\end{align}
Subtracting \eqref{firstex} from \eqref{secondex}, we get
\begin{align*}
0&=  \bar B^{k,Y}-\Gamma^{k}(h)+  \Big[(h(X_{s-}+x)-h(X_{s-}))\frac{k(x)}{x}-k(h(X_{s-}+x)-h(X_{s-}))\Big]\star \mu^X\\
&
+  \Big[k(h(X_{s-}+x)-h(X_{s-}))-  (h(X_{s-}+x)-h(X_{s-}))\frac{k(x)}{x}\Big]\star (\mu^X-\nu^X)\\
&=  \bar B^{k,Y}-\Gamma^{k}(h)
+  \Big[k(h(X_{s-}+x)-h(X_{s-}))-  (h(X_{s-}+x)-h(X_{s-}))\frac{k(x)}{x}\Big]\star \nu^X,
\end{align*}
that provides 
\begin{align*}
\bar B^{k,Y}
&=  \Gamma^{k}(h)
-  \Big[k(h(X_{s-}+x)-h(X_{s-}))-  \frac{(h(X_{s-}+x)-h(X_{s-}))}{x}k(x)\Big]\star \nu^X.
\end{align*}
\item[(ii)]Assume moreover that  $h \in C^2$, and that  $X$ is a semimartingale  with characteristics $(B^k, C, \nu)$.  Then  it is possible to express the characteristic $\bar B^k$ of $Y$ explicitly in terms of the  characteristics of $X$. In particular, this involves a Lebesgue integral with respect to the bounded variation process $B^k$.
As a matter of fact, for every $f\in C^2 \cap C_b^0$, $f(Y)$ is a special semimartingale. By Theorem \ref{T: equiv_mtgpb_semimart} applied to $(f\circ h)(X)$,  the predictable bounded  variation part of $f(Y)$ is given by 
\begin{align}\label{firstexsem}
&(f\circ h)(X_0) + \frac{1}{2} \int_0^{\cdot}  (f\circ h)''(X_s) \,dC^X_s+  \int_0^{\cdot}  (f\circ h)'(X_s) \,d B_s^{k,X}\notag\\
&+  \int_{]0,\cdot]\times \R} ((f\circ h)(X_{s-} + x) -(f\circ h)(X_{s-})-k(x)\,(f\circ h)'(X_{s-}))\,\nu^X(ds\,dx)\notag\\
&=f(Y_0) + \frac{1}{2} \int_0^{\cdot}  [f''(h(X_{s}))(h'(X_s))^2+ f'(h(X_{s}))h''(X_s)] \,dC^X_s+  \int_0^{\cdot} f'(h(X_{s})) h'(X_s)\,d B_s^{k,X}\notag\\
&+  \int_{]0,\cdot]\times \R} [(f\circ h)(X_{s-} + x) -(f\circ h)(X_{s-})-k(x)\,f'(h(X_{s-}))h'(X_{s-})]\,\nu^X(ds\,dx).
\end{align} 
On the other hand, again by Theorem \ref{T: equiv_mtgpb_semimart} applied to $f(Y)$,  and recalling \eqref{Cbar} and \eqref{nubar}, the  process above is equal to 
\begin{align}\label{secondexsem}
&
f(Y_0) + \frac{1}{2} \int_0^{\cdot}  f''(Y_s) \,dC^Y_s+ \int_0^{\cdot}  f'(Y_s) \,d \bar B_s^{k,Y}\notag\\
&+ \int_{]0,\cdot]\times \R} (f(Y_{s-} + y) -f(Y_{s-})-k(y)\,f'(Y_{s-}))\,\nu^Y(ds\,dy)\notag\\
&= f(Y_0) + \frac{1}{2} \int_0^{\cdot}  f''(h(X_{s})) \,(h'(X_{s}))^2\,dC^X_s+ \int_0^{\cdot}  f'(h(X_s)) \,d\bar B_s^{k,Y}\notag\\
&+ \int_{]0,\cdot]\times \R} [k(x)\,h'(X_{s-})-k(h(X_{s-}+x)-h(X_{s-}))]\,f'(h(X_{s-}))\,\nu^X(ds\,dy)\notag\\
&+ \int_{]0,\cdot]\times \R} [(f \circ h)(X_{s-}+x) -(f\circ h)(X_{s-})-k(x)\,f'(h(X_{s-}))h'(X_{s-})]\,\nu^X(ds\,dx).
\end{align}
Subtracting \eqref{secondexsem} from  \eqref{firstexsem}
 we get
\begin{align}\label{charactLim}
& \frac{1}{2} \int_0^{\cdot}   f'(h(X_{s}))h''(X_s) \,dC^X_s
+  \int_0^{\cdot} f'(h(X_{s})) [h'(X_s)d B_s^{k,X}-d \bar B_s^{k,Y}]\notag\\
&- \int_{]0,\cdot]\times \R} [k(x)\,h'(X_{s-})-k(h(X_{s-}+x)-h(X_{s-}))]\,f'(h(X_{s-}))\,\nu^X(ds\,dy)=0.
\end{align}
Define now the unit partition $\chi: \R \rightarrow \R$ as the smooth function $\chi(a)$ equal to $1$ if $a \leq -1$, equal to $0$ if $a \geq 0$, 
and  such that  $\chi(a) \in [0,\,1]$ for $a \in (-1, 0)$. 
	Set 
	\begin{align}
	\chi_N(x):=\chi(|x|- (N+1))=
	\left\{
	\begin{array}{ll}
	1 \quad \textup{if}\,\, |x| \leq N\\
	0 \quad \textup{if}\,\, |x| \geq N+1\\
	\in [0,\,1]\quad 
	\textup{otherwise}.
	\end{array}
	\right.\label{chiN}
	\end{align}
	 Notice that 
	$\chi_N(x)$  is a smooth function.
We apply 	\eqref{charactLim} with $f = f_N$, where $f_N(0)=0$ and $f'_N(x) = \chi_N(x)$. 
Letting $N \rightarrow \infty$ in \eqref{charactLim}, we get 
	\begin{align*}
\bar B^{k,Y}&= \frac{1}{2} \int_0^{\cdot}  h''(X_s) \,dC^X_s+  \int_0^{\cdot}  h'(X_s)\,d B_s^{k,X}\notag\\
&- \int_{]0,\cdot]\times \R} [k(x)\,h'(X_{s-})-k(h(X_{s-}+x)-h(X_{s-}))]\,\nu^X(ds\,dx).
\end{align*}
\end{itemize}
\end{remark}

\section{Generalized   martingale problems}
\subsection{Stochastic calculus related to martingale problems} \label{S:4.1}
Let $(\Omega, \F)$ be a measurable space   and $\P$ a probability measure. 
Suppose that $X$ is  a weakly finite quadratic variation process, with canonical filtration $\mathbb F^X=(\mathcal F^X_t)$, and that, for every $v$ belonging to some linear dense subspace $\mathcal D^{\mathcal S}$  of $C^{0,1}([0,T] \times \R)$, $v(\cdot, X)$ is a weak  Dirichlet process.  
The theorem below provides  necessary and sufficient conditions under which $v(\cdot, X)$ is weak Dirichlet for every $v\in C^{0,1}([0,T] \times \R)$.

In the sequel we will make use of the following property for the couple $(\mathcal D^S, X)$. 
\begin{hypothesis}\label{newH}
For every $v \in \mathcal D^S$, $Y^v:=v(\cdot, X)$ is a  weak Dirichlet process, with unique continuous local martingale component $Y^{v,c}$. 
\end{hypothesis}
 \begin{remark}\label{R:4.4}
Let $\mathcal D^S\subseteq C^{0,1}([0,T]\times \R)$, and $X$ be a c\`adl\`ag process satisfying \eqref{int_small_jumps}. If Hypothesis \ref{newH} holds true for $(\mathcal D^S, X)$, then $\mathcal D^S$ is contained in the set $\mathcal D$ introduced in Remark \ref{R:linmapGamma}. As a matter of fact, since  $\mathcal D^S$   is  contained in  $C^{0,1}$, then Hypothesis \ref{H:3.7-3.8}  holds true, see Proposition \ref{R:th_chainrule}.  

\end{remark}

\begin{theorem}\label{T:new4.4}
Let $\mathcal D^{\mathcal S}$ be a dense subspace of $C^{0,1}([0,T] \times \R)$, and  $X$ be a weakly finite quadratic variation process.  
The following are equivalent. 
\begin{enumerate}
\item $X$ is a weak Dirichlet process.
\item 
\begin{itemize}
	\item [(i)]$v \mapsto Y^{v,c}$, $ \mathcal {D}^{\mathcal S}\rightarrow \mathbb D^{ucp}$, 	is continuous in zero.
	\item [(ii)] Hypothesis \ref{newH} holds true for $(\mathcal D^S, X)$. 
\end{itemize}
	\item $v(t, X_t)$ is a weak Dirichlet process for every $v \in C^{0,1}([0,T] \times \R)$.
\end{enumerate}
\end{theorem}
\proof
3. $\Rightarrow$ 1. This follows taking $v \equiv \textup{Id}$. 

\noindent 1. $\Rightarrow$ 2. By  Theorem  \ref{P:3.10}, for every $v \in C^{0,1}([0,T]\times \R)$, $v(\cdot, X)$ is a weak Dirichlet process, and  
$$
	Y^{v,c} = Y_0 + \int_{0}^{\cdot}\partial_x v(s, X_{s})\,dX^c_s.
$$
By Problem 2.27, Chapter 3,  in \cite{ks}, this implies the continuity stated in item (i). Moreover, item (ii) trivially holds. 

\medskip 

\noindent 2. $\Rightarrow$ 3. 
By item (ii), for every $v \in \mathcal D^S$, $v(\cdot, X)$ is a weak Dirichlet process, with unique continuous martingale component $Y^{v,c}$. Since by item (i) $v \mapsto Y^{v,c}$, $ \mathcal {D}^{\mathcal S}\rightarrow \mathbb D^{ucp}$, 	is continuous, it extends continuously to $C^{0,1}$. 
We will denote in the same way the extended operator. Since the space of continuous local martingales is closed under the u.c.p. convergence, $Y^{v,c}$ is a continuous local martingale for every $v \in C^{0,1}([0,T]\times \R)$.

We denote $A^v:=v(\cdot, X)-Y^{v,c}$, for every $v \in C^{0,1}([0,T]\times \R)$.  It remains to prove that  $A^v$ is a martingale orthogonal process, namely that, for every continuous local martingale $N$,   
\begin{equation}\label{cov_Y-Yc,N}
	[v(\cdot, X_\cdot), N)]=[Y^{v,c}, N].
\end{equation}

\medskip 

\noindent \emph{Step a)}.
Equality  \eqref{cov_Y-Yc,N} holds true for every $v \in \mathcal D^{\mathcal S}$, since $v(\cdot, X_\cdot)$ is a weak Dirichlet process, and therefore $v(\cdot, X_\cdot) - Y^{v,c}$ is a martingale orthogonal process, see Proposition \ref{P:uniqdec}.

\medskip

\noindent \emph{Step b)}. Let $(\varepsilon_n)$ be a sequence converging to zero.  We need to show that (recall Proposition A.3 in \cite{BandiniRusso1})
\begin{equation}\label{Cov_approx}
	\int_0^\cdot [v(s + \varepsilon_n, X_{s + \varepsilon_n})-v(s, X_s)] (N_{s + \varepsilon_n}- N_s) \frac{ds}{\varepsilon_n} \rightarrow [Y^{v,c}, N] \quad\textup{u.c.p. as }\,n \rightarrow \infty. 
\end{equation}
Indeed, for this, it is enough to show the existence of a subsequence, still denoted by $(\varepsilon_n)$, such that \eqref{Cov_approx} holds.
Since $N$ is a martingale, $[N, N]$ exists, so that (again by Proposition A.3 in \cite{BandiniRusso1})
$$
[N,N]_{\varepsilon_n}:= 
	\int_0^\cdot  (N_{s + \varepsilon_n}- N_s)^2 \frac{ds}{\varepsilon_n} \rightarrow [N, N] \quad\textup{u.c.p. as }\,n \rightarrow \infty. 
$$
 By extraction of subsequence, we can assume that previous convergence holds uniformly almost surely.

Let us then prove \eqref{Cov_approx}. To this end, we introduce the maps 
\begin{align*}
	T_n : \,C^{0,1}([0,T] \times \R) &\rightarrow \mathbb D^{ucp}\\
	v&\mapsto \int_0^\cdot [v(s + \varepsilon_n, X_{s + \varepsilon_n})-v(s, X_s)](N_{s + \varepsilon_n}- N_s)\frac{ds}{\varepsilon_n}.
\end{align*}
They are linear and continuous, $\mathbb D^{ucp}$ is an $F$-space in the sense of \cite{ds}, Chapter 2.1.

Suppose  that we exhibit a metric $d_{ucp}$ related to $\mathbb D^{ucp}$ such that,  
\begin{equation}\label{FinalToprove}
\textup{for every fixed}  \,\,v \in C^{0,1}([0,T] \times \R),\quad 	T_n(v) \,\,\, \textup{is bounded in }\,\mathbb D^{ucp}.
\end{equation} 
By Step a), for every $v \in \mathcal D^{\mathcal S}$ (dense subset of $C^{0,1}([0,T] \times \R)$) we already know that
\begin{equation}\label{Tn_conv}
	T_n(v) \rightarrow [Y^{v,c}, N] \quad  \textup{u.c.p. as } \,n \rightarrow \infty. 
\end{equation}
 Then Banach-Steinhaus theorem would imply the existence of a linear and continuous map $T : \,C^{0,1}([0,T] \times \R) \rightarrow \mathbb D^{ucp}$ such that
$T_n (v) \rightarrow T(v)$ u.c.p. for all  $v \in C^{0,1}([0,T] \times \R)$. 
The map $v \mapsto [Y^{v,},N]$  is continuous by Proposition \ref{P:App}. Then by \eqref{Tn_conv},  $T(v) \equiv [Y^{v,c}, N]$.
 
 \medskip
 
\noindent \emph{Step c)}. It remains  to show \eqref{FinalToprove}. Let $v \in C^{0,1}([0,T] \times \R)$. We have 
$T_n(v) = T_n^1(v) + T_n^2(v)$,
where 
\begin{align*}
	T_n^1(v)&:= \frac{1}{\varepsilon_n}\int_0^\cdot [v(s + \varepsilon_n, X_{s + \varepsilon_n})-v(s + \varepsilon_n, X_s)](N_{s + \varepsilon_n}- N_s)ds,\\
	T_n^2(v)&:=\frac{1}{\varepsilon_n}\int_0^\cdot [v(s + \varepsilon_n, X_{s})-v(s, X_s)](N_{s + \varepsilon_n}- N_s)ds.
\end{align*}
By similar  arguments as in the proof of Proposition 3.10 in \cite{gr}, 
$
T_n^2(v)\rightarrow 0$ u.c.p. as $n \rightarrow \infty$.
Indeed, since $v$ is continuous, 
\begin{align*}
	\frac{1}{\varepsilon_n^2}\int_0^T  \Big|\int_{r-\varepsilon_n}^{r} [v(s + \varepsilon_n, X_{s})-v(s, X_s)]ds\Big|^2 d[N, N]_r
\end{align*}
converges to zero, so that 
\begin{align}\label{KSconv}
	\int_0^\cdot \frac{dN_r}{\varepsilon_n} \int_{r-\varepsilon_n}^{r} [v(s + \varepsilon_n, X_{s})-v(s, X_s)]ds
\end{align}
converges u.c.p.  to zero by Problem  2.27, Chapter 3, in \cite{ks}. By stochastic Fubini theorem, we observe that the processes in \eqref{KSconv} are equal to $(T_n^2(v))$ up to a sequence of processes converging u.c.p.  to zero. 
Concerning $(T_n^1(v))$, we have 
\begin{align*}
	T_n^1(v)(t)&:= \frac{1}{\varepsilon_n}\int_0^t ds \int_0^1 da\, \partial_x v(s + \varepsilon_n, X_{s}+ a(X_{s + \varepsilon_n}-X_{s}))(X_{s + \varepsilon_n}- X_s)  (N_{s + \varepsilon_n}- N_s).
\end{align*}

\medskip

\noindent \emph{Step d)}. 
We choose 
$
d_{ucp}(X^1, X^2)=\E\Big[\sup_{t \leq T}|X_t^1 - X_t^2|\wedge 1\Big]
$
if $X^1, X^2 \in \mathbb D^{ucp}$. 
Since $T^2_n(v)$ is a converging sequence, then it is necessarily bounded. To prove that $T_n(v)$ is bounded it remains to prove the same for $T^1_n(v)$. To this end, 
let $\kappa >0$. We need to show the existence of $\delta$ such that
\begin{equation}\label{est_distance}
	d_{ucp}(\delta T_n^1(v), 0)< \kappa \quad \forall n. 
\end{equation}
Let $\mathbb K^X(\omega)$ be the random set introduced in \eqref{KX}. 
Now, introducing the (finite) random variables
\begin{align*}
 \tilde \Lambda(\omega) &:= \sup_{s \in [0,\,T], \, x \in 
	\mathbb K^X(\omega)}|\partial_x v|(s,x),\\
 \Lambda(\omega) &:=\tilde\Lambda(\omega) \, \sup_n \Big[\int_0^T(N_{s + \varepsilon_n}(\omega)- N_s(\omega))^2 ds\Big]^{1/2}, 
\end{align*}
we get 
\begin{align*}
	\sup_{t \leq T} |T_n^1(v)|&\leq  \Big[\frac{1}{\varepsilon_n}\int_0^T(X_{s + \varepsilon_n}- X_s)^2 ds\,  \frac{1}{\varepsilon_n}\int_0^T(N_{s + \varepsilon_n}- N_s)^2 ds\Big]^{1/2}\tilde\Lambda\\
	&\leq \Big[\frac{1}{\varepsilon_n^2}\int_0^T(X_{s + \varepsilon_n}- X_s)^2 ds\Big]^{1/2}  \Lambda.
\end{align*}
Since $X$ is of weakly finite quadratic variation, we can introduce  $M >0$ such that 
\begin{align*}
\P\Big(\frac{1}{\varepsilon_n^2}\int_0^T(X_{s + \varepsilon_n}- X_s)^2 ds >M^2\Big)\leq \frac{\kappa}{4}, \qquad 
\P(|\Lambda|>M) \leq \frac{\kappa}{4}, \qquad \forall n.
\end{align*}
Now, setting 
$$
\Omega_{M,n}:= \Big(\Big\{\frac{1}{\varepsilon_n^2}\int_0^T(X_{s + \varepsilon_n}- X_s)^2 ds >M^2\Big\} \cup \{|\Lambda|>M\}\Big)^c,
$$
we notice that $\P(\Omega_{M,n}^c) \leq \frac{\kappa}{2}$. 
We have 
\begin{align*}
	\E\Big[\sup_{t \leq T} \delta |T_n^1(v)| \wedge 1\Big]= \E\Big[1_{\Omega_{M,n}^c}\sup_{t \leq T} \delta |T_n^1(v)| \wedge 1\Big]+\E\Big[1_{\Omega_{M,n}}\sup_{t \leq T} \delta |T_n^1(v)| \wedge 1\Big]
\end{align*}
is bounded by 
\begin{align*}
	\frac{\kappa}{2} + \E\Big[1_{\Omega_{M,n}}\sup_{t \leq T} \delta |T_n^1(v)| \wedge 1\Big]\leq \frac{\kappa}{2}  + \E[\delta M^2 \wedge 1] \leq \frac{\kappa}{2} + \delta M^2.
\end{align*}
Formula \eqref{est_distance} follows by choosing $\delta$ so that $\delta M^2 < \frac{\kappa}{2}$.
\endproof
\begin{remark}
In Section  \ref{S:mtgpb} we will introduce a suitable notion of  
path-dependent martingale problem with respect to some operator $\mathcal A$ and domain $\mathcal D_{\mathcal A}$. %
In particular, if $X$ is a solution to the aforementioned martingale problem, $v(t, X_t)$ is a special semimartingale for every $v \in \mathcal D_{\mathcal A}$. 
In that case, if $\mathcal D_{\mathcal A}$ is a dense subspace in $C^{0,1}$, Theorem \ref{T:new4.4} (with $\mathcal D^{S}=\mathcal D_{\mathcal A}$)  will contribute to obtain a (weak Dirichlet) decomposition for $v(t, X_t)$ to every  $v \in C^{0,1}$.
In many irregular situations,  the  identity function does not belong to $\mathcal D_{\mathcal A}$ but only to $C^{0,1}$; this   allows among others  to get a  sort of Doob-Meyer type decomposition for the process $X$ itself. 
\end{remark}

The result below reformulates point 2.(i) of Theorem \ref{T:new4.4} in terms of the map $\Gamma^k$ in Theorem \ref{T:2.11}.
\begin{proposition}\label{T:3.30_bis}
Let $\mathcal D^{\mathcal S}$ be a dense subspace of $C^{0,1}([0,T] \times \R)$.  Let $X$ be a weakly finite quadratic variation process. 
Assume that Hypothesis \ref{newH} holds true for $(\mathcal D^S, X)$. 

Then, 
$v \mapsto Y^{v,c}$, $ \mathcal {D}^{\mathcal S}\subseteq C^{0,1}([0,T] \times \R) \rightarrow \mathbb D^{ucp}$, 	is continuous in zero 
 if and only if the map $\Gamma^k$ in Theorem \ref{T:2.11} restricted to $\mathcal D^{\mathcal S}$  is continuous in zero with respect to the $C^{0,1}$-topology. 
  	\end{proposition}

\proof 
The equivalence property follows by Lemma \ref{L:contmtg_d} and formula \eqref{dec_wD1} in Theorem \ref{T:2.11}.
\endproof

We end this section relating the brackets of some martingales associated with a semimartingale of the type  $v(\cdot, X)$, where $X$ is a c\`adl\`ag process, to the map $\Gamma(v)$  in Theorem \ref{T:2.11bis}. 
\begin{proposition}\label{P:3.34}
Let 
 $X$ be a c\`adl\`ag process satisfying \eqref{int_small_jumps}. 
Let $v \in C^{0,1}([0,T]\times \R)$ and bounded  such that $Y:=v(\cdot, X)$ is a  semimartingale, with unique local martingale component $Y^{c}$.  Then, 
 $Y$ is a special semimartingale with unique martingale part $N$ satisfying 
	\begin{align}\label{bracketN_bis}
\langle 
N, 
N\rangle 
&= \Gamma(v^2)-2 \int_0^\cdot  v(s, X_{s-})d \Gamma_s(v) 	- \sum_{0 \leq s \leq \cdot} \left|\Delta \Gamma_s(v)
 	\right|^2. 
\end{align}
Moreover,
\begin{align}
\langle Y^{c}, Y^{c}\rangle  &=\Gamma(v^2) -2\int_0^\cdot  v(s, X_{s-})   d \Gamma_s(v)
-[v(s,X_{s-}
 	+ x) -v(s,X_{s-})]^2 
 	\star\nu^{X, \P}\label{bracketYc}\\
	&= \Gamma^c(v^2)-2 \int_0^\cdot  v(s, X_{s-}) d \Gamma_s^c(v)  - [v(s, X_{s-}+x)-v(s, X_{s-})]^2 \star \nu^{X,\P,c},  \label{Ycn_sigmaA_bis}
\end{align}
where 
$\nu^{X, \P,c}$ (resp. $\Gamma^c(v)$) is the continuous part of $\nu^{X, \P}$
(resp.  of  $\Gamma(v)$), 
where $\Gamma(v)$ is the process  appearing in Theorem \ref{T:2.11bis}. 
\end{proposition}
\begin{remark}
	Formula \eqref{bracketN_bis} implies in particular that if $\Gamma(v)$ and $\Gamma(v^2)$ are continuous, then $\langle  N,  N\rangle$ is continuous as well. 
\end{remark}
\begin{remark}
By Theorem \ref{T: equiv_mtgpb_semimart}, $Y^{v^2}$ is a semimartingale. On the other hand, $Y^{v}:=v(\cdot, X)$  and $Y^{v^2}:=v^2(\cdot, X)$ are special weak Dirichlet processes by Theorem \ref{T:spw} and Remark \ref{R:boundspecfirst}, being $v,v^2\in C^{0,1}$ and bounded. This implies that $Y^{v}$  and $Y^{v^2}$ are special semimartingales. 
\end{remark}
\noindent \emph{Proof of Proposition \ref{P:3.34}}.
We first  prove identity \eqref{Ycn_sigmaA_bis}. 
By Theorem \ref{T:2.11bis},  we get that
	\begin{align}
 	N_t &:= 
 	Y_t -Y_0 - \Gamma_t(v),\label{N1tA}\\
 	\bar N_t &:= Y_t^2 - Y_0^2
 	-  \Gamma_t(v^2)\label{N2tA}
 	\end{align}
 	are local martingales under $\P$. 
 	 Applying the integration by parts formula, and  taking into account  	 \eqref{N1tA} and \eqref{N2tA}, we get 
 	   	\begin{align}\label{bracketYA}
[Y, Y]_t &=Y_t^2- Y_0^2 -2 \int_0^t Y_{s-} d Y_s\notag\\
&=  \Gamma_t(v^2) -2\int_0^t  v(s, X_{s-})   d \Gamma_s(v)  + \bar N_t - 2 \int_0^t v(s,X_{s-}) d N_s.	
 	\end{align}
 	We now show that 
  \begin{align}
&[v(s,X_{s-}
 	+ x) -v(s,X_{s-})]^2 
 	\star\nu^{X, \P}
 	\in \mathcal{A}_{\textup{loc}}^+.\label{int_v_Aloc}
  \end{align}
  As a matter of fact, 
    \begin{align}
&[v(s,X_{s-}
 	+ x) -v(s,X_{s-})]^2 
 	\star\nu^{X, \P}\label{squarenu}
 	\\
 	&= [v(s,X_{s-}
 	+ x) -v(s,X_{s-})]^2 \frac{k^2(x)}{x^2}
 	\star\nu^{X, \P}+[v(s,X_{s-}
 	+ x) -v(s,X_{s-})]^2\frac{x^2 -k^2(x)}{x^2} 
 	\star\nu^{X, \P}.\notag
  \end{align}
The first term in the right-hand side of \eqref{squarenu} belongs to $\mathcal{A}_{\textup{loc}}^+$ by Proposition \ref{P:2.10}.
On the other hand, let $c>0$ such that $k(x)=x$ on $[-c,c]$. Then the second term in the right-hand side of \eqref{squarenu} belongs to $\mathcal{A}_{\textup{loc}}^+$ by 
Lemma \ref{L:c} and the fact that $v$ is bounded.

 At this point, 	by \eqref{QVC} and Proposition \ref{P:2.14},
\begin{align}\label{brackY}
&[Y, Y]= \langle Y^{c}, Y^{c}\rangle+    \sum_{s \leq \cdot} |\Delta Y_s|^2=\langle Y^{c}, Y^{c}\rangle +  [v(s,X_{s-}+x)-v(s,X_{s-})]^2 \star \mu^X\\
&=\langle Y^{c}, Y^{c}\rangle +  [v(s,X_{s-}+x)-v(s,X_{s-})]^2 \star \nu^{X, \P}
+  [v(s,X_{s-}+x)-v(s,X_{s-})]^2 \star (\mu^X-\nu^{X, \P}).\notag
\end{align}
Since   \eqref{bracketYA} and \eqref{brackY} provide two decompositions of the same special semimartingale $[Y,Y]$,  we get \eqref{bracketYc}.
Denoting by $C$ the right hand-side of \eqref{bracketYc}, we remark that  $C$ is a continuous bounded variation process. 
Therefore, taking into account Remark \ref{R:3items}-1., \eqref{bracketYc} implies  \eqref{Ycn_sigmaA_bis}.

On the other hand, 
Theorem \ref{T:2.11bis}  implies the unique decomposition \eqref{N1tA}
with
 \begin{align}
 	 	N_t
&:=  Y^{c}_t -Y_0  + [v(s,X_{s-} + x) -v(s,X_{s-})] 
 	\star(\mu^X-\nu^{X, \P})_t.\label{barNtA}
 \end{align}
Now we notice that, by formula \eqref{int_v_Aloc},  
 the stochastic integral appearing in the right-hand side of \eqref{barNtA} is a locally square integrable martingale, 
and its predictable bracket yields
$$
[v(s,X_{s-}
 	+ x) -v(s,X_{s-})]^2 
 	\star\nu^{X, \P}\notag\\
 	- \sum_{0 \leq s \leq \cdot} \Big|\int_{\R}  [v(s,X_{s-}
 	+ x) -v(s,X_{s-})] 
 	\nu^{X, \P}(\{s\}\times dx)\Big|^2,
$$
see the end of Section \ref{SPrelim}. 
Since that purely discontinuous martingale  is martingale orthogonal, 
\begin{align}\label{predbracketNbarA}
\langle 
N, 
N\rangle 
&=\langle Y^{c}, Y^{c}\rangle  +    [v(s,X_{s-}
 	+ x) -v(s,X_{s-})]^2 
 	\star\nu^{X, \P}\notag\\
 	&- \sum_{0 \leq s \leq \cdot} \Big|\int_{\R}  [v(s,X_{s-}
 	+ x) -v(s,X_{s-})] 
 	\nu^{X, \P}(\{s\}\times dx)\Big|^2.
\end{align}
Identity \eqref{bracketN_bis} follows by plugging \eqref{bracketYc} in \eqref{predbracketNbarA}. 
 \qed

\subsection{Definition and main properties of the martingale problem}\label{S:mtgpb}
Given $\eta \in D_{-}(0,\,T)$ (resp. $\zeta \in D(0,\,T)$), we will use the notation 
\begin{align*}
	\eta^{t}(s) := 
	\left\{
	\begin{array}{ll}
	\eta(s) \quad \textup{if}\,\,s <t,\\
	\eta(t)\quad \textup{if}\,\,s \geq t.
	\end{array}
	\right.
\end{align*}

We denote by $C^{NA}(D_{-}(0,\,T); B(0,T))$ the subspace  of $F \in C(D_{-}(0,\,T); B(0,T))$ such that $F(\eta)(s)= F(\eta^s)(s)$ for every $\eta \in D_{-}(0,\,T)$ and $s \in [0,T]$. From now on, for simplicity, we will write $F(s,\eta):=F(\eta)(s)$.

 We consider the following hypothesis for a triplet $(\mathcal D, \Lambda, \gamma)$. 
\begin{hypothesis}\label{H:nonant}
 $\mathcal D\subseteq C^{0,1}([0,T] \times \R)$, $\Lambda: \mathcal D \rightarrow C^{NA}(D_{-}(0,\,T); B(0,T))$ is a linear map.  $ \gamma:  [0,\,T] \times D_{-}(0,\,T)\rightarrow \R$ is such that 
for every  $\eta \in  D_{-}(0,\,T)$, $\gamma(\cdot, \eta)$ is of bounded variation, and fulfills the non-anticipating property, 
 i.e., for every $\eta \in D_{-}(0,\,T)$, 
 $
\gamma(t,\eta) = \gamma(t, \eta^{t})$.
\end{hypothesis}

\begin{definition}\label{D:newA}
Fix $N \in \N$. 
Let $\mathcal D_{\mathcal A}\subseteq C^{0,1}([0,T] \times \R)$,  $\Lambda_i: \mathcal D_{\mathcal A} \rightarrow C^{NA}(D_{-}(0,\,T); B(0,T))$  and  $ \gamma_i:  [0,\,T] \times D_{-}(0,\,T) \rightarrow \R$,  such that  $(\mathcal D_{\mathcal A}, \Lambda_i, \gamma_i)$ fulfill Hypothesis \ref{H:nonant} for all $i=1,..., N$.
For every  $v \in \mathcal {D}_{\mathcal A}$, $\eta \in D_{-}(0,\,T)$, we set
\begin{equation}\label{newA}
(\mathcal{A} v)(ds,\eta):= \sum_{i=1}^N (\Lambda_i v)(s,\eta)\,\gamma_i(ds, \eta).
\end{equation}
\end{definition} 
\begin{remark}
We have the decomposition $(\mathcal{A} v)(ds,\cdot)=  (\mathcal{A} v)((ds)^c,\cdot)+ (\mathcal{A} v)(\Delta s,\cdot)$, with
\begin{align}
(\mathcal{A} v)(\Delta s,\cdot)&:= \sum_{i=1}^N (\Lambda_i v)(s,\cdot)\, \gamma_i(\Delta s, \cdot),\label{newAjump}\\
(\mathcal{A} v)((ds)^c,\cdot)&:= \sum_{i=1}^N (\Lambda_i v)(s,\cdot)\,\gamma_i^c(ds, \cdot).\label{Acont}
\end{align}
\end{remark}

\normalcolor

\begin{definition}\label{D:mtpb1}
Let $\mathcal{A}$ and $\mathcal D_{\mathcal A}$  as  in Definition \ref{D:newA}. 
A process $X$  is said to solve the martingale problem (under a probability $\P$) 
with respect  to 
$\mathcal A$, $\mathcal D_{\mathcal A}$ and $x_0 \in \R$, if
\begin{itemize}
\item [(i)] condition \eqref{int_small_jumps} holds under $\P$;
\item[(ii)] 
for any  
$v \in \mathcal {D}_{\mathcal A}$ and bounded, 
\begin{align}\label{mtg_pb_general}
M^v_t := v(t,X_{t}) - v(0,x_0) - \int_0^{t} (\mathcal{A} v)(ds,X^{-}) 
\end{align}
is an $(\mathcal F^X_t)$-local martingale. 
\end{itemize}
\end{definition}
\begin{remark}
Let $\mathcal{A}$ and $\mathcal D_{\mathcal A}$  as  in Definition \ref{D:newA}, and set
	\begin{equation}
	\mathcal X = \{
	M^v
	:\,\,v \in \mathcal D_{\mathcal A}
	\},\label{chi_ex}
	\end{equation}
	with $M^v$ 
	 defined in \eqref{mtg_pb_general}. 
	Then $(X,\P)$ solves the martingale problem with respect to $\mathcal A$, $\mathcal D_{\mathcal A}$ and  $x_0$ if and only if 
 $\P$ is solution of the martingale problem in Definition  1.3, Chapter III, in \cite{JacodBook}
  associated to 
 $(\Omega, \mathcal F, \mathbb F)$, 
 $	\mathcal X$ in \eqref{chi_ex}, $ \mathcal H=\{B \in \mathcal F: \exists B_0 \in \mathcal B(\R) \textup{ such that } B= \{\omega \in \Omega: \omega(0) \in B_0\}\}$ and $\P_H $ corresponds   to $ \delta_{x_0}$  in the sense that, for any $B \in \mathcal F$, $\P_H(B) = \delta_{x_0}(B_0)$ with $B_0=\{\omega(0) \in \R: \,\, \omega \in B\}$. 
\end{remark}

\begin{definition}[Existence]
 We say that the martingale problem related to $\mathcal A$, $\mathcal D_{\mathcal A}$ and $x_0$ meets existence if there exists  a couple  $(X,\P)$ on some $(\Omega, \mathcal F)$, solution to the martingale problem in the sense of Definition \ref{D:mtpb1}.
\end{definition}

\begin{definition}[Uniqueness]
We say that the martingale problem related to $\mathcal A$, $\mathcal D_{\mathcal A}$ and $x_0$ admits uniqueness if, 
	given two  spaces $(\Omega_1, \mathcal F_1)$ and $(\Omega_2, \mathcal F_2)$, and two solutions $(X_1, \P_1)$ and $(X_2, \P_2)$ to the martingale problem in the sense of Definition \ref{D:mtpb1}, then $\mathcal L(X_1)|_{\P_1} = \mathcal L(X_2)|_{\P_2}$.
\end{definition}

\begin{remark}\label{3.2}
Let $(\check \Omega, \check {\mathcal F})$ be the canonical space,  and $\mathbb F$ be the canonical filtration of the canonical process $\check X$.
 $(X, \P)$ is a solution the the martingale problem  related to $\mathcal A$, $\mathcal D_A$ and  $x_0$, 
   if and only if  $(\check X, \mathcal L(X))$  is a solution to the martingale problem related to  $\mathcal A$, $\mathcal D_A$ and  $x_0$.
   \end{remark}

\subsection{Time-homogeneous towards time-inhomogeneous martingale problem}\label{S:mtpb_hom}

\begin{definition}\label{D:mtpb_hom}
Let $\mathcal D_{\mathcal L}\subseteq C^0_b$ and  $\mathcal L: \mathcal D_{\mathcal L} \rightarrow C^{NA}(D_{-}(0,\,T); B(0,T))$ \textcolor{blue}{be a linear map}.
We say that $(X, \P)$ fulfills the time-homogeneous martingale problem with respect to $\mathcal {D}_{\mathcal L}$, $\mathcal L$ and $x_0$, if for any  $f \in \mathcal {D}_{\mathcal L}$ and bounded, 
the process
\begin{align}\label{mtg_pb_timehom}
M^f_t := f(X_{t}) - f(x_0) - \int_0^{t} (\mathcal{L} f)(s,X^{-})  ds
\end{align}
is an $(\mathcal F^X_t)$-local martingale  under $\P$.
 From here on we adopt this notation.
\normalcolor
\end{definition}

%

Let $C_{BUC}^{NA}(D_{-}(0,\,T); B(0,T))$ be the set of functions  $G \in C^{NA}(D_{-}(0,\,T); B(0,T))$ bounded and uniformly continuous on closed balls $B_M \subset D_{-}(0,T)$ of radius $M$.
$C_{BUC}^{NA}(D_{-}(0,\,T); B(0,T))$ is a Fr\'echet space equipped  with the distance generated by the seminorms
$$
\sup_{ \eta \in B_M} ||G(\eta)||_\infty, \quad M \in \N.
$$
\begin{theorem}\label{T:passage}
Assume   that  $D_{-}(0,\,T)$ is equipped with the metric topology of the uniform convergence on  closed balls.
Let  $\mathcal D_{\mathcal L}$ be a Fr\'echet space,  topologically included in $C^0_b$, and   equipped with some metric $d_{\mathcal L}$. Let 
 $\mathcal L: \mathcal D_{\mathcal L} \rightarrow C^{NA}_{BUC}(D_{-}(0,\,T); B(0,T))$ be a \textcolor{blue}{linear} continuous map. 
 $(X, \P)$	fulfills the time-homogeneous martingale problem in Definition \ref{D:mtpb_hom} with respect to ${\mathcal D}_{\mathcal L}$, $\mathcal L$ and $x_0$, 
if and only if $(X, \P)$	fulfills the time-inhomogeneous martingale problem in Definition \ref{D:mtpb1} with respect to $x_0$,
\begin{equation}\label{D:DAnew}
	{\mathcal D}_{\mathcal A}:=C^1([0,\,T]; {\mathcal D}_{\mathcal L}), 
\end{equation}
and 
\begin{equation}\label{Adistrib1}
(\mathcal{A} v)(dt,\eta):=    (\partial_t v(t,\eta)+ (\mathcal L v(t, \cdot))(t,\eta))d t, \quad v \in \mathcal D_{\mathcal A}, \,\,\eta \in D_{-}(0,\,T).
\end{equation}
\begin{remark}
\begin{itemize}
\item[(i)]	
A typical example of metric $d_\shl$ comes from the graph topology.
Let $\mathcal D_{\mathcal L} \subseteq C^0_b$,  and  $\mathcal L: \mathcal D_{\mathcal L} \rightarrow C^{NA}_{BUC}(D_{-}(0,\,T); B(0,T))$ be a measurable map. Assume that $\mathcal D_{\mathcal L}$ is equipped with the graph topology of $\mathcal L$:  $v_n \rightarrow 0$ in $\mathcal D_{\mathcal L}$ if and only if $v_n \rightarrow 0$ in $C_b^0$ and   $\mathcal L v_n \rightarrow 0$ in $C^{NA}_{BUC}(D_{-}(0,\,T); B(0,T))$. Then $\mathcal L$ is obviously a  continuous map.
\item [(ii)] ${\mathcal D}_{\mathcal A}$  in \eqref{D:DAnew} is constituted by bounded functions whose  derivative in time  is also bounded. 
\item [(iii)]
Since $\mathcal L$ is continuous, there exists a constant $C$ such that 
for every closed ball $B_M \subset D_{-}(0,T)$ of radius $M$, 
$$
\sup_{t \in [0,T]}\sup_{\eta \in B_M}||(\mathcal L v(t, \cdot))(t,\eta)||_{\infty} \leq C.
$$
\end{itemize}
\end{remark}

\end{theorem}

\proof
The \emph{if} implication is trivial choosing $v$ not depending on time. 

 Let us now prove the \emph{only if} implication. We need to show that, for every $v \in  {\mathcal D}_{\mathcal A}$, 
 \begin{align}\label{mtgpb_L}
 	M^v_t :=v(t, X_t) -v(0,x_0)- \int_0^t \left(\partial_s v(s, X_s) + (\mathcal{L} v(s, \cdot))(s,X^{-})\right)ds
 \end{align}
 is a local martingale. 
If  $v = f \in \mathcal D_{\mathcal L}$ then \eqref{mtgpb_L} is obviously a  local martingale. Suppose that $v(t,x)= a(t) f(x)$, with $a \in C^1(0,\,T)$, and  $f \in \mathcal D_{\mathcal L}$. By integrating by parts, we get 
\begin{align}\label{intbyparts}
	v(t, X_t)&= a(0) f(x_0) + \int_0^t a'(s) f(X_s) ds + \int_0^t a(s)df(X_s) \notag\\
	&=a(0) f(x_0) + \int_0^t  a'(s) f(X_s)ds  + \int_0^t a(s)(\mathcal L f)(s,X^{-})ds+ \int_0^t a(s)dM_s^f\notag\\
	&=v(0, x_0) + \int_0^t (\mathcal A v)(ds, X^{-}) + M^v_t, 
\end{align}
where $M^v$ is a local martingale.  
 So \eqref{mtgpb_L} is a  local martingale when $u \in \hat {\mathcal D}_{\mathcal A}$,  with $\hat {\mathcal D}_{\mathcal A}$ the set in Lemma \ref{L:hatD}.

 Let now $v \in \mathcal D_{\mathcal A}$. By Lemma \ref{L:hatD}, there is a
 sequence $v_n \in \hat {\mathcal D}_{\mathcal A}$, 
 and 
 such that $v_n \underset{n \rightarrow \infty}{\rightarrow} v$ with respect to $C^1([0,\,T]; \mathcal D_{\mathcal L})$.
By \eqref{intbyparts}  we have that 
\begin{align}\label{uneq}
	v_n(t, X_t)=v_n(0, x_0) + \int_0^t (\partial_s v_n(s, X_{s-})+ (\mathcal L v_n(s, \cdot))(s,X^{-}))ds + M^{v_n}_t, 
\end{align}
where $M^{v_n}$ is a local martingale. 
%
%
Since the maps $(t, \eta) \mapsto (\mathcal L v_n(t, \cdot))(t, \eta)$ 
 converge to $(t, \eta) \mapsto (\mathcal L v(t, \cdot))(t, \eta)$
 uniformly on compacts in $[0,T] \times D_{-}(0,T)$,  and $\partial_s v_n \underset{n \rightarrow \infty}{\rightarrow} \partial_s v$ in $C([0,T] \times \R)$, it follows that $M^{v_n}\underset{n \rightarrow \infty}{\rightarrow} M^{v}$ u.c.p.

It remains to show that $M^v$  is a local martingale.
Since $X_{s-}$ is a c\`agl\`ad process, it is locally bounded, see Remark \ref{R:pred}-1. So, for every $\ell >0$, we define 
$\tau^\ell := \inf\{t \in [0,T] :\,\,|X_{t-}| \geq \ell\},
$
with the usual convention that $\inf \emptyset = +\infty$. Clearly, $\tau^\ell \uparrow + \infty$ a.s.
Then, on $\Omega_\ell:= \{\tau^\ell \leq T\}$, 
$\sup_{s \leq T}|X_{(s \wedge \tau^\ell)^-}|\leq \ell$ a.s., 
and 
$\sup_{s \leq T}||(X^{^-})^{(s \wedge \tau^\ell)}||_\infty\leq \ell$ a.s.

 It is thus enough to show that, for every $\ell$, $(M^v)^{\tau^\ell}$ is a  martingale.
 Now, by \eqref{uneq},
 \begin{align*}
	 (M^{v_n})^{\tau^\ell}_t = v_n(t\wedge \tau^\ell, X_{t\wedge \tau^\ell})-v_n(0, x_0) - \int_0^{t \wedge \tau^\ell}(\partial_s v_n(s, X_{s-})+(\mathcal L v_n(s, \cdot))(s,X^{-})ds, 
\end{align*}
so 
\begin{align}\label{numR}
	\sup_{t \leq T}|(M^{v_n})^{\tau^\ell}_t|&\leq 2\sup_{\underset{n}{t \in [0,T],\,|x|\leq \ell}}|v_n(t,x)|+ T\sup_{\underset{n}{t \in [0,T],\,|x|\leq \ell}}|(\partial_t v_n)(t,x)|\notag\\
	&+ T \sup_{\underset{n}{t \in [0,T],\,||\eta||_\infty\leq \ell}} |(\mathcal L v_n(t, \cdot))(t,\eta)|=:R,
\end{align}
where we recall that  $C^{NA}_{BUC}(D_{-}(0,\,T); B(0,T))$ is equipped with the metric topology of the uniform convergence on closed balls.
On the other hand, recalling that $M^{v_n}\rightarrow M^{v}$ u.c.p. and therefore  $(M^{v_n})^{\tau^\ell}\rightarrow (M^{v})^{\tau^\ell}$ u.c.p., by \eqref{numR} we have in particular that 
$\sup_{t \leq T}|(M^{v})^{\tau^\ell}_t|\leq R$.
 It follows that $(M^v)^{\tau^\ell}$ is a  martingale.
This concludes the proof.
\qed

The following  is a  consequence of Theorem \ref{T:passage} in the Markovian context. 
\begin{theorem}\label{R:passage}
Let $\mathcal D_{\mathcal L_M}$ be a Fr\'echet space,  topologically included in $C^0_b$, and  equipped with some metric $d_{\mathcal L_M}$. Let   
 $\mathcal L_M: \mathcal D_{\mathcal L_M} \rightarrow C^0$ be a continuous map. 
Set $\mathcal D_{\mathcal L}:=\mathcal D_{\mathcal L_M}$  equipped with the same topology of  $\mathcal D_{\mathcal L_M}$, 
\begin{equation}\label{D:DAnewMark}
	{\mathcal D}_{\mathcal A}:=C^1([0,\,T]; {\mathcal D}_{\mathcal L_M})
\end{equation}	
 and, for every $f \in \mathcal D_{\mathcal L_M}$, $v \in {\mathcal D}_{\mathcal A}$, $\eta \in D_{-}(0,\,T)$, 
\begin{align}
	(\mathcal  L f) (t,\eta) := (\mathcal L_M f)( \eta_t), \quad 
	(\mathcal  A v) (dt,\eta) := 
	(\partial_t v(t,\eta_t) +(\mathcal L_M v(t,\cdot))(\eta_t))d t.\label{4.16}
\end{align}
Then the following are equivalent.
\begin{itemize}
\item [(i)]  
 $(X, \P)$	fulfills the time-homogeneous martingale problem in Definition \ref{D:mtpb_hom} with respect to ${\mathcal D}_{\mathcal L}$, $\mathcal L$ and $x_0$. 
 \item [(ii)] For any $f \in \mathcal D_{{\mathcal L}_M}$, 
\begin{align}\label{mtg_pb}
M^f &:= f(X_{\cdot}) - f(x_0) - \int_0^{\cdot} {\mathcal L}_M f(X_{s}) ds
\end{align}
is an $(\mathcal F_t^X)$-local martingale  under $\P$.
\item[(iii)]  
$(X, \P)$	fulfills the time-inhomogeneous martingale problem in Definition \ref{D:mtpb1} with respect to $\mathcal D_{\mathcal A}$, $\mathcal A$ and $x_0$.
\end{itemize}	
\end{theorem}
\proof
By 
\cite{BandiniRusso_DistrDrift}, $\mathcal L: \mathcal D_{\mathcal L} \rightarrow C^{NA}_{BUC}(D_{-}(0,\,T); B(0,T))$ is continuous, where $C^{NA}_{BUC}(D_{-}(0,\,T); B(0,T))$  is equipped with the topology of uniform convergence on closed balls.  
Items (i) and (ii) are equivalent by construction. 
Items (i) and (iii) equivalent by Theorem \ref{T:passage}.  
\endproof

\subsection{Weak Dirichlet characterization of the solutions to the martingale problem} 
In the sequel we prove that  stochastic calculus framework in Section \ref{S:4.1}  fits in the concrete examples of martingale problems. 

Suppose that $(X,\P)$ is  a solution to the martingale problem in Definition \ref{D:mtpb1} related to $\mathcal A$, $\mathcal {D}_{\mathcal A}$, and $x_0$. 
First of all, we set $\mathcal D^S= \mathcal D_{\mathcal A}$. By definition of the martingale problem, for every $v \in \mathcal D^{ S}$ and bounded, $v(\cdot, X)$ is a special semimartingale. In addition, by Proposition \ref{R:th_chainrule},  Hypothesis \ref{H:3.7-3.8}  holds true being $v\in C^{0,1}([0,T] \times \R)$. We can therefore apply   Theorem \ref{T:2.11bis}, that provides decomposition \eqref{Spec_weak_chainrule1} and  $\Gamma(v)$ in \eqref{GammaY}. By the uniqueness of the decomposition of the special semimartingale, we get that 
\begin{equation}\label{Gammav}
	\Gamma(v) = \int_0^\cdot (\mathcal A v)(ds, X^{-}) 
\end{equation}
and 
\begin{equation}\label{Gammavk}
\Gamma^k(v) =\int_0^\cdot (\mathcal A v)(ds, X^{-})  - (v(s, X_{s-}+x)-v(s, X_{s-})) \frac{x - k(x)}{x}\star \nu^{X, \P}.
\end{equation}
Taking into account Remark \ref{R:3items}-1., we get in particular that, for every $v \in \mathcal D_{\mathcal A}$, 
\begin{equation}\label{rel_jumps}
(\mathcal{A} v)(\Delta t,X^{-}) = \int_{\R}(v(t, X_{t-}+x)-v(t, X_{t-})) \, \nu^{X, \P}(\{t\} \times dx), \quad t \in [0,\,T]. 
 \end{equation}

\begin{corollary}\label{C:3.30bis}
Let $\P$ be a probability on $(\Omega, \mathcal F)$, and 	 $X$ be a c\`adl\`ag process with weakly  finite quadratic variation. 
	Let  $ \mathcal {D}_{\mathcal A}$ be a dense subset of $C^{0,1}([0,T] \times \R)$, and 
$\mathcal{A} $ be the operator   in Definition \ref{D:newA}. 
For every $v \in \mathcal {D}_{\mathcal A}$, set  $Y^v_t=v(t, X_t)$, and  denote by $Y^{v,c}$ its unique continuous martingale component. 	If  $(X, \P)$ is a solution to the martingale problem in Definition \ref{D:mtpb1} related to $\mathcal D_{\mathcal A}$, $\mathcal A$ and $x_0$, 
then the following are equivalent. 
\begin{itemize}
	\item [(i)] 
$v \mapsto Y^{v,c}$, $ \mathcal {D}_{\mathcal A}\subseteq C^{0,1}([0,T] \times \R)\rightarrow \mathbb D^{ucp}$, 	is continuous in zero. 
\item[(ii)] The map 
 $$
v \mapsto  \int_0^\cdot (\mathcal A v)(ds, X^{-})  - (v(s, X_{s-}+x)-v(s, X_{s-})) \frac{x - k(x)}{x}\star \nu^{X, \P}
 $$
  is continuous in zero with respect to the $C^{0,1}$-topology.
 \item [(iii)] $X$ is a weak Dirichlet process.  
  \end{itemize} 
	\end{corollary}
\proof
The equivalence of items (i) and (iii) comes from Theorem  \ref{T:new4.4}.
The equivalence of items (i) and (ii) follows from Proposition \ref{T:3.30_bis} together with \eqref{Gammavk}. 
\endproof

\begin{corollary}\label{C:3.34_bis}
Let  $ \mathcal {D}_{\mathcal A}$ be a dense subset of $C^{0,1}([0,T] \times \R)$, and 
$\mathcal{A} $ be the operator   in Definition \ref{D:newA}. 
	Let $(X, \P)$ is a solution to the martingale problem in Definition \ref{D:mtpb1}  with respect to  $\mathcal D_{\mathcal A}$,  $\mathcal A$, and  $x_0 \in \R$. 
Then, for every  $v \in \mathcal {D}_{\mathcal A}$ and bounded,  the unique martingale part $N$ of the special semimartingale $Y_t:=v(t, X_t)$ is given by 	\begin{align}\label{bracketN}
\langle 
N, 
N\rangle 
&= \int_0^\cdot (\mathcal A v^2)(ds,X^{-})-2 \int_0^\cdot  v(s, X_{s-}) (\mathcal A v)(ds,X^{-})
 	- \sum_{0 \leq s \leq \cdot} \left|(\mathcal{A} v)(\Delta s,X^{-})
 	\right|^2. 
\end{align}
Moreover,
\begin{align}\label{Ycn_sigmaA}
	\langle Y^{v,c}, Y^{v,c} \rangle &= \int_0^\cdot (\mathcal A v^2)(ds,X^{-})-2 \int_0^\cdot  v(s, X_{s-}) (\mathcal A v)(ds,X^{-})  - [v(s, X_{s-}+x)-v(s, X_{s-})]^2 \star \nu^{X,\P},  
\end{align}
where $Y^{v,c}$ is the  unique continuous martingale part of $Y^v=v(\cdot, X)$,  and   
 $(\mathcal A v)(\Delta s,\cdot)$
is introduced in \eqref{newAjump}. 
\end{corollary}
\proof
The result follows from Proposition \ref{P:3.34} together with \eqref{Gammav}.
\endproof

\normalcolor


%
%


%
%
\subsection{Examples of martingale problems}

\label{SSExamples}

\subsubsection{Semimartingales}
\label{E:4.18}
	Let $X$ be an adapted  c\`adl\`ag  real  semimartingale with characteristics $(B^k, C, \nu)$, with  decomposition $\nu(\omega, ds\,dx) = \phi_s(\omega, dx)d\chi_s(\omega)$. For further notations we refer to Section \ref{A:D}. 
$X$ verifies condition \eqref{int_small_jumps},  since  it is a semimartingale, see Remark \ref{R:equiv}.  

Then $(X, \P)$ is a solution of the martingale problem  in Definition \ref{D:mtpb1}  with respect to $\mathcal A$, $\mathcal D_{\mathcal A}$ and $x_0$, with $\mathcal D_{\mathcal A}$  the set of functions in $C^{1,2}_b$ restricted to $[0,T] \times \R$ and, for every $v \in \mathcal D_{\mathcal A}$, 
\begin{align*}
(\mathcal A v)(ds, \eta)&:= \partial_s v(s, \eta_s)ds + \frac{1}{2} \partial_{xx}^2 v(s, \eta_s) \,d(C\circ \eta)_s +  \partial_x v(s, \eta_s) \,d (B^{k}\circ \eta)_s\notag\\
&+  \int_{\R} (v(s, \eta_s + x) -v(s, \eta_s)-k(x)\,\partial_x v(s, \eta_s))\,\phi_s(\eta, dx)d\chi_s(\eta), \quad \eta \in D_{-}(0,\,T). 
\end{align*}
This follows by Remark \ref{R:D2}.


\subsubsection{Weak Dirichlet processes derived from semimartingales}
\label{E:4.17}
	Let $X$ be a weak Dirichlet process with characteristics $(B^k, C, \nu)$, with  $\nu(\omega, ds\,dx) = \phi_s(\omega, dx)d\chi_s(\omega)$. For further notations we refer to Section \ref{S:genchar}. 
	Assume that there exists $h\in C^{0,1}$, $h(t, \cdot)$ bijective and  $h(t, X_t)$ is a semimartingale with characteristics $(\bar B^k, \bar C, \bar \nu)$. Being $h^{-1} \in C^{0,1}$, we have 
	$$
	\Delta X_s = h^{-1}(s, Y_{s-} + \Delta Y_s)-h^{-1}(s, Y_{s-})= \Delta Y_s \int_0^1 \partial_x h^{-1}(s, Y_{s-} + a \Delta Y_s)da,  
	$$
	 so condition \eqref{int_small_jumps} for $Y$ (see previous subsection) implies the same for $X$.  

Then by Example \ref{E:4.18},   $(X, \P)$ is a solution of the martingale problem  in Definition \ref{D:mtpb1}  with respect to $\mathcal A$, $\mathcal D_{\mathcal A}$ and $x_0$, with $\mathcal D_{\mathcal A}$ the set of functions $v \in C^{0,1}_b$ such that $v \circ h^{-1} \in C^{1,2}$, restricted to $[0,T] \times \R$,  and, for every $v \in \mathcal D_{\mathcal A}$, 
 	\begin{align*}
&(\mathcal A v)(ds, \gamma)\\
&= \partial_s (v \circ h^{-1})(s, h(s,\gamma_s))ds + \frac{1}{2} \partial_{xx}^2 (v \circ h^{-1})(s, h(s,\gamma_s))\,(\partial_x h(s,\gamma_s))^2\,d(C \circ \gamma)_s\\
& +  \partial_x (v \circ h^{-1})(s, h(s,\gamma_s)) \,d (\bar B^{k} \circ h(\cdot, \gamma))_s\notag\\
&+(v(s, \gamma_s +x) -v (s, \gamma_s)-k(h(s, \gamma_s +x)-h(s, \gamma_s))\,\partial_x h^{-1} (s, h(s,\gamma_s))\,\partial_x v(s, \gamma_s))\, \phi_s(\gamma, dx)d\chi_s(\gamma)
\end{align*}
for every $\gamma \in  D_{-}(0,\,T)$.



\subsubsection{Discontinuous Markov processes with distributional drift}\label{S:433}
In  \cite{BandiniRusso_DistrDrift} we study existence and uniqueness for a time-homogeneous martingale problem with  distributional drift in a discontinuous Markovian framework.

In this section, we will consider a fixed $\alpha \in [0,1]$. If  $\alpha \in ]0,1[$, 
 $C^\alpha_{\textup{loc}}$ denotes  the space of functions locally $\alpha$-H\"older continuous. 
By $C^{1+\alpha}_{\textup{loc}}$ we  will denote the functions in $C^1$ whose derivative is  $\alpha$-H\"older continuous.
For convenience,  we set  $C^0_{\textup{loc}}:=C^0$,  $C^1_{\textup{loc}}:=C^1$,  $C^2_{\textup{loc}}:=C^2$.

Let $k \in \mathcal K$ be a continuous function. Let $\beta=\beta^k:  \R \rightarrow \R$ and  $\sigma:  \R \rightarrow \R$ be  continuous functions, with $\sigma$ not vanishing at zero. We consider formally the PDE operator of the type
\begin{equation}\label{Lbeta}
 L
  \psi = \frac{1}{2}\sigma^2  \psi'' + \beta' \psi' 
\end{equation}
in the sense introduced by \cite{frw1, frw2}.
 For a mollifier $\phi \in \mathcal S(\R)$ with $\int_\R \phi(x) dx =1$, we set 
 $$
\phi_n(x) := n \,\phi(n x), \quad  \beta'_n := \beta' \ast \phi_n,  \quad  \sigma_n := \sigma \ast \phi_n.
$$
  \begin{hypothesis}\label{H:Sigma}
\begin{enumerate}
\item We assume the existence of the function 
\begin{equation}\label{Sigma}
\Sigma(x) := \lim_{n \rightarrow \infty} 2 \int_0^x \frac{\beta'_n}{\sigma_n^2} (y) dy
\end{equation}
in $C^0$, independently from the mollifier. 
\item The function $\Sigma$ in \eqref{Sigma}
 is upper and lower bounded, and belongs to $C^{\alpha}_{\textup{loc}}$.
\end{enumerate} 
\end{hypothesis}


The following  proposition and definition  are given in \cite{frw1}, see respectively   Proposition 2.3 and the definition  at page 497.
 \begin{proposition}
 \label{P:equiv} Hypothesis \ref{H:Sigma}-1 is equivalent to ask that 
 there is a solution $h \in C^1$ to $L h=0$ such that $h(0) =0$ and 
 \begin{equation}\label{h'}
h'(x) := e^{-\Sigma(x)},  \quad x \in \R.
\end{equation}
In particular, $h'(0) =1$, and  $h'$ is strictly positive so 
	the inverse function by $h^{-1}: \R \rightarrow \R$ is well defined. 

 \end{proposition}
Let $\mathcal D_{L}$ be the set of functions $f\in C^1$ such that there is $\phi \in C^1$ with 
\begin{equation}\label{Lf_bis}
	f' = e^{-\Sigma} \phi.
\end{equation}
 For any $f \in \mathcal D_{L}$, we set
 \begin{equation}\label{Lf}
 	 L f= \frac{\sigma^2}{2} (e^{\Sigma} f')' e^{-\Sigma}.
 \end{equation}
This defines without ambiguity $L : \mathcal D_{L} \subset C^1 \rightarrow C^0$.
We also define
\begin{align}
	\mathcal{D}_{\mathcal L_M}&:= 
\mathcal{D}_{L} \cap  C^{1+ \alpha}_{\textup{loc}} \cap  C_b^0,
\label{barD}
\end{align}
equipped  with the graph topology of $L
$, the natural topology of $C^{1+ \alpha}_{\textup{loc}}$ and the uniform convergence topology. 

Then we  consider  
a transition kernel $Q(\cdot,dx)$ from
  $(\R, \mathcal{B} (\R))$   into $(\R,\mathcal{B} (\R))$, with $Q(y,\{0\})=0$,
satisfying the following condition.  
	\begin{hypothesis}\label{H:Kmeas}
For all
$B \in \mathcal B(\R)$,	 the map 	
$
y \mapsto \int_{B} (1 \wedge |x|^{1+\alpha}) \, Q(y,dx)  
$
is bounded and the measure-valued map
$y \mapsto  (1 \wedge |x|^{1+\alpha}) \, Q(y,dx)=: \tilde Q(y, dx)$
is continuous in the total variation topology.
\end{hypothesis}

For every $f\in \mathcal{D}_{\mathcal L_M}$,  we finally introduce the operator 
%
%
\begin{equation}\label{Ldistrib}
 \mathcal L_M f(y) := L f(y)
+  \int_{\R \setminus 0} (f(y + x) -f(y)
-k(x)f'(y)) Q(y,dx). 
\end{equation} 
Under Hypothesis \ref{H:Kmeas}, the operator above  takes values in 
	 $C^0$, see \cite{BandiniRusso_DistrDrift}. 
	 
	 In \cite{BandiniRusso_DistrDrift} we study the following Markovian martingale problem.	 
\begin{definition}\label{D:mtpb_hom_Mark}
We say that $(X, \P)$ fulfills the time-homogeneous Markovian martingale problem with respect to $\mathcal {D}_{\mathcal L_M}$ in \eqref{barD}, $\mathcal L_M$ in \eqref{Ldistrib}  and $x_0 \in \R$, if for any  $f \in \mathcal {D}_{\mathcal L_M}$, 
the process
\begin{align*}
f(X_{\cdot}) - f(x_0) - \int_0^{\cdot} \mathcal{L}_{M} f(X_{s}) ds
\end{align*}
is an $(\mathcal F^X_t)$-local martingale  under $\P$.
\end{definition}

\begin{remark}\label{R:twomartgpb}	 
Under Hypotheses \ref{H:Sigma} and \ref{H:Kmeas}, 	in \cite{BandiniRusso_DistrDrift} we provide existence and uniqueness for the  Markovian martingale problem in Definition \ref{D:mtpb_hom_Mark}.	
  Moreover, the solution $(X, \P)$ is a finite quadratic variation process, with  $\nu^{X,\P}(ds\, dx)= Q(X_{s-},dx)ds$.  
\end{remark}  


We can state the following theorem. 

\begin{theorem}\label{T:thmpassage}
	Assume that Hypotheses \ref{H:Sigma} and \ref{H:Kmeas} hold true. If $(X, \P)$ is a solution to the martingale problem in Definition \ref{D:mtpb_hom_Mark}, then it is a solution to the martingale problem in Definition \ref{D:mtpb1} with respect to $x_0$, 
\begin{equation}\label{DAm}
	\mathcal D_{\mathcal A}:=C^1([0,\,T]; {\mathcal D}_{\mathcal L_M}),
\end{equation}
and
\begin{align}\label{Adistrib2}
&(\mathcal{A} v)(dt,\eta):=    \partial_t v(t,\eta_{t})dt\\
 &+L v(t,\eta_{t})dt
+  \int_{\R} (v(t,\eta_{t} + x) -v(t,\eta_{t})
-k(x)\partial_y v(t,\eta_{t})) Q(\eta_{t},dx)d t, \quad v \in \mathcal D_{\mathcal A},\,\,\eta  \in D_{-}(0,\,T).\notag
\end{align}
Moreover, that martingale problem meets uniqueness.  
\end{theorem}
\begin{remark}\label{R:densityDA}
By \cite{BandiniRusso_DistrDrift}, $\mathcal D_{\mathcal L_M}$ in \eqref{barD} is dense in $C^1$.
Then, by Lemma  \ref{L:density},
  $\mathcal D_{\mathcal A}
  $ in \eqref{DAm}
	 is dense in $C^{0,1}([0,T] \times \R)$. 
\end{remark}
\noindent\emph{Proof of Theorem \ref{T:thmpassage}.}
A solution $(X, \P)$ to the martingale problem in Definition \ref{D:mtpb_hom_Mark} exists by Remark \ref{R:twomartgpb}. 
By  the equivalence (ii)-(iii) in  Theorem \ref{R:passage},  $(X, \P)$	fulfills the time-inhomogeneous martingale problem in Definition \ref{D:mtpb1} with respect to $x_0$, $\mathcal D_{\mathcal A}$ in \eqref{DAm} and $\mathcal A$ in \eqref{Adistrib2}.
 
Concerning uniqueness, given a solution $(X, \P)$ to the martingale problem by Definition \ref{D:mtpb1}, we know that,  for every  $v \in \mathcal D_{\mathcal A}$, $M^v$ is a  local martingale. In particular this holds for $v$ not depending on time, which implies the validity of the martingale problem not depending on time in the sense of Definition \ref{D:mtpb_hom_Mark}, for which we have uniqueness, see Remark \ref{R:twomartgpb}. 
\qed


Below we discuss some other properties of the solution to our martingale problem. 
We evaluate first the quadratic variation of the martingale component of a function of $X$ belonging to  $\mathcal D_{\mathcal{A}}$.  We can in particular   apply 
Corollary \ref{C:3.34_bis} to the case of the operator $\mathcal{A}$ in \eqref{Adistrib2}.

\begin{proposition}\label{L:418}	
Let $(X, \P)$ be a solution to the martingale problem in Definition \ref{D:mtpb_hom_Mark}.
Then, for every  $v \in \mathcal D_{\mathcal A}$ in \eqref{DAm} we have
\begin{equation}\label{Ycn_sigma}
	\langle Y^{v,c}, Y^{v,c} \rangle_t = \int_0^t \sigma^2(X_s) (\partial_x v(s, X_s))^2 ds,  
\end{equation}
where     $Y^{v, c}_t$ denotes the  unique continuous martingale part of $v(t, X_t)$.
In particular,  $v \mapsto Y^{v,c}$, $\mathcal D_{\mathcal A} \subseteq C^{0,1}([0,T] \times \R)\rightarrow \mathbb D^{ucp}$, 	is continuous in zero.
\end{proposition}
\noindent \emph{Proof of Proposition \ref{L:418}.}
By Theorem \ref{T:thmpassage},  $(X, \P)$ is a solution to the martingale problem in Definition  \ref{D:mtpb1} with respect to $x_0$, $\mathcal D_{\mathcal A}$ in \eqref{DAm} and $\mathcal A$ in \eqref{Adistrib2}. 

Assume that \eqref{Ycn_sigma} holds. 
Let $v_n \in \mathcal D_{\mathcal A}$ such that  $v_n \rightarrow 0$ in $C^{0,1}([0,T] \times \R)$. By  \eqref{Ycn_sigma}, $\langle Y^{v_n,c}, Y^{v_n,c} \rangle_T$ converges to zero in probability. 
Then, by Problem 5.25 in \cite{ks}, Chapter 1, 
it follows  that,  
\begin{equation}\label{supYfn}
Y^{v_n, c} \rightarrow 0 \quad \textup{u.c.p.}\,\ \textup{as}\,\,n \rightarrow \infty. 
\end{equation}
It remains to  prove \eqref{Ycn_sigma}. 
Take	 $v \in \mathcal D_{\mathcal A}$.  Then formula 
\eqref{Ycn_sigmaA} 
 in  Corollary \ref{C:3.34_bis} applied to  $\mathcal A$ in \eqref{Adistrib2} (taking into account Remark \ref{R:twomartgpb})
 yields
\begin{align*}
	\langle Y^{v,c}, Y^{v,c} \rangle_t &= \int_0^t (\mathcal A v^2 - 2 v \mathcal A v)(ds,X^{-})  - \int_{]0,t]\times \R}[v(s,X_{s-}+x)-v(s,X_{s-})]^2  Q(X_{s-},dx)ds\\
	&= \int_0^t(\partial_s v^2(s,X_s)- 2 v(s,X_s) \partial_s v(s,X_s))ds+\int_0^t (L v^2 - 2 v L v)(s,X_s) ds\\
	&+ \int_{]0,t]\times \R} (v^2(s,X_{s-} + x) -v^2(s,X_{s-})-2 k(x)\,  v(s,X_{s-})\, \partial_x v(s,X_{s-})\,) Q(X_{s-},dx)ds\\
	&- 2\int_{]0,t]\times \R}v(s,X_s)  (v(s,X_{s-} + x) -v(s,X_{s-})-k(x)\,\partial_x v(s,X_{s-})\,) Q(X_{s-},dx)ds\\
	& -  \int_{]0,t]\times \R}[v(s,X_{s-}+x)-v(s, X_{s-})]^2 Q(X_{s-},dx)ds\\
	&=  \int_0^t (L v^2 - 2 v L v)(s,X_s) ds = \int_0^t \sigma^2(X_s) (\partial_x v(s,X_s))^2 ds,
\end{align*}
where the latter equality  follows from the fact that  $L v^2 = 2 v L v + (\sigma \partial_x v)^2$, see  
Propositions 2.10 in \cite{frw1}.
\qed

\begin{corollary}\label{C:newcor}
 Let $(X, \P)$ be a solution to the martingale problem in Definition \ref{D:mtpb_hom_Mark}.
\begin{itemize}
	\item [(i)]For every $v \in C^{0,1}$, $Y^v:=v(\cdot, X)$ is a weak Dirichlet process. In particular, $X$ is a weak Dirichlet process. 
	\item [(ii)]\eqref{Ycn_sigma} holds for every $v \in C^{0,1}$. 
\end{itemize}
\end{corollary}
\proof
By Theorem \ref{T:thmpassage},  $(X, \P)$ is a solution to the martingale problem in Definition  \ref{D:mtpb1} with respect to $x_0$, $\mathcal D_{\mathcal A}$ in \eqref{DAm} and $\mathcal A$ in \eqref{Adistrib2}.

(i)
By Remark \ref{R:densityDA}, 
  $\mathcal D_{\mathcal A}
  $ in \eqref{DAm}
	 is dense in $C^{0,1}([0,T] \times \R)$. 	
Thanks to Proposition  \ref{L:418},  $v \mapsto Y^{v,c}$ is a continuous map. Since $X$ is a finite quadratic variation process (see Remark \ref{R:twomartgpb}),  we can  apply Corollary \ref{C:3.30bis}, which states that $X$ is a weak Dirichlet process.  Theorem \ref{P:3.10} concludes the proof of item (i).

(ii) Let $v \in C^{0,1}$. Since $\mathcal D_{\mathcal A}$ is dense in $C^{0,1}$ (see Remark \ref{R:densityDA}), there exists a sequence $(v_n)\in \mathcal D_{\mathcal A}$ converging to $v$ in $C^{0,1}$. Since $v \mapsto  Y^{v,c}$ is continuous, $Y^{v_n, c}\rightarrow Y^{v,c}$, u.c.p.  By Proposition \ref{P:App}, $\langle Y^{v_n, c}, Y^{v_n, c}\rangle \rightarrow \langle Y^{v, c}, Y^{v, c} \rangle$. 
The result follows from  Proposition \ref{L:418} since $v \mapsto \int_0^t \sigma^2(X_s) (\partial_x v(s, X_s))^2 ds$ is continuous.
\endproof

\begin{theorem}\label{T:G}
Let $(X, \P)$ be a solution to the martingale problem in Definition \ref{D:mtpb_hom_Mark}.
 Then there exists an $(\mathcal F_t)$-Brownian motion $W^X$ such that
	\begin{align}\label{decompX}
X&= x_0+ \int_0^\cdot \sigma(X_s) dW^X_s + \int_{]0,\cdot]\times \R} 
k(x)\,(\mu^X(ds\,dx)- Q(X_{s-},dx)d s)
 +\lim_{n \rightarrow \infty}\int_0^\cdot  L f_n(X_{s})ds\notag\\
&+ \int_{]0,\cdot]\times \R} (x - k(x))\mu^X(ds\,dx), 
\end{align}
for every sequence $(f_n)_n \subseteq \mathcal D_{\mathcal L_M}$   such that $f_n \underset{n \rightarrow \infty}{\rightarrow} Id$ in $C^{1}$. The limit appearing in \eqref{decompX} holds in the u.c.p. sense. 
\end{theorem}
\proof
By Theorem \ref{T:thmpassage},  $(X, \P)$ is a solution to the martingale problem in Definition  \ref{D:mtpb1} with respect to $x_0$, $\mathcal D_{\mathcal A}$ in \eqref{DAm} and $\mathcal A$ in \eqref{Adistrib2}. 
%
By Corollary \ref{C:genchar} we have 
\begin{align*}
X= X^{ c} + 
k(x)\,\star (\mu^X- \nu^{X, \P}) + \Gamma^{k}(Id) +  (x - k(x))\star \mu^X,    
\end{align*}
where $\Gamma^k$ is the operator defined in Theorem \ref{T:2.11}. By Proposition \ref{R:th_chainrule} and Remark \ref{R:linmapGamma}, $\Gamma^k$ is well defined in particular on $C^{0,1}$. 
By Proposition \ref{L:418}, we can apply Proposition \ref{T:3.30_bis} with $\mathcal {D}^S = \mathcal D_{\mathcal A}$, which yields that $v \mapsto \Gamma^{k}(v)$ restricted  to  $\mathcal D_{\mathcal A}$ is continuous. Since  $\mathcal D_{\mathcal A}$ is dense in  $C^{0,1}$,
 $\Gamma^{k}(v)$ is the continuous extension of the map  defined in \eqref{Gammavk}.

 At this point we evaluate $\Gamma^k(Id)$. Take $(f_n)_n \subseteq  \mathcal D_{\mathcal L_M}$  such that  $f_n \rightarrow Id$ in $C^{1}$.
By \eqref{Gammavk} together with \eqref{Adistrib2}, we get 
\begin{align*}
\Gamma^{k}(Id)&= \lim_{n \rightarrow \infty}\Gamma^{k}(f_n)\\
&=\lim_{n \rightarrow \infty}\Big(\int_0^\cdot (\mathcal A f_n)(ds, X^{-})  - \int_{]0,\cdot]\times \R}(f_n( X_{s-}+x)-f_n( X_{s-})) \frac{x - k(x)}{x}Q(X_{s-},dx)d s\Big)\\
& =\lim_{n \rightarrow \infty}\Big(\int_0^\cdot L f_n(X_{s-})ds
+  \int_{]0,\cdot]\times \R} (f_n(X_{s-} + x) -f_n(X_{s-})
-k(x)\,f_n'(X_{s-})\,) Q(X_{s-},dx)d s\\
&- (f_n( X_{s-}+x)-f_n( X_{s-})) \frac{x - k(x)}{x}\star \nu^{X, \P}\Big)\\
&=\lim_{n \rightarrow \infty}\int_0^\cdot L f_n(X_{s})ds.
\end{align*}
In order to get \eqref{decompX} it remains to identify $X^{c}$. 
Setting 
$
W^X_t := \int_0^t \frac{1}{\sigma(X_s)}d X^{c}_s 
$,
we have
$$
\langle W^X, W^X\rangle_t= \int_0^t \frac{1}{\sigma(X_s)}d \langle X^{c}, X^{c}\rangle_s. 
$$
On the other hand, Corollary \ref{C:newcor}-(ii) with $v \equiv \textup{Id}$ yields
$
\langle X^{c}, X^{c}\rangle=\int_0^\cdot \sigma^2(X_s) ds,
$
which implies that 
$\langle W^X, W^X\rangle_t= t$,  
and by L\'evy representation theorem, $W^X$ is a Brownian motion. Finally, by construction we have
$X^{c} = \int_0^\cdot \sigma(X_s) dW_s$. 
\endproof

\subsubsection{Continuous path-dependent SDEs with distributional drift}

Let  $\sigma$,   $\beta$ be continuous real functions and  $L: \mathcal D_L \rightarrow C^0$ in \eqref{Lf},
with $\mathcal D_L$ the  subset of $C^1$ introduced just before. We only suppose item 1 of Hypothesis \ref{H:Sigma}. 
Let $G^d:[0,\,T] \times  D_{-}(0,\,T)\rightarrow \R$ be a  Borel functional,  uniformly continuous on closed balls, 
and define
$$
\tilde G^d(t,\eta) = \frac{G^d(t,\eta)}{\sigma(\eta(t))}, \quad (t, \eta) \in [0,T] \times D_{-}(0,T).
$$
We suppose moreover that $\tilde G$ is bounded.
\normalcolor
We set 
 $G:[0,\,T] \times  C(0,\,T)\rightarrow \R$ as the restriction of $G^d$. 
	In \cite{ORT1_PartI} one investigates the martingale problem related to a path-dependent SDE of the type
\begin{equation}\label{eq1}
	d X_t = \sigma(X_t) dW_t + (\beta'(X_t) + G(t, X^t))dt.
\end{equation}
We set $ \mathcal D_{L,b}:=  \mathcal D_L \cap C^0_b$.
\begin{definition}\label{D:ORT}
$(X, \P)$ is solution to the (non-Markovian) path-dependent martingale problem related to \eqref{eq1} an initial condition $x_0$ if, for every $f \in \mathcal D_{L,b}$, 
\begin{equation}\label{ContdistrM}
M^f :=f(X_t) - f(x_0) - \int_0^t ((Lf) (X_s)+ f'(X_s) G(s, X^s)) ds	
\end{equation}
is an $(\mathcal F^X)$-local martingale under $\P$. 
\end{definition}
\begin{proposition} \label{P432}
  Let  $(X, \P)$ be a solution to the martingale problem
  in the sense of Definition \ref{D:ORT}. Then 
  %
  $X$ is necessarily a continuous process.
\end{proposition}
\proof
Let $h$ be the function introduced in Proposition \ref{P:equiv}, and set $Y=h(X)$. 
For every $\phi \in C^2$, we denote by  $L^0$   the classical PDE operator 
$
L^0\phi(y) =\frac{(\sigma h')^2(h^{-1}(y))}{2} \phi''(y)
$.
By Propositions  2.13 in \cite{frw1},
 $\phi \in \mathcal D_{L^0}(=C^2)$ if and only if $\phi \circ h \in \mathcal D_{L}$. 	This implies that $\phi \in C^2_b$ if and only if  $\phi \circ h \in \mathcal D_{L,b}$. Moreover, 
 		$L (\phi \circ h) = (L^0 \,\phi)\circ h$
 		for every $\phi \in C^2$. 

Since  $(X, \P)$	fulfills the time-homogeneous martingale problem in Definition \ref{D:ORT}, 
   for every $f \in  \mathcal {D}_{L,b}$,
\begin{align*}
f(X_{\cdot}) - f(x_0) - \int_0^t ((Lf) (X_s)+ f'(X_s) G(s, X^s)) ds	
\end{align*}
is an $(\mathcal F_t^X)$-local martingale  under $\P$.
Setting  $y_0 = h^{-1}(x_0)$, this yields that, for every  $\tilde f \in C^2_b$, 
$$
\tilde f(Y_\cdot) - \tilde f(y_0)-\frac{1}{2}\int_0^{\cdot} (\sigma h')^2(h^{-1}(Y_{s}))  \tilde f''(Y_{s}) ds -\int_0^{\cdot}h' (h^{-1}(Y_s))G(s, h^{-1}(Y^s))\tilde f'(Y_{s})ds
$$
is an $(\mathcal F_t^Y)$-local martingale  under $\P$.
It follows from Theorem \ref{T: equiv_mtgpb_semimart} that $Y$ is a  semimartingale with characteristics 
$B = \int_0^\cdot b(s,\check Y) ds$, $ C = \int_0^\cdot c(s,\check Y) ds$, $\nu(ds\,dz) = 0$, where $b(s,\eta):=h' (h^{-1}(\eta(s)))G(s, h^{-1}(\eta^s))$  and 
			$
			c(s,\eta) := (\sigma h')^2(h^{-1}(\eta(s)))  
			$. Consequently $\mu^Y(ds\,dy)=0$, so $Y=h(X)$ 
		  is necessarily a continuous process, and the same holds for $X$.     
\endproof

\begin{lemma}\label{L:densityC0b}
	$ \mathcal D_{L,b}$ is dense in $ \mathcal D_L$ equipped  with its graph topology. In particular, $ \mathcal D_{L,b}$ is dense in $C^1$. 
\end{lemma} 
\proof 
We consider the sequence $(\chi_N)$  introduced in \eqref{chiN}. 
Let $f \in \mathcal D_L$. 
We define a sequence $(f_N)\subset C^1$  such that  $f_N(0)=f(0)$ and  $f'_N=\chi_N f'$.
$f_N \in \mathcal D_L$ since $f_N' = (\phi \chi_N)e^{-\Sigma}$ and $\phi \chi_N \in C^1$, where $\phi$ has been defined in \eqref{Lf_bis}. Now, each $f_N$ is a  bounded function since $f_N'$ has compact support. Clearly, $f_N \rightarrow f$ in $C^1$. 
Moreover, making use of \eqref{Lf} we get 
\begin{align*}
L f_N = \frac{\sigma^2}{2}(e^\Sigma \chi_N f')' e^{-\Sigma}=\frac{\sigma^2}{2}(\phi \chi_N)' e^{-\Sigma} \rightarrow \frac{\sigma^2}{2}\phi' e^{-\Sigma}= Lf \quad \textup{in}\,\,C^0.
\end{align*}
This concludes the proof.
\endproof

\begin{corollary}
	Existence and uniqueness of  a solution to the martingale problem in Definition \ref{D:ORT} holds true.
\end{corollary}
\proof
By Theorem 4.23  in \cite{ORT1_PartI}, there is a solution $(X,\P)$
in the sense of  Definition \ref{D:ORT}
which is even continuous. This shows existence.

Concerning uniqueness, let  $(X,\P)$ be a solution
of the martingale problem in the sense of  Definition \ref{D:ORT}.
By Proposition \ref{P432}, $f(X)$ is necessarily continuous for every $f \in  \mathcal {D}_{L,b}$.
Moreover 
 the process in \eqref{ContdistrM} is a martingale 
 also  for every $f \in \mathcal D_{L}$
not necessarily bounded.
Indeed, let $f \in \mathcal D_L$. By Lemma \ref{L:densityC0b} there is
a sequence $f_N \in \mathcal D_{L,b}$ converging to $f$ in $\mathcal D_{L}$.
This implies that $M^{f_N}$ converges to $M^f$ u.c.p.    
We remark that the space of continuous local martingales is closed with respect to the u.c.p. convergence topology so $M^f$ is again a continuous local
martingale.
The conclusion follows by
Proposition 4.24 in  \cite{ORT1_PartI} which
states uniqueness in the framework of continuous processes.
\endproof

\begin{proposition}\label{P:4.35}
	Let $(X, \P)$ be a solution to the martingale problem in Definition \ref{D:ORT}. Then $(X, \P)$ is a solution to the martingale problem in Definition \ref{D:mtpb1} with  respect to $x_0$,
	\begin{equation}\label{DAG}
	{\mathcal D}_{\mathcal A}:=C^1([0,\,T];\mathcal D_{L,b})	
	\end{equation}
	and
\begin{equation}\label{AG}
		(\mathcal A v)(ds , \eta) = (\partial_s v(s, \eta(s)) + L v(s, \eta(s)) + G^d(s, \eta^s) \partial_x v(s, \eta(s)))ds.
\end{equation}
\end{proposition}
\proof
We apply Theorem \ref{T:passage} with $\mathcal D_{\mathcal L}:=\mathcal D_{L,b}$ equipped with its graph topology. 
\endproof
\begin{remark}
  Consider the martingale problem in Definition \ref{D:mtpb1}  with respect to  $\mathcal D_{\mathcal{A}}$ in \eqref{DAG}, $\mathcal{A}$ in \eqref{AG}, and $x_0 \in \R$.  
Replacing $G^d$ in \eqref{AG} with another Borel extension of $G$ one gets the same solution to the martingale problem. 

\end{remark}

Proceeding analogously as for the proof of Theorem \ref{T:G} we can prove the following result.
\begin{proposition}
Let $(X, \P)$ be a solution to the martingale problem in Definition \ref{D:ORT}.
Then we have the following. 
\begin{itemize}
	\item [(i)]$X$ is a weak Dirichlet process.
	\item [(ii)]There exists an $(\mathcal F_t)$-Brownian motion $W^X$ such that
	\begin{align}\label{decompX2}
X&= x_0+ \int_0^\cdot \sigma(X_s) dW^X_s + \int_0^\cdot  G^d(s, X^s) ds+\lim_{n \rightarrow \infty}\int_0^\cdot  L f_n(X_{s})ds, 
\end{align}
for every sequence $(f_n)_n \subseteq \mathcal D_{L}$  such that $f_n \underset{n \rightarrow \infty}{\rightarrow} Id$ in $C^{1}$. The  limit in \eqref{decompX2} holds in the u.c.p. sense. 
\end{itemize}
\end{proposition}

\subsubsection{The PDMPs case}
\label{S:The PDMPs case}
	Let $X$ be a  piecewise deterministic Markov process (PDMP)
generated by a marked point process $(T_n, \zeta_n)$, where $(T_n)_n$ are increasing random times such that 
	$ 
	T_n \in ]0,\,\infty[, 
	$ 
	where either there is a finite number of 
 times  $(T_n)_n$  or $\lim_{n \rightarrow \infty} T_n = + \infty$, 
	and $\zeta_n$ are random variables in $[0,1]$.
	We will follow the notations in \cite{Da-bo}, Chapter 2, Sections 24 and 26. 	The behavior of the PDMP $X$ is described by a triplet of 
	local characteristics  $(h,\lambda,Q)$: 
	$h: ]0,\,1[ \rightarrow \R$  is a Lipschitz continuous function, 
	$\lambda: ]0,1[ \rightarrow \R$ is a measurable function such that 
	$\sup_{x \in ]0,1[}|\lambda(x)| < \infty$, 
	and $Q$ is a transition probability  measure on $[0,1]\times\mathcal{B}(]0,1[)$. 	Some other technical assumptions  are specified in the over-mentioned reference, that we do not recall here.
	Let us denote by $\Phi(s,x)$  the unique solution of $
	{g}'(s) = h(g(s))$, $g(0)= x$. Then 
		\begin{equation} \label{X_eq}
			X(t)= 
	\left\{ 
	\begin{array}{ll} 
	\Phi(t,x),\quad t \in [0,\,T_{1}[\\ 
	\Phi(t-T_n, \zeta_n),\quad t \in [T_n,\,T_{n+1}[,\,\, n \in \N,  
	\end{array} 
	\right. 
	\end{equation}  
	and, for any $x_0 \in [0,1]$,  verifies the equation (provided the second integral in the right-hand side is well-defined)
	\begin{equation}\label{PDP_dynamic_estended} 
	X_t=x_0 + \int_{0}^t h(X_s)\,ds + \int_{]0,\,t]\times \R}x\,\mu^X(ds\,dx)
	\end{equation} 
	with 
\begin{equation}\label{muPDMPS}
	\mu^X(ds\,dx) 
	=\sum_{n \geq 1} 1_{\{\zeta_{n}\in ]0,1[\}} \delta_{( T_n,\,\zeta_n- \zeta_{n-1})}(ds\,dx).  
\end{equation}
Moreover, we introduce the predictable process counting the number of jumps of  $X$ from the  boundary of its domain:
	\begin{equation}\label{p_ast} 
	p^{\ast}_t = \sum_{0 <s \leq t}\,1_{\{X_{s-} \in \{0,1\}\}}.
	\end{equation} 
The knowledge of $(h,\,\lambda,\,Q)$ completely specifies the law of $X$, see Section 24 in \cite{Da-bo}, and also Proposition 2.1 in \cite{BandiniPDMPs}. In particular, let $\P$ be the unique probability measure under which  
	the compensator of $\mu^X$ has the form 
	\begin{equation}\label{nuPDPs} 
	\nu^X(ds\, dx) = \tilde Q(X_{s-},\,dx)\,(\lambda(X_{s-})\,ds + d p^{\ast}_s), 
	\end{equation} 
	where   $\tilde Q(y, dx)= Q(y, y+ dx)$, and $\lambda$ is trivially extended to $[0,1]$ by the zero value. 
	Notice that $X$ is a finite variation process, so \eqref{int_small_jumps} holds true. 
According to Theorem 31.3 and subsequent Section 31.5 in \cite{Da-bo}, for every measurable absolutely continuous function $v:\R_+ \times \R \rightarrow \R$ such that $(v(t, X_{t-}+x)-v(t, X_{t-})) \star\mu^X \in \mathcal A_{\textup{loc}}^+$,  
\begin{align*}
&v(t, X_t)- v(0, x_0)\\
&- \int_{]0,t]} \Big(\partial_s v(s, X_{s-})+   h(X_{s-})\partial_x v(s, X_s)+\lambda(X_{s-})\int_{\R}(v(s, X_{s-}+x)-v(s, X_{s-}))\tilde Q(X_{s-}, dx)\Big)ds\\
&- \int_{]0,t]\times \R}(v(s, X_{s-}+x)-v(s, X_{s-}))\tilde Q(X_{s-}, dx)dp^\ast_s
\end{align*}
is an $(\mathcal F^X)$-local martingale under $\P$. 
Therefore, $(X, \P)$ solves the martingale problem in Definition \ref{D:mtpb1}  with respect  to $\mathcal A$, $\mathcal D_A:=C^{1}([0,\,T] \times \R)$ and $x_0$, 
with 
$$
	(\mathcal  A v)(ds, \eta):= \Lambda_1 v(s,\eta_{s-})\,\gamma_1(ds, \eta_{s-})+\Lambda_2 v(s,\eta_{s-})\,\gamma_2(ds, \eta_{s-}), \quad v \in \mathcal D_A,\,\eta \in D(0,\,T),
$$	
with, for any $y \in \R$, $\gamma_1(ds, y)=ds$, $\gamma_2(ds,y)= dp^\ast_s(y)$, and  
\begin{align*}
	\Lambda_1 v(s,y)&= \partial_s v(s, y) +  h(y)\partial_y v(s, y)+
	\lambda(y)\int_{\R} (v(s, y+x)-v(s, y))\,\tilde Q(y, dx),\quad y \in ]0,1[,\\
	\Lambda_2 v(s,y)&= \int_{\R} (v(s, y+x)-v(s, y))\,\tilde Q(y, dx), \quad y \in \{0,1\}.	
\end{align*}

\appendix 
\renewcommand\thesection{Appendix} 
\section{} 
\renewcommand\thesection{\Alph{subsection}} 
\renewcommand\thesubsection{\Alph{subsection}}

\subsection{Some technical results on the (weak) finite quadratic variation}
\setcounter{equation}{0}
 \setcounter{theorem}{0}
\begin{proposition}\label{P:ucpequivbrackets}
	Let $Y=(Y(t))_{t \in [0,T]}$ and $X=(X(t))_{t \in [0,T]}$  be respectively a c\`adl\`ag and  a continuous process.
	Then 
	$$
	[X,Y]_\varepsilon^{ucp}(t) = C_\varepsilon(X,Y)(t) + R(\varepsilon, t)
	$$
	with
	$
	R(\varepsilon, t)\underset{\varepsilon \rightarrow 0}{\rightarrow 0 } \,\,\textup{u.c.p.}
	$ 
\end{proposition}
\proof
See Proposition A.3 in \cite{BandiniRusso1}.
\endproof

\begin{lemma}\label{L:3.15}
	Let  $G_n : C(0,\,T) \rightarrow \R$, $n \in \N$,  be a sequence of functions  such that
	\begin{itemize}
		\item [(i)] $\sup_n ||G_n||_{var} \leq M \in [0,\,+ \infty[$,
		\item [(ii)] $G_n \underset{n \rightarrow \infty}{\rightarrow} 0$ uniformly. 
	\end{itemize}
		Then, for every  $g: [0,\,T] \rightarrow \R$ c\`agl\`ad,
		 $
		 \int_0^\cdot g \, d  G_n \underset{n \rightarrow \infty}{\rightarrow} 0$, uniformly.
\end{lemma}
\proof
For every $n \in \N$, let us define  the operator
	$T_n: D_{-}(0,\,T)\rightarrow C(0,\,T)$,
	$g  \mapsto \int_0^\cdot g\, d G_{n}$.
 We denote by $\mathcal E_{-}(0,\,T)$ the linear space  of c\`agl\`ad step functions of the type
$\sum_i c_i \one_{]a_i, b_i]}$.
We first notice that, 
by (ii), if $g \in \mathcal E_{-}(0,\,T)$, then 
$T_n(g) \underset{n \rightarrow \infty}{\rightarrow} 0$, uniformly.

On the other  hand, $\mathcal E_{-}(0,\,T)$ is dense in $ D_{-}(0,\,T)$, see Lemma 1, Chapter 3, in \cite{Bil99} (that lemma is written for c\`adl\`ag function; however the same follows for c\`agl\`ad functions since the time reversal of a c\`adl\`ag function is  c\`agl\`ad).  

Let now $g \in D_{-}(0,\,T)$. Since $g$ is bounded, we have
$
\sup_{s \in [0,\,T]} |g(s)| \leq m	
$
for some constant $m$. 
 We get 
$
||T_n (g)||_\infty \leq M m. 
$
The conclusion follows by the   Banach-Steinhaus theorem, see e.g.  Chapter 1.2 in \cite{ds}. 
\endproof

\begin{proposition}\label{P: weakfinconv}
	Let $X$ be an $\mathbb F$-weak Dirichlet with weakly finite quadratic variation. Let $g$ be a c\`agl\`ad  process, and $N$ be a continuous $\mathbb F$-local martingale. Then
	$$
	\int_0^t g(s) (X_{s+\varepsilon}-X_s)(N_{s+\varepsilon}-N_s)\frac{ds}{\varepsilon} \,\,\underset{\varepsilon \rightarrow0}{\rightarrow} \,\,\int_0^t g(s)\, d[X^c, N]_s,  \quad t \in [0,T], \quad \textup{u.c.p}.
	$$
\end{proposition}
\proof
For $\varepsilon >0$, we set 
$F_\varepsilon(t) :=  C_{\varepsilon}(X,N)(s)$.
By Proposition \ref{P:uniqdec}, $X= X^c+ A$ with $A$ a martingale orthogonal process. Therefore, recalling Proposition \ref{P:ucpequivbrackets}, 
$$
F_\varepsilon(t)\underset{\varepsilon \rightarrow0}{\rightarrow} F(t):= [X^c, N]_t, \quad \textup{u.c.p}.
$$ 
Let  $\varepsilon_n$ be a sequence converging to zero as $n \rightarrow \infty$. It is sufficient to show the existence of a subsequence, still denoted by $\varepsilon_n$, such that 
\begin{equation}\label{gFconv}
	\int_0^\cdot g(s) \,d F_{\varepsilon_n}(s) \,\,\underset{n \rightarrow \infty}{\rightarrow} \,\,\int_0^\cdot g(s)\, d F(s), \quad \textup{u.c.p}.	
\end{equation}
By extracting a sub-subsequence, there exists a null set $\mathcal N$ such that 
\begin{equation}\label{ConvFeps}
F_{\varepsilon_n}(t) \underset{n \rightarrow \infty}{\rightarrow} F(t), \quad \textup{uniformly for all}\,\, \omega \notin \mathcal N.
\end{equation}
We remark that $N$ has also weakly finite quadratic variation.  Let $\kappa >0$. For any $\ell > 0$, we denote by $\Omega_{n, \ell}$ the subset of $ \omega \in \Omega$ such that 
\begin{align}\label{estquadr}
	&\int_0^T (X_{(s+\varepsilon_n)\wedge T}(\omega)-X_{s}(\omega))^2\frac{ds}{\varepsilon_n} + \int_0^T (N_{(s+\varepsilon_n)\wedge T}(\omega)-N_{s}(\omega))^2\frac{ds}{\varepsilon_n} \leq \ell,\notag\\
&	\langle X^c, X^c\rangle_T + \langle N, N\rangle_T\leq \ell.
\end{align}
In particular, we can choose $\ell$ such that $\P(\Omega_{n, \ell}^c) \leq \kappa $. Moreover, since $g$ is locally bounded, without restriction of generality we can take
\begin{equation}\label{estg}
\sup_{s \in [0,\,T]} |g(s)| \leq \ell, \quad \forall \omega \in \Omega_{n, \ell}. 	
\end{equation}
Collecting \eqref{estquadr} and \eqref{estg}, on $\Omega_{n, \ell}$ we have 
\begin{align}\label{estell}
	&\sup_{0 \leq t \leq T} \left|\int_0^t g(s) \,d  F_{\varepsilon_n}(s)\right| \leq \ell  
	||F_{\varepsilon_n}||_{\textup{var}}\notag\\
	&\leq \ell \sqrt{\int_0^T (X_{s+\varepsilon_n}(\omega)-X_s(\omega))^2\frac{ds}{\varepsilon_n}} \sqrt{\int_0^T (N_{s+\varepsilon_n}(\omega)-N_s(\omega))^2\frac{ds}{\varepsilon_n}}\notag\\
	&= \ell \sqrt{\int_0^T (X_{(s+\varepsilon_n) \wedge T}(\omega)-X_s(\omega))^2\frac{ds}{\varepsilon_n}} \sqrt{\int_0^T (N_{(s+\varepsilon_n)\wedge T}(\omega)-N_s(\omega))^2\frac{ds}{\varepsilon_n}}\leq \ell^2,\notag\\
	&\sup_{0 \leq t \leq T} \left|\int_0^t g(s) \,d  F(s)\right| \leq \ell 
	||F||_{\textup{var}}\leq \ell \sqrt{\langle X^c, X^c\rangle_T
	\langle N, N\rangle_T
	}\leq \ell^2.
\end{align}
Let us come back to prove \eqref{gFconv}. We set 
$$
\chi_n(g) := \sup_{0 \leq t \leq T} \Big| \int_0^t g(s) \,d  F_{\varepsilon_n}(s) - \int_0^t g(s)\, d F(s) \Big|.
$$
Let $K >0$.
Using \eqref{estell}, together with  Chebyshev's inequality, we get
\begin{align}\label{ineqchin}
	\P\left(\chi_n(g)> K\right) 
	&\leq \P(\Omega_{n, \ell}^c)+\P\left(\left \{\chi_n(g)> K\right\} \cap \Omega_{n, \ell} \cap \mathcal N^c\right)\notag\\
	&\leq \kappa +\P\left(\left \{\chi_n(g)> K\right\} \cap \Omega_{n, \ell} \cap \mathcal N^c\right)\notag\\
	&= \kappa +\P\left(\left \{\chi_n(g) \wedge 2 \ell^2 > K\right\} \cap \Omega_{n, \ell} \cap \mathcal N^c\right)\notag\\
		&= \kappa +\P((\chi_n(g) \wedge 2 \ell^2)\one_{ \Omega_{n, \ell}\cap \mathcal N^c} > K)\notag\\
	&\leq \kappa +\frac{\E\left[(\chi_n(g) \wedge 2 \ell^2) \one_{ \Omega_{n, \ell}\cap \mathcal N^c} \right]}{K^2}.
\end{align}
To prove that previous expectation goes to zero as $n$ goes to infinity, by using Lebesgue's dominated  convergence theorem, it remains to show that
\begin{equation}\label{chiconv}
 \chi_n(g)\one_{ \Omega_{n, \ell}\cap \mathcal N^c}  \underset{n \rightarrow \infty}{\rightarrow} 0 \quad \textup{a.s.}
\end{equation}
To this end, let us set $G_n:= \one_{\Omega_{n,\ell}}(F_n-F)$. Since 
by  \eqref{estquadr} we have 
$
||G_n||_{\textup{var}} \leq \frac{\ell}{2}$, 
 the convergence in \eqref{chiconv} follows from  Lemma \ref{L:3.15}.
 Consequently, by \eqref{ineqchin}, 
$\limsup_{n \rightarrow \infty} \P\left(\chi_n(g)> K\right) 
	\leq \kappa$.
Since $\kappa>0$ is arbitrary,  the proof is concluded.
\normalcolor
\qed

	\begin{proposition}\label{P:App}
		Let $(M^n(t))_{t \in [0,T]}$ (resp. $(N^n(t))_{t \in [0,T]}$) be a sequence of continuous local martingales, converging u.c.p. to $M$ (resp. to $N$). Then 
	\begin{align*}
	[M^n,N^n] \underset{n\rightarrow \infty}{\longrightarrow} [M,N]\quad \textup{u.c.p}.
\end{align*}
\end{proposition}
In order to prove Proposition \ref{P:App} we first give a technical result. 
\begin{lemma}\label{L:D6}
	Let $A, \delta >0$. Let $M$ be a continuous local martingale vanishing at zero. We have 
\begin{equation}\label{Eq:lemma}
	\P\Big([M,M]_T \geq A\Big) \leq \frac{\E\Big[ \max_{t \in [0,T]}|M_t|^2\wedge \delta^2 \Big]}{A} + \P\Big(\max_{t \in [0,T]}|M_t| \geq \delta\Big). 
\end{equation}
\end{lemma}
\proof
We bound the left-hand side of \eqref{Eq:lemma} by
$$
\P\Big([M,M]_T \geq A, \max_{t \in [0,T]}|M_t|\leq \delta\Big) + 
\P\Big(\max_{t \in [0,T]}|M_t|\geq \delta\Big):= I_1 + I_2. 
$$
Let $\tau:= \inf\{s \in [0,T]: \,\,|M_s|\geq \delta\}$. 
We notice that on $\Omega_0 :=\{\omega \in \Omega: \max_{t \in [0,T]}|M_t(\omega)|\leq \delta\}$ we have $M=M^\tau$. Therefore, by the definition of covariation, $[M, M]=[M^\tau, M^\tau]$ on $\Omega_0$, so that 
$$
I_1 = \P\Big([M^\tau,M^\tau]_T \geq A, \max_{t \in [0,T]}|M_t|\leq \delta\Big).
$$ 
Using Chebyshev and Burkholder-Davis-Gundy inequalities, 
we get that there is a constant $c>0$ such that 
\begin{align*}
	I_1 
	&\leq \P\left([M^\tau,M^\tau]_T \geq A\right)
	\leq \frac{\E[[M^\tau, M^\tau]_T]}{A}
	\leq \frac{c\,\E\Big[\sup_{t \in [0,T]}|M_t^\tau|^2\Big]}{A}
	= \frac{c\,\E\Big[\sup_{t \in [0,T]}|M_t^\tau|^2 \wedge\delta^2\Big]}{A}\\
	&\leq  \frac{c\,\E\Big[\sup_{t \in [0,T]}|M_t|^2 \wedge\delta^2\Big]}{A}.
\end{align*}
\qed

\noindent \emph{Proof of Proposition \ref{P:App}.}
Obviously, we can take $M\equiv N \equiv 0$. 
By polarity arguments, it is enough to suppose $(M^n) \equiv (N^n)$ and to prove that $[M^n, M^n]_T \rightarrow 0$ in probability.  Let $\varepsilon>0$, and let $N_0$ such that, for every $n \geq N_0$, 
$
\P(\sup_{t \in [0,T]}|M_t^n|\geq 1 )\leq \varepsilon.
$
Let $A>0$. For $n \geq N_0$, Lemma \ref{L:D6} gives 
\begin{equation*}
	\P\left([M^n,M^n]_T \geq A\right) \leq \frac{\E\Big[ \max_{t \in [0,T]}|M^n_t|^2\wedge 1\Big]}{A} + \varepsilon. 
\end{equation*}
By taking $n \rightarrow \infty$ we get 
\begin{equation*}
	\limsup_{n \rightarrow \infty}\P\left([M^n,M^n]_T \geq A\right) \leq  \varepsilon 
\end{equation*}
and the result follows from the arbitrariness of $\varepsilon$. 
\qed

\subsection{Some  results on densities of martingale problems' domains}
\setcounter{equation}{0}
\setcounter{theorem}{0}

We have the following density results.
\begin{lemma}\label{L:C1}
	Let $M >0$ and $E$ be a topological vector  $F$-space with some metric $d$. We introduce the distance  
	$$
	\bar d(e_1, e_2) := \sup_{t \in [0,M]}d(t e_1, t e_2), \quad e_1, e_2 \in E. 
	$$
	Then $d$ and $\bar d$ are equivalent in the following sense: if $\varepsilon >0$, there is $\delta >0$ such that 
	\begin{equation}\label{dequiv1}
		d(e_1, e_2)\leq \delta \Rightarrow \bar d(e_1, e_2) \leq \varepsilon, \quad e_1, e_2 \in E,
	\end{equation}
	and 
	\begin{equation}\label{dequiv2}
		\bar d(e_1, e_2)\leq \delta \Rightarrow  d(e_1, e_2) \leq \varepsilon, \quad e_1, e_2 \in E.
	\end{equation}
\end{lemma}
\proof
\eqref{dequiv2} is immediate since $d \leq \bar d$. 

\eqref{dequiv1} follows since we can easily prove that $e \mapsto \bar d(e,0)$ is continuous. This follows because $(t, e) \mapsto t e$ and therefore $(t,e) \mapsto d(te, 0)$ is continuous.
\endproof


\begin{definition}
	An $F$-space $E$ is said to be generated by a countable sequence of seminorms  $(||\cdot||_{\alpha\in \N})$ if it can be equipped with the distance 
	\begin{equation}\label{distFspace}
		d_E(x,y):= \sum_\alpha \frac{||x-y||_\alpha}{1+||x-y||_\alpha}2^{-\alpha}.
	\end{equation}
\end{definition}
The space $C^0([0,T]; E)$ will be equipped with the distance 
$
d(f,g) := \sup_{t \in [0,T]}d_E(f(t),g(t)).
$
\begin{remark}
	$C^1$, $C^2$, $C^0$, and $\mathcal D_L \cap C^\alpha_\textup{loc} \cap C_b^0$ (see \eqref{barD}) are $F$- spaces  generated by a countable sequence of seminorms.
\end{remark}
\begin{definition}\label{D:D4}
Let $E$ be an $F$-space  generated by a countable sequence of seminorms.
	Let $d_E$ be the distance in \eqref{distFspace} related to $E$. 
	A function $f : [0,T] \rightarrow E$ is said to be $C^1([0,T];E)$ if there exists $f':[0,T] \rightarrow E$ continuous such that 
	\begin{equation}\label{C1def}
		\lim_{\varepsilon \rightarrow 0}\sup_{t \in [0,T]}d_E\Big (\frac{f(t + \varepsilon)- f(t)}{\varepsilon}, f'(t)\Big )=0.
	\end{equation} 
\end{definition}
We remark that by Definition \ref{C1def} $f'$ is continuous.
We are not aware about Bochner integrability for functions taking values in $E$. For this we make use of Riemann integrability. The following lemma makes use of classical arguments, which exploits the fact that $f$ is uniformly continuous. 
\begin{lemma}\label{L:C5}
	Let $f :[0,T] \rightarrow B$ be a continuous function, where $B$ is a seminormed space. We denote by 
	$$
	s_n(f)(t) = \sum_{k=0}^{2^n-1}f(k t 2^{-n}) 2^{-n}, \quad t \in [0,T], 
	$$
	the Riemann sequence related to the dyadic partition. Then $(s_n(f))$ is Cauchy in $C^0([0,T]; B)$, namely
	$$
	\sup_{t \in [0,T]} ||(s_n(f)- s_m(f))(t)||_B\rightarrow 0 \quad \textup{as}\,\,n,m \rightarrow \infty. 
	$$  
\end{lemma}
 \begin{remark}\label{R:C6}
Let $E$ be an $F$-space  generated by a countable sequence of seminorms  $(||\cdot||_{\alpha\in \N})$.
Let $d_E$ be the distance in \eqref{distFspace} related to $E$.	Let $f_n, f_m$ (resp. $f$) be   sequences of functions (resp. a function) from $[0,T]$ to $E$.
	\begin{itemize}
	\item[(i)]$\sup_{t \in [0,T]} d_E(f_n(t) , f(t)) \underset{n \rightarrow \infty}{\rightarrow} 0$ is equivalent to $\sup_{t \in [0,T]} ||f_n(t)-f(t)||_\alpha\underset{n \rightarrow \infty}{\rightarrow} 0$;
\item[(ii)]$\sup_{t \in [0,T]} d_E(f_n(t) , f_m(t)) \underset{n,m \rightarrow \infty}{\rightarrow} 0$ is equivalent to $\sup_{t \in [0,T]} ||f_n(t)-f_m(t)||_\alpha \underset{n,m \rightarrow \infty}{\rightarrow} 0$.
	\end{itemize}
	\end{remark}
	\begin{remark}[Riemann integral]\label{R:RiemannIntegral}
	Let $E$ be an $F$-space  generated by a countable sequence of seminorms  $(||\cdot||_{\alpha\in \N})$. 
	 By Remark \ref{R:C6}, 
	$C^0([0,T];E)$ is an $F$-space.
	Let $f\in C^0([0,T]; E)$. By Lemma \ref{L:C5}, $(s_n(f))$ is Cauchy with respect to all seminorms $||\cdot||_{\alpha\in \N}$. By Remark \ref{R:C6}, 
	the sequence $(s_n(f))$  is Cauchy with respect to $d$. 
	Being $C^0([0,T];E)$  an $F$-space, it is complete, and $(s_n(f))$
	converges to an element in $C^0([0,T]; E)$ that we denote by $\int_0^\cdot f(s) ds$. For $a, b \in [0,T]$ we write $\int_a^b f(s) ds=\int_0^b f(s) ds-\int_0^a f(s) ds$.  It is not difficult to show that Riemann integral satisfies the following properties.

	\begin{enumerate}
	\item 	$\int_a^b f(t) dt =\int_{a-h}^{b-h} f(t+h) dt$, $h\in \R$ (by definition of the integral).
	\item $\lim_{\varepsilon \rightarrow 0}\frac{1}{\varepsilon}\int_t^{t+\varepsilon} f(s) ds \rightarrow f(t)$ in $C^{0}([0,T];E)$, $\lim_{\varepsilon \rightarrow 0}\frac{1}{\varepsilon}\int_{t-\varepsilon}^{t} f(s) ds \rightarrow f(t)$ in $C^{0}([0,T];E)$ (by Remark \ref{R:C6}-(i) and the fact that $f$ is uniformly continuous).
	\item If $f \in C^0([0,T];E)$, then $f(\cdot)=f(0) + \lim_{\varepsilon \rightarrow 0}\frac{1}{\varepsilon}\int_0^{\cdot} (f(s+ \varepsilon)-f(s)) ds $ (by items 1. and 2. above).
	\item  If $f \in C^1([0,T];E)$, then $f(\cdot)=f(0) + \int_0^{\cdot} f'(s) ds $ (by item 3.  above and Definition \ref{D:D4}).
	\end{enumerate}

	\end{remark}

\begin{lemma}\label{L:hatD}
Let $\mathcal D_{\mathcal L}$ be a topological vector  $F$-space generated by a countable sequence of seminorms, equipped with the metric $d_{\mathcal L}(=d_E)$ in \eqref{distFspace} related to $\mathcal D_{\mathcal L}(=E)$.

We suppose that $\mathcal D_{\mathcal L}$ is topologically embedded in  $C^0$.  Let  
 $\mathcal L: \mathcal D_{\mathcal L} \rightarrow C^0([0,\,T]\times D_{-}(0,\,T))$ be a continuous map. 
	Let  $\hat {\mathcal D}_{\mathcal A}$ be the subspace  of $ {\mathcal D}_{\mathcal A}:=C^1([0,\,T]; \mathcal D_{\mathcal L})$ 
	 constituted by  functions of the type
\begin{equation}\label{formun}
u(t,x)= \sum_{k} a_k(t) u_k(x), \quad a_k \in C^1(0,\,T), \,\,u_k \in \mathcal D_{\mathcal L}. 
\end{equation}
Then $\hat {\mathcal D}_{\mathcal A}$ is dense in ${\mathcal D}_{\mathcal A}$ equipped with the metric $d_{\mathcal A}$
 governing the following convergence: \\$u_n \rightarrow 0$ in ${\mathcal D}_{\mathcal A}$ if 
\begin{align*}
	u_n \rightarrow 0 \,\, \textup{in}\,\,C^0([0,\,T]; \mathcal D_{\mathcal L})\,\,\textup{and}\,\,
	\partial_t u_n \rightarrow 0 \,\, \textup{in}\,\,C^0([0,\,T]; \mathcal D_{\mathcal L}).
\end{align*}
\end{lemma}
\proof
Let  $u \in {\mathcal D}_{\mathcal A}$. 
 We have to prove that  there is a 
 sequence $u_n \in \hat {\mathcal D}_{\mathcal A}$ such that $u_n \rightarrow u$ in ${\mathcal D}_{\mathcal A}$.
We denote $v := u'$, i.e. the time derivative.
We divise the proof into two steps.

\noindent \emph{First step: approximation of $v $.} 
In this step we only use the fact that $E= \mathcal D_{\mathcal L}$ is an $F$-space.

Notice that $v  \in C^0([0,\,T]; \mathcal D_{\mathcal L})$.  Let $\varepsilon >0$. Since $v$ is uniformly continuous, by Lemma \ref{L:C1}, there exists $\delta >0$ such that 
\begin{equation}\label{C3}
 |t-s|\leq\delta \Rightarrow \bar d_{\mathcal L}(v(t), v(s))< \frac{\varepsilon}{2}, 
\end{equation}
where $\bar d_{\mathcal L}$ is the metric related to $d_{\mathcal L}$ introduced in Lemma \ref{L:C1} with $M=1$. 
We consider a dyadic partition of $[0,\,T]$ given by  $t_k =2^{-n} kT$,
 $k \in \{0,..., 2^n\}$, $n \in \N$.
 We define the open recovering of $[0,\,T]$ given by 
\begin{equation}
U_k^n =
\left\{\begin{array}{l}
[t_0, t_1[\quad k=0,\\
]t_{k-1}, t_{k+1}[ \cap[0,T]\quad k \in \{1,..., 2^{n}\}.
\end{array}\right.
\end{equation}
We also introduce a smooth partition of the unit $\varphi_k^n$, $k \in \{0,..., 2^{n}\}$. In particular, 
$\sum_{k=0}^{2^n} \varphi_k^n =1$, $ \varphi_k^n \geq 0$, and $ \supp \varphi_k^n \subset U_k^n$.
We define 
\begin{align*}
	v_n(t)&:= \sum_{k=0}^{2^n}v(t_k) \varphi_k^n(t).
\end{align*}
 We notice that, if $t \in [t_{k-1}, t_k[$, then $v_n(t) =v(t_{k-1}) \varphi_{k-1}^n(t)+v(t_k) \varphi_k^n(t)$ and $v(t) =v(t) \varphi_{k-1}^n(t)+v(t) \varphi_k^n(t)$. 
Since $d_{\mathcal L}$ is a homogeneous distance, using the triangle inequality, we have
\begin{align}\label{bard}
d_{\mathcal L}(v_n(t), v(t))&= d_{\mathcal L}((v(t_{k-1})-v(t)) \varphi_{k-1}^n(t)+(v(t_k)-v(t)) \varphi_k^n(t),0)\notag\\
&\leq d_{\mathcal L}((v(t_{k-1})-v(t)) \varphi_{k-1}^n(t),0)+d((v(t_k)-v(t)) \varphi_k^n(t),0)\notag\\
&\leq d_{\mathcal L}(v(t_{k-1}) \varphi_{k-1}^n(t),v(t)) \varphi_{k-1}^n(t))+d_{\mathcal L}(v(t_k) \varphi_k^n(t),v(t) \varphi_k^n(t))\notag\\
&\leq \bar d_{\mathcal L}(v(t_{k-1}),v(t))+ \bar d_{\mathcal L}(v(t_k),v(t)).
\end{align}
 Then we choose $N$ such that $2^{-N} T\leq \delta$. Recalling \eqref{C3},  we obtain from \eqref{bard} that 
 \begin{align*}
n > N \Rightarrow \sup_{t \in [0,T]} d_{\mathcal L}(v_n(t), v(t))
&\leq \varepsilon.
\end{align*}
This shows that $v_n \rightarrow v=\partial_t u \,\, \textup{in}\,\,C^0([0,\,T]; \mathcal D_{\mathcal L})$.


\noindent \emph{Second step: approximation of $u$.}
Now we set 
\begin{align*}
	u_n(t)&
	:=
	 u(0) + \sum_{k=0}^{2^n}v(t_k) \bar \varphi_k^n(t), 
\end{align*}
where  $\bar \varphi_k^n(t) := \int_0^t  \varphi_k^n(s) ds$.
We remark that  
$
		u_n(t) = u(0) + \int_0^t v_n(s) ds
$
in the sense of  Remark \ref{R:RiemannIntegral}. 
	We have to show that 
	$
		\sup_{t \in [0,T]} d_{\mathcal L}(u_n(t), u(t))
	$
	converges to zero. To this end, by Remark \ref{R:C6},  it is enough to show that, for every fixed $\alpha$,  
	\begin{equation}\label{convde}
	\sup_{t \in [0,T]} ||u_n(t)-u(t)||_\alpha\underset{\varepsilon \rightarrow 0}{\rightarrow} 0.	
	\end{equation}
	Since  $u$  is uniformly continuous with respect to $d_{\mathcal L}$, obviously it has the same property with respect to the seminorm $\alpha$. 
	 For every $t \in [0,T]$, 
	\begin{align*}
	||u_n(t)-u(t)||_\alpha  \leq 	\int_{0}^t ||v_n(s)-v(s)||_\alpha ds \leq T \sup_{s \in [0,T]} ||v_n(s)-v(s)||_\alpha \leq C \sup_{s \in [0,T]}\,d_{\mathcal L}(v_n(s),v(s)) 
	\end{align*}
	that converges to zero by the first step of the proof.
	This implies \eqref{convde} and  concludes the proof. 
\endproof

\begin{lemma}\label{L:density}
Let $\mathcal D_{\mathcal L_M}$ be a topological vector $F$-space
generated by a countable sequence of seminorms, and let ${\mathcal D}_{\mathcal A}$ be defined in \eqref{D:DAnewMark}. If  $\mathcal D_{\mathcal L_M}$ is dense in $C^1$, then ${\mathcal D}_{\mathcal A}$ is dense in $C^{0,1}$. 
\end{lemma}
\proof
We start by noticing that   $C^{0,1}=C^0([0,\,T]; C^1)$. 
Then we divise the proof into two steps. 
\noindent \emph{First step.} $C^0([0,\,T]; {\mathcal D}_{\mathcal L_M})$ is dense in $C^0([0,\,T]; C^1)$. 

In this step we only use the fact that $E:= C^1$ is an $F$-space. 
Let $d_E$ be the distance in \eqref{distFspace} related to $E$. 
Let $f \in C^0([0,\,T]; C^1)$. Let $\delta >0$. 
We need to show the existence of $f^\varepsilon \in C^0([0,\,T]; \mathcal D_{\mathcal L_M})$ such that 
$$
d_E(f(t), f^\varepsilon(t))\leq \delta, \quad \forall t \in [0,\,T].
$$
Let $0=t_1 < t_1<\cdots < t_n=T$ be a subdivision with mesh $\varepsilon$. Since $\mathcal D_{\mathcal L_M}$ is dense in $C^1$, using Lemma \ref{L:C1},  for every $i=1,..., n$, there is $f^\varepsilon_{t_i} \in  \mathcal D_{\mathcal L_M}$ such that 
\begin{equation}\label{bard2}
\bar d_E(f(t_i), f^\varepsilon_{t_i}) \leq \frac{\delta}{6},
\end{equation}
where $\bar d_E$ is the metric in Lemma \ref{L:C1} with $M=1$.
The candidate now is 
\begin{align*}
	f^\varepsilon(t) = f^\varepsilon_{t_i}+ \frac{t-t_{i-1}}{t_i - t_{i-1}}(f^\varepsilon_{t_i}-f^\varepsilon_{t_{i-1}}), \quad t \in [t_i, t_{i+1}[.
\end{align*}
To prove that $f^\varepsilon$ is a good approximation of $f$, we define $f^\pi: [0,\,T] \rightarrow C^1$ as 
$$
f^\pi(t)= f(t_{i-1}) + \frac{t-t_{i-1}}{t_i - t_{i-1}}(f(t_i)-f(t_{i-1})), \quad t \in [t_i, t_{i+1}[.
$$
We evaluate the difference between $f^\pi$ and $f^\varepsilon$. 
Setting $A:=\frac{t-t_{i-1}}{t_i - t_{i-1}}\in [0,\,1]$, taking into account the homogeneity of $d$ and the triangle inequality, we have  
\begin{align*}
d_E(f^\pi(t),f^\varepsilon(t))&=d\Big(f(t_{i-1}) + A(f(t_i)-f(t_{i-1})),f^\varepsilon_{t_{i-1}}+A(f^\varepsilon_{t_i}-f^\varepsilon_{t_{i-1}})\Big)\\
	&\leq d_E(f(t_{i-1}),f^\varepsilon_{t_{i-1}})+ \bar d_E(f(t_i),f^\varepsilon_{t_{i}})+\bar d_E(f(t_{i-1}),f^\varepsilon_{t_{i-1}})\leq 3 \frac{\delta}{6}= \frac{\delta}{2},
\end{align*}
where in the latter inequality we have used \eqref{bard2}. 
This shows that 
$$
\sup_{t \in [0,\,T]} d_E(f^\pi(t), f^\varepsilon(t)) \leq \frac{\delta}{2}.
$$
%
\emph{Second step}. $C^1([0,\,T]; {\mathcal D}_{\mathcal L_M})$ is dense in $C^0([0,\,T]; {\mathcal D}_{\mathcal L_M})$.

In this step we set $E:= {\mathcal D}_{\mathcal L_M}$. Let $f \in C^0([0,\,T]; {\mathcal D}_{\mathcal L_M})$, and set 
$
f_\varepsilon(t):=\frac{1}{\varepsilon}	\int_{t- \varepsilon}^t f(s) ds$, $t \in [0,T]$.
This integral is well-defined, 
and converges to $f$ by    Remark \ref{R:RiemannIntegral}-2.
\endproof

	\subsection{Some further technical results}
	\setcounter{equation}{0}
\setcounter{theorem}{0}

\begin{proposition}\label{P:new}
	Condition \eqref{int_small_jumps} is equivalent to ask that 
	$$
	  (1 \wedge |x|^2) \star \mu^X \in \mathcal A_{\textup{loc}}^+.
	$$
\end{proposition}
\proof
We start by noticing that
\begin{align*}
	  (1 \wedge |x|^2) \star \mu^X &= \int_\R \one_{|x| >1 }\,\mu^X(ds\,dx)+\int_\R \one_{|x| \leq1 }x^2\,\mu^X(ds\,dx)\\
	  &=\sum_{s \leq \cdot }\one_{\{|\Delta X_s|>1\}}+\sum_{s \leq \cdot} |\Delta X_s|^2 
		\,\one_{\{|\Delta X_s| \leq 1\}}.
\end{align*}
The process $\sum_{s \leq \cdot }\one_{\{|\Delta X_s|>1\}}$ is locally bounded having bounded jumps (see Remark \ref{R:pred}-2.), 
 and therefore it is locally integrable. 
On the other hand, since there is almost surely a finite number of jumps above one,  condition \eqref{int_small_jumps} is equivalent to 
\begin{equation}\label{Sfinite}
		S:=\sum_{s \leq \cdot} |\Delta X_s|^2 
		\,\one_{\{|\Delta X_s| \leq 1\}} 
		< \infty \,\,\,\,\textup{a.s.}
\end{equation}
Being $S$  a process with bounded jumps, \eqref{Sfinite} is equivalent to say that $S$ is locally bounded and therefore locally integrable. This concludes the result.
\endproof

In the sequel of the  section we consider  a c\`adl\`ag  process $X$  satisfying condition \eqref{int_small_jumps}.
	\begin{lemma}\label{L:c}
	For all $c>0$, 
	$$
	\one_{\{|x| >c\}}\star \mu^X \in \mathcal{A}_{\textup{loc}}^+.
	$$
	\end{lemma}
	\proof
It is enough to prove the result for $c \leq 1$. We have 
	\begin{align*}
		\one_{\{|x|>c\}}\star \mu^X = \frac{1}{c^2} c^2 \, \one_{\{|x|>c\}}\star \mu^X \leq \frac{1}{c^2} (|x|^2 \wedge 1) \star \mu^X
	\end{align*}
	and the result follows by \eqref{int_small_jumps} and Proposition \ref{P:new}. 
	\endproof

	\begin{lemma}\label{L:3.32}
Let $v:\R_+ \times \R \rightarrow \R$ be a locally bounded function. Then for all  $0 <a_0 \leq a_1$, 
$$
 |v(s,X_{s-}+x)-v(s,X_{s-})| \,\one_{\{a_0 <|x| \leq a_1\}} \star \mu^X \in \mathcal{A}^+_{\textup{loc}}.
 $$
\end{lemma}
\proof
Without restriction of generality, we can take $a_0 < 1 \leq a_1$. 
Since $X_{s-}$ is c\`agl\`ad, it is locally bounded, and therefore we can consider the localizing sequence $(\tau_n)_n$
 such that, for every $n \in \N$, $\tau_n=\inf\{s \in \R_+:\,\,|X_{s-}| \leq n\}$. Then, for every $n \in \N$, on $[0,\,\tau_n]$,  
\begin{align*}
\one_{[0,\,\tau_n]}(s)	|v(s,X_{s-}+x)-v(s,X_{s-})| \,\one_{\{a_0 <|x| \leq a_1\}} \star \mu^X
&\leq 2 \sup_{y \in (-n-a_1, n+a_1)}|v(s, y)|\,\one_{\{a_0 <|x| \leq a_1\}}\star \mu^X\\
&\leq 2 \sup_{y \in (-n-a_1, n+a_1)}|v(s, y)|\,\one_{\{|x|> a_0\}}\star \mu^X
\end{align*}
that belongs to $\mathcal{A}^+$ by \eqref{int_small_jumps} and Lemma \ref{L:c}.
\endproof
	We now recall the following fact, which constitutes a generalization  of 	Proposition 4.5 (formula (4.4) and considerations below) in \cite{BandiniRusso1}, obtained replacing   $ \one_{\{|x|\leq 1\}}$ with  $\frac{k(x)}{x}$, with  $k\in \mathcal K$. 
\begin{proposition}\label{P:2.10}
Let     $v:\R_+ \times \R \rightarrow \R$  be  a function of class $C^{0,1}$. Then, for every $k \in \mathcal K$,  
\begin{align*}
&|v(s,X_{s-} + x)-v(s,X_{s-})|^2\,\frac{k^2(x)}{x^2}\star \mu^X \in \mathcal{A}^+_{\rm loc}.
\end{align*}
		In particular, the process
$
	(v(s,X_{s-} + x)-v(s,X_{s-}))\frac{k(x)}{x}\,\star (\mu^X- \nu^X)
$
is a square integrable purely discontinuous local martingale.
\end{proposition}

	We continue by giving a basic lemma. 
	
\begin{lemma}\label{L:contmtg_d}
Let $k \in \mathcal K$. 
	The maps
	\begin{align*}
	(i)\quad &v \mapsto (v(s, X_{s-}+x)-v(s, X_{s-})) \frac{x-k(x)}{x}\star \mu^X=:D^{v},\\
	(ii)\quad &v \mapsto (v(s, X_{s-}+x)-v(s, X_{s-})) \frac{k(x)}{x}\star (\mu^X- \nu^X)=:M^{v,d}, 
	\end{align*}
	from   $C^{0,1}$ with values in $\mathbb D^{ucp}$ are well-defined and continuous.
\end{lemma}
\proof
Let $T >0$. 

\noindent (i) For every $v \in C^{0,1}$, the process $D^v$  is defined pathwise.  Indeed, let $a_0 >0$ such that $k(x) = x$ for $x \in [-a_0, a_0]$. Then 
\begin{align}\label{supDv}
	&\sup_{t \in [0,\,T]} \Big|\int_{]0,\,t] \times\R}(v(s, X_{s-}(\omega)+x)-v(s, X_{s-}(\omega))) \frac{x-k(x)}{x} \mu^X(ds\,dx)\Big|\notag\\ 
	&\leq  C\sum_{0 < s\leq T}|\Delta v(s, X_{s
})|\one_{\{|\Delta X_s| >a_0\}}, 
	\end{align}
	for some constant $C= C(k)$. 
	Since  the number of jumps of $X$ larger than $a_0$ is finite, previous quantity is finite a.s. 
	
Let $v_\ell \rightarrow 0$ in $C^{0,1}$ as $\ell\rightarrow \infty$.
The continuity follows since replacing $v$ by $v_\ell$, \eqref{supDv} converges to zero a.s., taking into account that $v_\ell$ converges uniformly on compact sets.  

\medskip

\noindent (ii)  For every $v \in C^{0,1}$, 
$M^{v,d}$   is a square integrable local martingale by Proposition \ref{P:2.10}. Moreover, again by Proposition \ref{P:2.10},  taking $v=Id$, we have $|k(x)|^2 \star \nu^X \in \mathcal A^+_{\textup{loc}}$.
    
Let  
 $
 \tau_n^1 := \inf\{t \geq 0:\,\,|X_{t-} |\geq n\}
 $
  and $\tau_n^2 \uparrow \infty$ be an increasing sequence   of stopping times such that 
$\int_{]0,\tau^2_n] \times \R} |k(x)|^2 \nu^X(ds\,dx) \in \mathcal A^+$ and $M^{v,d}_{\tau^2_n \wedge \cdot} $ is a square integrable martingale. 
Take 
$\tau_n:= \tau_n^1 \wedge \tau_n^2$.
At this point, let us fix $\varepsilon >0$. We have 
\begin{align*}
	\P\Big(\sup_{t \in [0,\,T]}|M^{v,d}_t|> \varepsilon\Big)&= \P\Big(\sup_{t \in [0,\,T]}|M^{v,d}_t|> \varepsilon, \tau_n \leq T\Big) + \P\Big(\sup_{t \in [0,\,T]}|M^{v,d}_t|> \varepsilon, \tau_n > T\Big)\\
	&\leq \P\Big(\tau_n \leq T\Big)+\P\Big(\sup_{t \in [0,\,T]}|M^{v,d}_{\tau_n \wedge t} |> \varepsilon\Big).
\end{align*}
Let us now choose $n_0 \in \N$ in such a way that there exists  $\delta >0$ such that   $\P(\tau_n \leq T) \leq \delta$ for all $n \geq n_0$.
For $n \geq n_0$, applying the Chebyshev inequality, previous inequality gives 
\begin{align}\label{Mvest}
	\P\Big(\sup_{t \in [0,\,T]}|M^{v,d}_t|> \varepsilon\Big)
		&\leq \delta+\P\Big(\sup_{t \in [0,\,T]}|M^{v,d}_{\tau_n \wedge t} |> \varepsilon\Big)\leq \delta + \frac{\E\Big[\sup_{t \in [0,\,T]}|M^{v,d}_{\tau_n \wedge t}|^2\Big]}{\varepsilon^2}.
\end{align}
By Doob 
 inequality we have  
\begin{equation}\label{Doob}
\E\Big[\sup_{t \in [0,\,T]}|M^{v,d}_{\tau_n \wedge t}|^2\Big]\leq 
4 \E[|M^{v,d}_{\tau_n \wedge T}|^2]
=4\sper{\langle M^{v,d},M^{v,d}\rangle_{\tau_n \wedge T}}.
\end{equation}
Denoting by $[-m, m]$ the compact support of $k$, we have 
\begin{align}\label{Mvdb}
\langle M^{v,d},M^{v,d}\rangle_{\tau_n \wedge T} &\leq  \int_{]0,\,\tau_n \wedge T] \times\R} \int_0^1 da |\partial_x v(s, X_{s-}+ a x)|^2 |k(x)|^2 \nu^X(ds\,dx)\notag\\
&\leq \sup_{y \in [-(m+n), m+n], s \in [0,\,T]} |\partial_x v(s, y)|^2   \int_{]0,\tau_n \wedge t] \times \R}|k(x)|^2 \nu^X(ds\,dx). 
\end{align}
Collecting \eqref{Doob} and \eqref{Mvdb}, inequality  \eqref{Mvest} becomes 
\begin{align*}
	\P\Big(\sup_{t \in [0,\,T]}|M^{v,d}_t|> \varepsilon\Big)
	 \leq \delta + \frac{4}{\varepsilon^2}\sup_{y \in [-(m+n), m+n], s \in [0,\,T]} |\partial_x v(s, y)|^2\,\E\Big[   \int_{]0,\tau_n \wedge t] \times \R}|k(x)|^2 \nu^X(ds\,dx)\Big].
\end{align*}
Let us now show the continuity with respect to $v$. 
Let $v_\ell \rightarrow 0$ in $C^{0,1}$ as $\ell\rightarrow \infty$.
Previous estimate shows that 
\begin{align*}
	&\limsup_{\ell \rightarrow \infty} \P\Big(\sup_{t \in [0,\,T]}|M^{v_\ell,d}_t|> \varepsilon\Big)\leq  \delta. 
\end{align*}
By the arbitrariness of $\delta$ this shows that $(M^{v_\ell, d})$ converges to zero u.c.p. 
\endproof

\begin{lemma}\label{L:D4}
	Let $k_1, k_2 \in \mathcal K$. Then $|k_1(x)-k_2(x)| \star \nu^X 
	\in \mathcal{A}^+_{\textup{loc}}
	$. 
\end{lemma}
\proof
Let $a_0$ such that $k_1(x) = k_2(x) = x$ on $[-a_0, a_0]$. Then 
\begin{align*}
	|k_1(x)-k_2(x)| \star \nu^X &= |k_1(x)-k_2(x)|\one_{\{|x| >a_0\}} \star \nu^X\leq (||k_1||_{\infty} + ||k_2||_{\infty})\one_{\{|x| >a_0\}} \star \nu^X
\end{align*}
that belongs to $\mathcal{A}^+_{\textup{loc}}$ by Lemma \ref{L:c}. 
\endproof

\begin{lemma}\label{R:B}
Let $v:\R_+ \times \R \rightarrow \R$ be a locally bounded function.
	Condition \eqref{G1_cond} is equivalent to 
	\begin{align*}
&(v(s,X_{s-} + x)-v(s,X_{s-}))\,\frac{k(x)}{x}\in {\mathcal G}^1_{\textup{loc}}(\mu^X), \qquad \forall k \in \mathcal K. 
\end{align*}
\end{lemma}
\proof
Let $k_1, k_2 \in \mathcal K$. It is enough to show that 
$\left|(v(s,X_{s-} + x)-v(s,X_{s-}))\,\frac{k_1(x)-k_2(x)}{x}\right|\star \mu^X\in \mathcal A_{\textup{loc}}^+$. 
Let   $a_0, a_1>0$ 
such that  $k_1(x) = k_2(x) = x$ on $|x| \leq a_0$ $k_1(x)=k_2(x) =0$ for $|x| > a_1$.  
We have 
\begin{align*}
&\Big|(v(s,X_{s-} + x)-v(s,X_{s-}))\,\frac{k_1(x)-k_2(x)}{x}\Big| \star \mu^X \\
&= \Big|(v(s,X_{s-} + x)-v(s,X_{s-}))\,\frac{k_1(x)-k_2(x)}{x}\Big| \one_{\{a_0 < |x| \leq a_1\}}\star \mu^X\\
&\leq \frac{||k_1||_{\infty} + ||k_2||_{\infty}}{a_0}\one_{\{a_0 < |x| \leq a_1\}}\star \mu^X 
\end{align*}
that belongs to $\mathcal A_{\textup{loc}}^+$ by 
Lemma \ref{L:3.32}.

  \subsection{Recalls on  discontinuous semimartingales and related Jacod's martingale problems}\label{A:D}
  
\setcounter{equation}{0}
 \setcounter{theorem}{0}
 
 We recall that a special semimartingale is a semimartingale $X$ which admits a decomposition $X=M+V$, where  $M$ is a local martingale and $V$ is a  finite variation and predictable process such that $V_0=0$, see Definition 4.21, Chapter I,  in \cite{JacodBook}. Such a decomposition is unique, and is called canonical decomposition of $X$, see respectively Proposition 3.16 and Definition 4.22, Chapter I, in \cite{JacodBook}. 
In the following we  set $\mathcal {\tilde K}:=\{k:\R \rightarrow \R \textup{ bounded: } k(x) =x \textup{ in a neighborhood of } 0\}$. 

Assume now that $X$ is a semimartingale 
with jump measure $\mu^X$.
	   Given $k \in \tilde{\mathcal K}$, the process
 $
 X^k= X- \sum_{s \leq \cdot} [\Delta X_s - k(\Delta X_s)]
 $
 is a special semimartingale with unique decomposition
\begin{equation}\label{dec_Xk_semimart}
	X^k=  X^c + M^{k,d} + B^{k, X}, 
\end{equation}
where  $M^{k,d}$ is  a purely discontinuous local martingale such that $M^{k,d}_0=0$,   $X^c$ is the unique continuous martingale part of $X$ (it coincides with the process $X^c$ introduced in  Proposition \ref{P:uniqdec}), and $B^{k, X}$ is a predictable process of bounded variation vanishing at zero. 

Let   now  $(\check \Omega, \check {\mathcal F}, \check {\mathbb F})$ be the canonical filtered space, and   $\check X$ the canonical process. 
According to Definition 2.6, Chapter II in \cite{JacodBook}, the characteristics of $X$ associated with $k \in \tilde {\mathcal  K}$
is  then the triplet $(B^k,C,\nu)$ on $(\Omega,  {\mathcal F}, {\mathbb F})$ such that 
the following items hold.
\begin{itemize}
	\item[(i)] $B^k$ is  $\check {\mathbb F}$-predictable, with finite variation on finite intervals, and $B_0^k=0$, i.e., $B^{k, X}= B^k \circ X$ is the process in \eqref{dec_Xk_semimart};
	\item[(ii)] $C$ is 
	a continuous process of finite variation with $C_0=0$,  i.e., $C^X:= C \circ  X=  \langle \check X^c, \check X^c\rangle$.
	\item[(iii)] $\nu$ is an $\check {\mathbb F}$-predictable random measure on $\R_+ \times \R$,  i.e.,    $\nu^X: = \nu \circ X$ is the compensator of  $\mu^X$.
	\end{itemize}
\begin{theorem}[Theorem  2.42, Chapter II, in \cite{JacodBook}]\label{T: equiv_mtgpb_semimart}
Let $X$ be an adapted  c\`adl\`ag process. 	Let $B^k$ be an  $\check{\mathbb F}$-predictable process, with finite variation on finite intervals, and $B_0^k=0$, 
 $C$ be
	an $\check{\mathbb F}$-adapted continuous process of finite variation with $C_0=0$, 	and  $\nu$ be an $\check{\mathbb F}$-predictable random measure on $\R_+ \times \R$.
There is equivalence between the two following statements.
\begin{itemize} \item [(i)] X is a real  semimartingale with characteristics $(B^k, C, \nu)$. \item [(ii)] For each bounded function $f$ of class $C^2$, the process 
\begin{align}\label{f:ito}
&f(X_{\cdot}) - f(X_0) - \frac{1}{2} \int_0^{\cdot}  f''(X_s) \,dC^X_s-  \int_0^{\cdot}  f'(X_s) \,d B_s^{k,X}\notag\\
&-  \int_{]0,\cdot]\times \R} (f(X_{s-} + x) -f(X_{s-})-k(x)\,f'(X_{s-}))\,\nu^X(ds\,dx)
\end{align} 
is a local martingale.
\end{itemize} 
\end{theorem}                                       
\begin{remark}\label{R:D2}
Assuming item (i) in Theorem \ref{T: equiv_mtgpb_semimart}, if $f$ is a bounded function of class $C^{1,2}$, formula \eqref{f:ito} can be generalized into  
\begin{align*}
&f(t,X_{t}) - f(0,X_0)- \int_0^{t}\frac{\partial }{\partial s}f(s, X_s)ds - \frac{1}{2} \int_0^{t}  \frac{\partial^2 }{\partial x^2}f(s, X_s) \,dC^X_s-  \int_0^{t} \frac{\partial}{\partial x}f(s, X_s) \,d B_s^{k,X}\notag\\
&-  \int_{]0,t]\times \R} \Big(f(s, X_{s-} + x) -f(s, X_{s-})-k(x)\,\frac{\partial }{\partial x}f(s, X_{s-})\Big)\,\nu^X(ds\,dx), \quad t \in [0,\,T].
\end{align*}
\end{remark}

\small
\paragraph{Acknowledgements.} 
 The work of the first named author was partially  supported by PRIN 2015  \emph{Deterministic and Stochastic Evolution equations}. 
The work of the second named author was partially supported by a public grant as part of the
{\it Investissement d'avenir project, reference ANR-11-LABX-0056-LMH,
  LabEx LMH,}
in a joint call with Gaspard Monge Program for optimization, operations research and their interactions with data sciences.

\addcontentsline{toc}{chapter}{Bibliography} 
\bibliographystyle{plain} 
\bibliography{../../BIBLIO_FILE/BiblioLivreFRPV_TESI} 

\end{document}